\documentclass[twoside, 11pt]{article}
\usepackage[english]{babel}

\usepackage{geometry}
\newlength{\defaultparindent} \AtBeginDocument{\setlength{\defaultparindent}{\parindent}}

\usepackage{fancyhdr}
\usepackage{amsmath}
\usepackage{amsthm}
\usepackage{amssymb}
\usepackage{mathrsfs}
\usepackage{accents}
\usepackage{bbm}
\usepackage{stmaryrd}
\usepackage{array}
\usepackage[all]{xy}
\usepackage{pict2e}
\usepackage{graphicx}
\usepackage{xcolor}
\usepackage{hyperref}
\usepackage[nottoc,notlot,notlof]{tocbibind}

\theoremstyle{definition}
\newtheorem{definition}{Definition}[subsection]
\newtheorem{remark}[definition]{Remark}
\newtheorem{example}[definition]{Example}

\theoremstyle{theorem}
\newtheorem{lemma}[definition]{Lemma}
\newtheorem{proposition}[definition]{Proposition}
\newtheorem{theorem}[definition]{Theorem}
\newtheorem{corollary}[definition]{Corollary}

\numberwithin{equation}{subsection}

\theoremstyle{definition}
\newtheorem*{definition*}{Definition}
\theoremstyle{theorem}
\newtheorem*{theorem*}{Theorem}

\theoremstyle{theorem}
\newtheorem{theoremintro}{Theorem}
\newtheorem{corollaryintro}[theoremintro]{Corollary}

\DeclareMathOperator{\ord}{ord}
\DeclareMathOperator{\res}{res}

\DeclareMathOperator{\Spec}{Spec}

\DeclareMathOperator{\Conf}{Conf}
\DeclareMathOperator{\Sym}{Sym}
\DeclareMathOperator{\Hom}{Hom}
\DeclareMathOperator{\Pic}{Pic}
\DeclareMathOperator{\pt}{part}
\DeclareMathOperator{\pr}{pr}
\DeclareMathOperator{\Frac}{Frac}
\DeclareMathOperator{\Id}{Id}
\DeclareMathOperator{\car}{char}
\DeclareMathOperator{\rk}{rk}
\DeclareMathOperator{\Eff}{Eff}
\DeclareMathOperator{\HD}{HD}

\newcommand{\longhookrightarrow}{\lhook\joinrel\longrightarrow}

\newcommand{\longtwoheadrightarrow}{\relbar\joinrel\twoheadrightarrow}

\newcommand{\C}{\mathbb{C}}
\newcommand{\projsp}{\mathbb{P}}
\newcommand{\calB}{\mathcal{B}}
\newcommand{\calH}{\mathcal{H}}
\newcommand{\calL}{\mathcal{L}}
\newcommand{\calM}{\mathcal{M}}
\newcommand{\calP}{\mathcal{P}}
\newcommand{\calV}{\mathcal{V}}
\newcommand{\calW}{\mathcal{W}}
\newcommand{\scrC}{\mathscr{C}}
\newcommand{\scrS}{\mathscr{S}}
\newcommand{\scrT}{\mathscr{T}}

\newcommand{\scrST}{{\scrS_\scrT}}
\newcommand{\CmS}{{\scrC \smallsetminus |\scrS|}}
\newcommand{\Poinc}{{\mathcal{P}\mathrm{oinc}}}
\newcommand{\trivlb}{{\mathcal{O}_\scrC}}
\newcommand{\hatOCinf}{\widehat{\mathcal{O}_{\scrC_\kappa,\infty}}}

\newcommand{\inddd}{{\mathbbm{1}_{d=2}}}
\newcommand{\inddt}{{\mathbbm{1}_{d\geq3}}}
\newcommand{\indl}{{\mathbbm{1}_{l=0}}}

\title{Equidistribution of curves in hypersurfaces of low degree}
\author{Ariane BUREL}
\date{}

\begin{document}

\maketitle

\begin{abstract}
We prove that over $\C$, for any smooth projective curve $\scrC$ and any  smooth hypersurface $V \subset \projsp^n$ of degree $d$ such that $n \geq 2^d(d-1)$, Faisant's equidistribution principle holds for $(\scrC,V)$. In other words, for any reasonable set of jet conditions $W$, the moduli space $\Hom^e(\scrC,V|W)$ of degree $e$ maps $\scrC \to V$ satisfying conditions $W$ converges in the Grothendieck ring of complex varieties, after renormalisation, as $e \to \infty$. This provides the dimension and decomposition into irreducible components of $\Hom^e(\scrC,V|W)$ for~$e$ large enough. The strategy follows Bilu and Browning's motivic version of the circle method, and Hase-Liu's framework for genus $\geq1$.
\end{abstract}

\makeatletter
\newcommand{\blindnote}{\xdef\@thefnmark{}\@footnotetext} 
\makeatother
\blindnote{2020 {\it Mathematics Subject Classification.} 14H10, 14E18, 14J70, 11P55, 14G40.}

\tableofcontents

\newpage
\section{Introduction}

\subsection{Background}

Let $K$ be a field and $f_1, \dots, f_R \in K[x_0, \dots, x_n]$ non-singular homogeneous polynomials of common degree $d\geq2$ (we denote $\mathbf{f} = (f_1, \dots, f_R)$), which define an affine complete intersection $\tilde V_K \subset \mathbb{A}_K^{n+1}$, smooth outside zero. In the case $K = \mathbb{Q}$ Birch proved the following asymptotic count for the integral points of $\tilde V_K$.

\begin{theorem*}[\cite{Bir62}, section 7, theorem 1, simplified]
As long as $n \geq R(R+1)2^d(d-1)$, there exists $\varepsilon>0$ such that
\[
\#\left\{\left. \mathbf{x} \in \llbracket -B;B \rrbracket^{n+1} \,\right|\, \mathbf{f}(\mathbf{x})=0 \right\} \underset{B\to\infty}{=} B^{n-Rd+1}\left( \prod_v \lambda(v) + O(B^{-\varepsilon})\right) ,
\]
where the product is over all places $v$ of $\mathbb{Q}$ and the $\lambda(v)$ are local densities of solutions. If there are solutions at all places then this product is positive.
\end{theorem*}

He used the \emph{Hardy-Littlewood circle method}, an important tool from analytic number theory, based on the simple formula
\begin{equation} \label{orth_int}
\int_0^1 \exp(2\pi i \alpha N) \,\mathrm{d}\alpha =
\left\{ \begin{array}{ll}
1 & \text{if } N = 0 \\
0 & \text{otherwise}
\end{array} \right.
\end{equation}
for $N \in \mathbb{Z}$: this allows one to count the zeros of an integer expression $N$, and in particular to rewrite the cardinality above as the integral of an exponential sum. Note that more recently Rydin Myerson improved on the method and brought $n$ down to be linear in the number $R$ of equations (theorem 1.3 of \cite{RM18}, see also \cite{RM19}).

A more geometric version of the above consists in taking $K = k(t)$ with some base field $k$: the integral points we are counting become
\begin{equation} \label{intpt_affine_hyp}
\left\{\left. \mathbf{g} \in k[t]_{\leq e}^{n+1} \,\right|\, \mathbf{f}(\mathbf{g})=0 \right\}
\end{equation}
for some $e\geq0$, where $k[t]_{\leq e}$ is the space of polynomials of degree $\leq e$. Lee (\cite{Lee11}) established a direct analog of Birch's theorem for $K = \mathbb{F}_q(t)$ with $\car(\mathbb{F}_q) > d$ (this variant of the circle method is called `functional'). In the case where $\mathbf{f}$ has coefficients in $k$, one can consider the smooth projective complete intersection $V \subset \projsp_k^n$ defined by $\mathbf{f}$, and then (\ref{intpt_affine_hyp}) is closely related to the moduli space of degree $e$ morphisms\footnote{\label{f_coef_variable}More generally when $\mathbf{f}$ is allowed to have coefficients in $k(t)$, one finds instead the moduli space $\Hom_{\projsp^1}(\projsp^1, \mathscr{V})$ of sections of $\mathscr{V} \to \projsp^1$, where $\mathscr{V}$ is a model of $V$ over $\projsp^1$. For generalities on this space of sections, see \cite{Fai25_1}, paragraph 1.1.2.} from the projective line to $V$, or more precisely to its $k$-points
\begin{equation} \label{kpt_Hom}
\Hom^e(\projsp^1, V)(k) =
\left\{ \mathbf{g} \in k[t]_{\leq e}^{n+1}\smallsetminus\{0\} \,\left|\,
\begin{array}{l}
\max\deg g_i = e \\
\gcd(g_0, \dots, g_n) = 1 \\
\mathbf{f}(\mathbf{g}) = 0
\end{array} \right.\right\} / k^\times .
\end{equation}
These moduli spaces are of great interest and parametrize, up to automorphisms of $\projsp^1$, the rational curves in $V$. Let us now focus on the case $R=1$ of a hypersurface $V$, defined by $f = f_1$. Following the functional circle method and under the assumption\footnote{This is the assumption of \cite{BS23} after improving on \cite{BV17}.} $n \geq 2^{d-1}(2d-1)$, Browning, Vishe and Sawin (\cite{BV17, BS23}) established an asymptotic for the count of (\ref{kpt_Hom}) over $k = \mathbb{F}_q$, deducing that $\Hom^e(\projsp^1, V)$ was irreducible of the expected dimension over finite fields of characteristic $>d$ ; from there they could obtain the homologous result over $k=\C$ as stated below, thanks to a \emph{spreading out} argument.

\begin{theorem*}[\cite{BS23}, theorem 1.1]
Let $V \subset \projsp_\C^n$ be a smooth hypersurface of degree $d\geq3$ with $n \geq 2^d(d-1/2)$: for all $e\geq1$, $\Hom^e(\projsp^1, V)$ is irreducible of dimension $(n-d+1)e + n-1$.
\end{theorem*}

Let us make the notion of \emph{expected dimension} clear. If $\scrC$ is a smooth curve of genus $g$ and $X$ a Fano variety, let $\Hom^{\mathbf{e}}(\scrC, X)$ be the moduli space of morphisms $\scrC \to X$ of multidegree $\mathbf{e}$. Under weak assumptions, deformation theory (or more prosaically, analysing how many equations define the variety locally) establishes the lower bound
\begin{equation} \label{expect_dim}
\dim \Hom^\mathbf{e}(\scrC, X) \geq \omega_X^{-1} \cdot \mathbf{e} - (g-1) \dim X
\end{equation}
where $\omega_X^{-1}$ is the anticanonical bundle of $X$. This bound is often the actual dimension of the space, and is generally called its expected dimension. 

Back to the theorem above: with the same strategy, \cite{BVY24} extended the result to the case $R\geq2$ of complete intersections, and to morphisms with a finite number of prescribed values. Hase-Liu carried out analogous computations this time in positive genus (and with $R=1$), as follows.

\begin{theorem*}[\cite{HL24}, theorem 1, simplified]
Let $V \subset \projsp_\C^n$ be a smooth hypersurface of degree $d\geq2$ with $n \geq 2^d(d-1)+1$, and $\scrC$ be a smooth projective complex curve of genus $g\geq1$. For $e$ large enough, the moduli space $\Hom^e(\scrC,V)$ is irreducible of the expected dimension, that is
\[
(n-d+1)e - (n-1)(g-1) .
\]
\end{theorem*}

In order to achieve this he replaced the polynomials with global sections of line bundles over $\scrC$, and adapted the method to this more geometric setting; in particular this involves replacing arguments from geometry of numbers with vector bundles and slope theory. Let us also mention that again using the functional circle method, \cite{BS23} bounded the dimension of the singular locus of $\Hom^e(\projsp^1,V)$, and \cite{GHL24} showed that under some assumptions, the singularities of $\Hom^e(\scrC,V)$ are at worse terminal.

Taking a step back, the results of Birch and Lee mentioned first are instances of the Batyrev-Manin-Peyre principle (or \emph{BMP}), which conjectures that the number of rational points of height $\leq B$ on a nice enough Fano variety (after removing a potential accumulation subvariety) is equivalent to
\begin{equation}
C B^a (\log B)^{b-1}
\end{equation}
as $B\to\infty$, for certain constants $a>0$, $b\geq1$ and $C>0$ (\cite{FMT89, BM90, Pey95, BT98}). Peyre (\cite{Pey21}) and Faisant (\cite{Fai25_1}, paragraph 3.1.3) later introduced a geometric version of the BMP principle, concerned with the behaviour of the moduli space $\Hom^{\mathbf{e}}(\scrC, X)$ asymptotically in $\mathbf{e}$. Naturally this moduli space is not finite over a general field, and one needs to make sense of the notion of asymptotic behaviour in this context.

To allow for this let us quickly introduce, for any noetherian scheme $S$, the \emph{Grothendieck ring of varieties over $S$}, denoted $K_0\left(\mathrm{Var}_S\right)$ (see \cite{CLL16}). A generalisation of it will be presented with more details in paragraph \ref{parag_GR}. This object is first defined as an abelian group using generators and relations: the generators are the isomorphism classes of $S$-varieties, while the relations are \emph{cut-and-paste}: $X-Y-(X \smallsetminus Y)$ for any $S$-variety $X$ and closed subscheme $Y$ of $X$. The class of $X$ is written $[X]$. Then the product is simply the fibered one: $[X] \!\!\cdot [Y] := \left[ X \times_S Y \right]$. We single out the class $\mathbf{L} := \left[ \mathbb{A}_S^1 \right]$, that will allow us to renormalise expressions: by inverting it we obtain the ring $\mathscr{M}_S := K_0\left(\mathrm{Var}_S\right)\left[\mathbf{L}^{-1}\right]$. One then chooses a topology on this ring, the most natural example being the dimensional one. The associated completion is denoted $\widehat{\mathscr{M}}_S$. This ring allows for computations using the varieties directly, and this framework may be called \emph{motivic}. Furthermore Bilu introduced in \cite{Bil23} the concept of \emph{motivic Euler products} indexed by a variety, analogous to the usual Euler products (see paragraph \ref{MEP}).

Let $\scrC$ be a smooth projective curve and $X$ be a Fano variety, both over a field $k$. The geometric BMP principle states that the class $\left[\Hom^{\mathbf{e}}(\scrC, X)\right]$, after renormalisation by its expected dimension, should stabilise in $\widehat{\mathscr{M}}_k$ as a motivic Euler product of motivic volumes. More generally, Faisant formalized a principle of equidistribution for curves, in which only morphisms that satisfy some fixed \emph{jet conditions} are considered; this is also the `geometrisation' of a homonymous principle for rational points. Faisant's principle of equidistribution is as follows. Fix a  closed subscheme $\scrS = \sum_{p\in|\scrS|} l_p p$ of $\mathscr{C}$ of dimension $\leq0$, and for all $p\in|\scrS|$ take a constructible subset $W_p \subset \Hom(l_p p, X)$, where $l_p p$ is to be understood as a point of multiplicity $l_p$. The space of allowed restrictions is taken as $W = \prod_{p\in|\scrS|} W_p \subset \Hom(\scrS, V)$. To avoid accumulation phenomena on certain closed subschemes, choose a dense open~$U$ of $X$ appropriately. One denotes by
\[
\Hom^{\mathbf{e}}(\scrC, X | W)_U
\]
the moduli space of morphisms $\scrC \to X$ of multidegree $\mathbf{e}$, whose restriction to $\scrS$ lies in $W$, and that intersect $U$. Finally write $\Eff(X)$ for the effective cone of $X$ and $J(\scrC)$ for the Jacobian of $\scrC$.
\begin{definition*}[Principle of equidistribution -- \cite{Fai25_1}, definition 4.3, simplified\footnote{\label{equidistrib_model}
Faisant states the principle of equidistribution more generally for sections of a proper model $\mathscr{X} \to \scrC$. When $X$ is proper, his `model change' theorem 4.6 allows one to reduce to the case of the trivial model $X \times \scrC \to \scrC$, that is to morphisms $\scrC \to X$ as written here.}]
One says that $(X, \scrC)$ satisfies the principle of equidistribution for a given topology on $\mathscr{M}_k$ if there exists a $U$ such that for all $(\scrS, W)$ as above,
\begin{equation} \label{Hom_renorm}
\mathbf{L}^{-(\omega_X^{-1} \cdot \mathbf{e} - (g-1) \dim X)} \left[ \Hom^{\mathbf{e}}(\scrC, X | W)_U \right]
\end{equation}
converges in $\widehat{\mathscr{M}}_k$ to
\begin{align} \label{Fai_Tamagawa}
\frac{1}{\left(1-\mathbf{L}^{-1}\right)^r} \cdot \frac{\left[J(\scrC)\right]}{\mathbf{L}^g} \left(\prod_{x\in\scrC\smallsetminus|\scrS|} \left(1-\mathbf{L}^{-1}\right)^r \frac{[X]}{\mathbf{L}^{\dim X}}\right) \left(\prod_{p\in|\scrS|} \left(1-\mathbf{L}^{-1}\right)^r \frac{[W_p]}{\mathbf{L}^{l_p \dim X}}\right)
\end{align}
with $r = \rk\Pic(X)$, as $\mathbf{e} \in \Pic(X)^\vee \cap \Eff(X)^\vee$ goes arbitrarily far away from the boundary of $\Eff(X)^\vee$.
\end{definition*}
The first product is a motivic Euler product as mentioned above, while the second one is a regular finite product. Observe that this expression remains the same if we exchange the condition data $(\scrS, W)$ for another of the same `shape', hence the terminology of equidistribution. Naturally, properties of $\Hom^\delta(\scrC, V | W)_U$ can be deduced, including data about its dimension and irreducibility. The geometric BMP principle is the special case with no conditions: $\scrS = \varnothing$. For these two principles, results are known for certain varieties $X$. Faisant showed that equidistribution, with any $\scrC$, holds for equivariant compactifications of vector spaces over characteristic $0$ algebraically closed fields (\cite{Fai23}, using \cite{Bil23} and \cite{CLL16}); and for smooth projective split toric varieties over any field (\cite{Fai25_1, Fai25_2}, with the help of previous work from \cite{Bou09}). But more important for us is the following result of Bilu and Browning.

\begin{theorem*}[\cite{BB25}, theorem 1.3]
Smooth hypersurfaces $V \subset \projsp_\C^n$ of degree $d\geq3$ with $n \geq 2^d (d-1)$, together with $\scrC = \projsp^1$, satisfy the BMP principle (with $U=V$ and the weight topology).
\end{theorem*}

For this they developed a motivic version of the functional circle method, which means that the whole process is carried out over $\C$, and spreading out is completely bypassed. Here the dimensional topology does not suffice for the methods used, and a weight function on $\mathscr{M}_\C$, built on Hodge theory, takes it place. The convergence of the normalised class of $\Hom^e(\projsp^1,V)$ allows one to retrieve the dimension and irreducibility of $\Hom^e(\projsp^1, V)$ like in \cite{BS23}, but with a slightly looser condition between $n$ and $d$. On top of this, approximating the Grothendieck class of a variety provides more data than just these, including for instance a proportion of the coefficients of its Hodge-Deligne polynomial.

The main concern of the present article is to show that the same hypersurfaces satisfy the principle of equidistribution for any smooth projective curve, again with no accumulation phenomenon; see \textbf{theorem \ref{th_main_intro}} (we also include the case $d=2$). \textbf{Corollary \ref{corol_main_intro}} establishes the dimension of $\Hom^e(\scrC, V | W)$ and its decomposition into irreducible components when $e$ is large enough. This naturally generalises the result of \cite{BB25} discussed just above, as well as \cite{HL24}'s theorem 1 on $\Hom^e(\scrC,V)$ that we mentioned earlier. Furthermore Hase-Liu used spreading out, while our computations are motivic.

Before moving on towards the announced results, let us make some more remarks about previous works. First, let us go back to the article \cite{BVY24} that was mentioned towards the start of the introduction, and more precisely to its theorem 1.2, that treats morphisms with prescribed values. With the notations from the principle of equidistribution this corresponds to conditions $W_p = \{y_p\}$ of multiplicity $l_p = 1$. Our result therefore also generalises the case $R=1$ of this one.

A link can also be made with weak approximation. Observation 1.2 of \cite{Has10} states that for a smooth proper variety $X$ over $K(\scrC)$ that admits a regular model $\mathscr{X} \to \scrC$, weak approximation is satisfied as soon as for all $\scrS \subset \scrC$ closed of dimension $\leq0$ and all $\phi \in \Hom_\scrC(\scrS, \mathscr{X})$, there exists $e$ such that $\Hom_\scrC(\scrC, \mathscr{X} | \{\phi\}) \neq \varnothing$. In particular the principle of equidistribution is a vast generalisation of weak approximation (see footnote \ref{equidistrib_model}). For smooth projective hypersurfaces, weak approximation has already been proven by Hassett and Tschinkel for $n$ larger than an explicit fast-growing function of $d$ (theorem 1 of \cite{HT09}), then by de Jong and Starr as soon as $n \geq d(2d-1)-1$ (corollary 1.5 of \cite{dJS06}); our theorem \ref{th_main_intro} retrieves it for $n$ exponential in $d$.

Independently of the above, let us mention that in a dramatic and very recent improvement on the functional version of Birch's theorem, Hase-Liu exchanged the exponential relation between $n$ and $d$ for a quadratic one in \cite{HL26}. Drawing on ideas of Pugin and Sawin (\cite{Pug11, Saw26}), the method to control the minor arcs is completely changed: the exponential sums are bounded using the dimension of a related singular locus (Hase-Liu's theorem 7).

Finally, we have not yet mentioned purely geometric methods for the study of the moduli spaces of morphisms (apart from the expected dimension obtained by deformation theory). These, such as bend-and-break, usually only produce results for $X$ \emph{general}, although they may otherwise require weaker hypotheses than arithmetic methods. In the case of hypersurfaces $V$, improving on \cite{HRS04} and \cite{BK13}, Riedl and Yang proved the following.

\begin{theorem*}[\cite{RY19}, theorem 3.3, simplified]
As soon as $n \geq d+2$, $\Hom^e(\projsp^1, V)$ is of the expected dimension for general hypersurfaces $V \subset \projsp_\C^n$ of degree $d\geq1$.
\end{theorem*}

The condition on $n$ is much more optimal than what can be achieved with the circle method (be it the classical exponential assumption or Hase-Liu's new quadratic one), and the two approaches are complementary. It is also important to note that these geometric methods are, for now at least, incompatible with the adjunction of jet conditions such that the $W$ from the principle of equidistribution; which a priori leaves arithmetic methods as the sole strategy in this case. To our knowledge, no result in positive genus has been obtained with geometric methods either.

\subsection{Our framework}

We go over important notations introduced in the previous subsection, and introduce some additional ones. Let $V \subset \projsp_\C^n$ be a smooth projective hypersurface of degree $d\geq1$: it is given by the vanishing of a homogeneous polynomial $f \in \C[x_0, \dots, x_n]$ which is of degree $d$ and nonsingular outside of zero. Our goal is to study the spaces of curves in $V$ satisfying certain jet conditions: we assimilate such curves to morphisms from a given smooth projective curve $\scrC$ of genus $g$ over $\C$ to our hypersurface $V$, and consider the moduli spaces 
\[
\Hom^e(\scrC,V)
\]
of degree $e$ morphisms $\scrC \to V$, with $e\geq0$. Additionally let us fix $\scrS$ a subscheme of $\scrC$ of dimension at most $0$ (that is, either $\scrS=\varnothing$ or $\dim\scrS=0$), which will be the condition locus. We can also view it as an effective divisor on $\scrC$, and write it
\begin{equation} \label{def_scrS}
\scrS = \sum_{p\in|\scrS|} (m_p+1) p
\end{equation}
where $m_p \geq0$ is the depth of the point $p$ in $\scrS$. Here we wrote the multiplicity as $m_p+1$ to respect usual notations, but for practicality we will mostly use $l_p := m_p+1$ instead. Schematically a point with multiplicity is
\begin{equation}
l_p p = \Spec\left( \C[t]/t^{l_p} \right) .
\end{equation}
We also denote
\[
l := \sum_{p\in|\scrS|} l_p = \deg \scrS .
\]
The jet conditions we mentioned consist in imposing that morphisms $\scrC \to V$ restrict to $\scrS$ in a certain way. In order to formalise this clearly, we need more notation.

We write $\Pic^e(\scrC)$ for the subvariety of $\Pic(\scrC)$ made up of the degree $e$ classes, and $J(\scrC)$ for the Jacobian of the curve; one has
\[
\Pic^e(\scrC) \simeq \Pic^0(\scrC) = J(\scrC) ,
\]
which will allow us to consider objects equivalently over $\Pic^e(\scrC)$ or $J(\scrC)$. We denote, for $L$ a line bundle over $\scrC$, $L^{\otimes r}$ its $r$-th tensor power for $r\in\mathbb{Z}$, and
\[
P_L := H^0(\scrC,L)
\]
its space of global sections\footnote{We take inspiration from \cite{HL24} for denoting our spaces of sections, although in a slightly different way: our $P_L$ would be his $P_{e,\scrC}$ and our $P_{|Z}$ his $P_{e,Z}$.}. Via the cup product $P_{L_1} \times P_{L_2} \to P_{L_1\otimes L_2}$, we can see our polynomial $f$ as a map
\[
P_L^{n+1} \longrightarrow P_{L^{\otimes d}} ,
\]
which means that for any $\mathbf{u} = (u_0, \dots, u_n) \in P_L^{n+1}$ we will consider $f(\mathbf{u}) \in P_{L^{\otimes d}}$. If $Z$ is an effective divisor on $\scrC$, we define the space of \emph{local} sections 
\[
P_{|Z} := H^0(Z, \trivlb|_Z) .
\]
This time the cup product endows $P_{|Z}$ with a $\C$-algebra structure, and $f$ induces again a map 
\[
P_{|Z}^{n+1} \longrightarrow P_{|Z} .
\]
Because $P_{|Z}$ is a ring, its multiplicative group $\mathbb{G}_{m, P_{|Z}}$ can be considered. The following fact will be justified at the start of paragraph \ref{global_sect}.
\begin{remark} \label{sect_loc}
Actually for any $L$ whose base locus does not intersect $Z$, $H^0(Z, L|_Z)$ is isomorphic to $P_{|Z}$. In this case there is a restriction
\[
P_L \overset{\cdot|_Z}{\longrightarrow} H^0(Z, L|_Z) \simeq P_{|Z}
\]
that commutes with the cup product, and thus with $f$.
\end{remark}

\begin{remark}
Since $P_L$ and $P_{|Z}$ are finite-dimensional vector spaces, we can see them as affine varieties; in fact we will adopt this point of view by default, and when we want to see $P_{|Z}$ as a $\C$-algebra  (for instance) we will say it clearly.
\end{remark}

Recall that any line bundle of degree at least $2g$ is basepoint-free (\cite{Har77}, chapter IV, corollary 3.2). As a consequence, for
\begin{equation} \label{cond_BPF}
e\geq2g
\end{equation}
there is a map
\begin{equation} \label{map_Hom_Pic}
\Hom^e(\scrC,V) \longrightarrow \Pic^e(\scrC) \simeq J(\scrC)
\end{equation}
whose fibre over a line bundle $L$ is
\begin{align} \label{def_Hom_L}
\Hom^L(\scrC,V) & = \left\{ \mathbf{u} \in P_L^{n+1} \,\left|\,
\begin{array}{l}
u_0, \dots, u_n \text{ generate } P_L \\
f(\mathbf{u}) = 0
\end{array} \right.\right\} /\mathbb{G}_m \nonumber \\
& = \left\{ \mathbf{u} \in P_L^{n+1} \,\left|\,
\begin{array}{l}
\forall x\in\scrC, \mathbf{u}(x)\neq0 \\
f(\mathbf{u}) = 0
\end{array} \right.\right\} /\mathbb{G}_m
\end{align}
(see remark \ref{constr_fibr_Hom} for a justification).
As for the space of morphisms $\scrS \to V$, it can be decomposed as
\[
\Hom(\scrS,V) = \prod_{p\in|\scrS|} \Hom(l_p p, V)
\]
with factors
\begin{equation} \label{Hom_lp}
\Hom(l_p p, V) = \left\{ \mathbf{w} \in P_{|l_p p}^{n+1} \,\left|\,
\begin{array}{l}
\mathbf{w}(p)\neq0 \\
f(\mathbf{w}) = 0
\end{array} \right.\right\} \left/ \mathbb{G}_{m, P_{|l_p p}} \right. ,
\end{equation}
where $\mathbb{G}_{m, P_{|l_p p}}$ acts by diagonal multiplication (for more insight on this description, see subsection \ref{jet_sch}). Given that $P_{|\scrS} = \prod_{p\in|\scrS|} P_{|l_p p}$ (both as affine varieties and as $\C$-algebras), we have
\begin{equation} \label{Hom_S}
\Hom(\scrS, V) = \left\{ \mathbf{w} \in P_{|\scrS}^{n+1} \,\left|\,
\begin{array}{l}
\forall p\in|\scrS|, \mathbf{w}(p)\neq0 \\
f(\mathbf{w}) = 0
\end{array} \right.\right\} \left/ \mathbb{G}_{m, P_{|\scrS}} \right. ,
\end{equation}
with again a diagonal action. We can now express the restriction
\[
\cdot|_\scrS : \Hom^L(\scrC,V) \longrightarrow \Hom(\scrS,V)
\]
using the composition
\[
P_L^{n+1}/\mathbb{G}_m \overset{\cdot|_\scrS}{\longrightarrow} P_{|\scrS}^{n+1}/\mathbb{G}_m \longtwoheadrightarrow P_{|\scrS}^{n+1} \big/ \mathbb{G}_{m, P_{|\scrS}} ,
\]
where the second arrow allows for multiplication by any invertible local section rather than just a constant. Now for every $p\in|\scrS|$ we fix a non empty \emph{constructible} subset
\[
W_p \subset \Hom(l_p p, V),
\]
that will be our authorised restrictions to $l_p p$, and set
\[
W := \prod_{p\in|\scrS|} W_p \subset \Hom(\scrS,V).
\]
We say a morphism $\scrC \to V$ satisfies the jet conditions $W$ if its restriction to $\scrS$ lies in $W$, and define
\begin{equation} \label{def_Hom_LW}
\Hom^L(\scrC,V|W) := \left\{\left. \varphi \in \Hom^L(\scrC,V) \,\right|\, \varphi|_\scrS \in W \right\}
\end{equation}
as the moduli space of such morphisms. As in (\ref{map_Hom_Pic}), for $e\geq2g$ there exists
\[
\Hom^e(\scrC,V|W) \subset \Hom^e(\scrC,V)
\]
mapping to $\Pic^e(\scrC) \simeq J(\scrC)$ with fibres $\Hom^L(\scrC,V|W)$ (again, see remark \ref{constr_fibr_Hom}). As mentioned earlier, this last moduli space is the subject of our study, and we are particularly interested in its \emph{asymptotic behaviour} as $e\to\infty$; or more precisely in that of its classes in the rings $K_0\left(\mathrm{Var}_{J(\scrC) \times W}\right)$ and $\mathscr{M}_{J(\scrC) \times W}$ that were introduced in the previous subsection. To complete what was said there, let us add that if $S \to T$ is a map of noetherian schemes, any class over $T$ can be viewed as over $S$ via base change (see (\ref{GR_functor})). Furthermore just like in \cite{BB25} the dimensional topology would not work for us so we make use of the weight function $w_S : \mathscr{M}_S \to \mathbb{Z} \cup \{-\infty\}$ (such that only zero has weight $-\infty$) that they construct for complex varieties $S$. Recall that the associated completion is denoted $\widehat{\mathscr{M}}_S$.

In order to do the computations we will first focus on \emph{naive}, affine versions of the spaces of morphisms. For this purpose denote by $\tilde W$ the inverse image of W via the projection $P_{|\scrS}^{n+1} \twoheadrightarrow P_{|\scrS}^{n+1} \big/ \mathbb{G}_{m, P_{|\scrS}}$; for all $p\in|\scrS|$ we will also call $\tilde W_p$ the inverse image of $W_p$ via $P_{|l_p p}^{n+1} \twoheadrightarrow P_{|l_p p}^{n+1} \big/ \mathbb{G}_{m, P_{|l_p p}}$, so that $\tilde W = \prod_{p\in|\scrS|} \tilde W_p$. For $\deg L \geq 2g$ we consider the naive version of the moduli space $\Hom^L(\scrC,V|W)$ from (\ref{def_Hom_L}),
\begin{equation} \label{def_ML}
M_L^W := \left\{ \mathbf{u} \in P_L^{n+1} \,\left|\,
\begin{array}{l}
f(\mathbf{u}) = 0 \\
\mathbf{u}|_\scrS \in \tilde W
\end{array} \right.\right\} .
\end{equation}
In subsection \ref{CM_framework} we will also construct a variety
\begin{equation} \label{fibr_Me}
\calM_e^W \longrightarrow \Pic^e(\scrC) \simeq J(\scrC)
\end{equation}
with fibres $M_L^W$. In remark \ref{expect_dim_naive} we explain that we expect $\calM_e^W$ to have $J(\scrC) \times \tilde W$-dimension
\begin{equation} \label{def_xi}
\xi_e := (n-d+1)e - n(g-1) - nl .
\end{equation}
As for $\Hom^e(\scrC,V|W)$, one expects the same $J(\scrC) \times W$-dimension as for $\calM_e^W / \mathbb{G}_m$, that is $\xi_e - 1 + l$ ($+l$ for the passage from $\tilde W$ to $W$, and $-1$ for the quotient by $\mathbb{G}_m$). In particular one retrieves the expected complex dimension
\[
(n-d+1)e - (n-1)(g-1) - ((n-1)l - \dim W)
\]
from (\ref{expect_dim}), in a version with conditions $W$. We will see in remark \ref{dim_jet} that $\Hom(l_p p,V)$ has dimension $(n-1)l_p$, and thus
\[
\dim \Hom(\scrS,V) = (n-1)l ;
\]
as a consequence, in the expected dimension above, $(n-1)l - \dim W$ is the codimension of $W$ in $\Hom(\scrS,V)$. 

\begin{remark}
Now is a good moment to stress that most important notations are consistent among the main references we use, with the noticeable exception of the dimension of the projective hypersurface, which is usually either $n-2$ or $n-1$. For us this dimension is $n-1$; as a consequence for instance, the $n$ from \cite{BB25} is our $n+1$, while the $n$ from \cite{HL24} is the same as ours.
\end{remark}

\subsection{Results}

\begin{theoremintro}[theorem \ref{th_naive}, simplified] \label{th_naive_intro}
Assume $d\geq2$, $n \geq 2^d(d-1)$ and
\[
e \geq \max(\theta_1^0, \theta_2^0, 1+\inddd) \quad\text{with}\quad
\left\{\begin{array}{l}
\theta_1^0 = 4g-\inddd+2l\\
\theta_2^0 = 2dg-2+l
\end{array}\right. .
\]
Then there exist $a_0<0, a_1, a_2, a_3$ that depend on $d$ and $n$, such that in $\widehat{\mathscr{M}}_{J(\scrC) \times \tilde W}$,
\[
\mathbf{L}^{-\xi_e} \left[\calM_e^W\right] = \mathbf{L}^{-l} \mathfrak{S} + \varepsilon_e
\]
where:
\begin{itemize}
\item $\mathbf{L}^{-l} \mathfrak{S}$ can, over $\C$, be written as the motivic Euler product of \emph{local densities}
\begin{equation} \label{intro_mep_local_dens}
\mathbf{L}^{-l} \mathfrak{S} = \prod_{x\in\CmS} \lim_{N\to\infty} \mathbf{L}^{-nN} \left[ \Lambda_N(x) \right] \in \widehat{\mathscr{M}}_\C
\end{equation}
with $\Lambda_N(x)$ jet schemes defined in (\ref{def_Lambda}), such that $w_\C \left(\mathbf{L}^{-l} \mathfrak{S} - 1\right) < 0$;

\item and $\varepsilon_e$ is an error term satisfying
\[
w_{J(\scrC) \times \tilde W}(\varepsilon_e) \leq a_0 e + a_1 g +a_2 l + a_3 .
\]
\end{itemize}
In particular,
\[
\mathbf{L}^{-\xi_e} \left[\calM_e^W\right] \underset{e\to\infty}{\longrightarrow} \mathbf{L}^{-l} \mathfrak{S} \in \widehat{\mathscr{M}}_{J(\scrC) \times \tilde W} .
\]
\end{theoremintro}

The complete version of this theorem actually holds over
\[
J(\scrC) \times \Sym^\mu(\CmS) \times \tilde W
\]
for any partition $\mu$, which allows us to impose additional vanishing conditions parametrized by $\Sym^\mu(\CmS)$. This is required to deduce theorem \ref{th_main_intro} from this one via motivic inclusion-exclusion, as will be made clear in subsection \ref{subsec_MIE}. In particular the $\xi_e$ we have here is the $\xi_{e,0}$ from (\ref{def_xi_mu}), and $\mathfrak{S}$ is the $\mathfrak{S}^\varnothing$ from lemma \ref{sing_cv} or corollary \ref{sing_locdens}. The coefficients $a_j$ represent the $\rho$ from (\ref{def_rho_max})
\footnote{For instance $a_0$ can be taken to be the maximum of the coefficients of the $\rho_i$ in $e$, etc; the resulting bound is however suboptimal since linear combinations do not commute with maxima. An analogous remark can be made for theorem \ref{th_main_intro} and corollaries \ref{corol_naive_intro} and \ref{corol_main_intro}.}.

One remark for the reader familiar with the circle method: $\mathfrak{S}$ is the singular series, while no singular integral appears for us. The reason is that we do not single out some point $\infty$ of $\scrC$, unlike what usually happens when $\scrC = \projsp^1$. In particular the expression of $\mathbf{L}^{-l}\mathfrak{S}$ as a motivic Euler product is indexed by all of $\CmS$, as opposed to lacking a factor at $\infty$.

\begin{corollaryintro}[corollary \ref{corol_naive}, simplified] \label{corol_naive_intro}
Assume again $d\geq2$ and $n \geq 2^d(d-1)$. There exist $b_1, b_2, b_3$ that depend on $d$ and $n$ such that for all
\[
e > b_1 g + b_2 l + b_3 ,
\]
\[
\dim_{J(\scrC) \times \tilde W} \calM_e^W = \xi_e = (n-d+1)e - n(g-1) - nl .
\]
In particular over $\C$ this becomes
\[
\dim \calM_e^W = (n-d+1)e - (n-1)(g-1) + 1 - \left(nl - \dim \tilde W\right) .
\]
Furthermore under the same condition on $e$, $\calM_e^W$ has the same number of irreducible components of maximal dimension as $W$ and $\tilde W$.
\end{corollaryintro}

Theorem \ref{th_naive_intro} is a step towards the following result on the moduli spaces of morphisms, which is a simplified reunion of propositions \ref{th_main_no_estimate} and \ref{error_mininfty}.

\begin{theoremintro}[Equidistribution of curves in hypersurfaces of low degree] \label{th_main_intro}
Assume either $d\geq3$ and $n \geq 2^d(d-1)$, or $d=2$ and $n \geq 2^{d+1}(d-1) = 8$, and suppose also that
\[
e \geq \max(\theta_1^0, \theta_2^0, 1+\inddd, \tau) \quad\text{with}\quad
\left\{\begin{array}{l}
\theta_1^0 = 4g-\inddd+2l\\
\theta_2^0 = 2dg-2+l \\
\tau = 4(d-1)g - 4 + \inddd
\end{array}\right. .
\]
Then there exist $a'_0<0, a'_1, a'_2, a'_3$ that depend on $d$ and $n$, such that in $\widehat{\mathscr{M}}_{J(\scrC) \times W}$ we have
\begin{equation} \label{equidistrib_base}
\mathbf{L}^{-(\xi_e-1+l)} \left[ \Hom^e(\scrC,V|W) \right] = \left(1-\mathbf{L}^{-1}\right)^{\#|\scrS|-1} \cdot \prod_{x\in\CmS} \left( \left(1-\mathbf{L}^{-1}\right) \frac{[V]}{\mathbf{L}^{n-1}} \right) + \mathcal{E}_e
\end{equation}
where $\mathcal{E}_e$ is an error term whose weight satisfies
\[
w_{J(\scrC) \times W}\left(\mathcal{E}_e\right) \leq a'_0 e + a'_1 g +a'_2 l + a'_3 .
\]
In particular the normalised space of morphisms converges to the first term on the right-hand side of (\ref{equidistrib_base}) as $e\to\infty$. Over $\C$ the identity becomes
\begin{align} \label{equidistrib_C}
& \mathbf{L}^{-((n-d+1)e - (n-1)(g-1))} \left[ \Hom^e(\scrC,V|W) \right] \nonumber \\
& = \left(1-\mathbf{L}^{-1}\right)^{-1} \frac{[J(\scrC)]}{\mathbf{L}^g} \cdot \prod_{x\in\CmS} \left( \left(1-\mathbf{L}^{-1}\right) \frac{[V]}{\mathbf{L}^{n-1}} \right) \cdot \prod_{p\in|\scrS|} \left( \left(1-\mathbf{L}^{-1}\right) \frac{[W_p]}{\mathbf{L}^{(n-1)l_p}} \right) \nonumber \\
& \quad + \mathbf{L}^{-g-(n-1)l} \mathcal{E}_e.
\end{align}
As a consequence, the equidistribution principle (\ref{Fai_Tamagawa}) holds for the couple $(V,\scrC)$, with respect to the weight topology, and without accumulation phenomenon ($U=V$).
\end{theoremintro}

Let us justify the passage from (\ref{equidistrib_base}) to (\ref{equidistrib_C}) right away. As we descend over $\C$, a factor $\left[J(\scrC) \times W\right]$ appears, that we want to renormalise. We multiply both sides by $\mathbf{L}^{-g-(n-1)l}$, hence the normalisations of the space of morphisms and of the error term in (\ref{equidistrib_C}). For the main term we now have factors $\mathbf{L}^{-g} [J(\scrC)]$ (which is normalised) and $\mathbf{L}^{-(n-1)l} [W]$. The latter can be combined with $\left(1-\mathbf{L}^{-1}\right)^{\#|\scrS|-1}$ to give
\[
\left(1-\mathbf{L}^{-1}\right)^{-1} \prod_{p\in|\scrS|} \left( \left(1-\mathbf{L}^{-1}\right) \frac{[W_p]}{\mathbf{L}^{(n-1)l_p}} \right) .
\]

As for the naive spaces, we have a corollary about the dimension and irreducible components (although better than corollary \ref{corol_naive_intro} on this second aspect).

\begin{corollaryintro}[corollaries \ref{corol_main_dim} and \ref{corol_main_irred}, simplified] \label{corol_main_intro}
Assume once again either $d\geq3$ and $n \geq 2^d(d-1)$, or $d=2$ and $n\geq8$. There exist $b'_1, b'_2, b'_3$ that depend on $d$ and $n$ such that as soon as
\[
e > b'_1 g + b'_2 l + b'_3 ,
\]
we have
\begin{align*}
\dim_{J(\scrC) \times W} \Hom^e(\scrC,V|W) & = \xi_e - 1 + l \\
& = (n-d+1)e - n(g-1) - 1 - (n-1)l .
\end{align*}
Over $\C$ this translates as
\begin{equation} \label{corol_main_dimC}
\dim \Hom^e(\scrC,V|W) = (n-d+1)e - (n-1)(g-1) - ((n-1)l - \dim W) .
\end{equation}
Under the same assumption on $e$, if $(\mathcal{W}_i)_i$ is a decomposition of $W$ into irreducible components, then $\big(\Hom^e(\scrC,V|\mathcal{W}_i)\big)_i$ is a decomposition of $\Hom^e(\scrC,V|W)$ into irreducible components.
\end{corollaryintro}

In particular $\Hom^e(\scrC,V|W)$ is irreducible if and only if $W$ is so.

Finally, the Hodge numbers $h^{p,q}(X)$ of a smooth proper complex variety $X$ can be assembled into a generating series called its \emph{Hodge-Deligne polynomial}
\[
\HD(X) := \sum_{p,q \geq0} (-1)^{p+q} h^{p,q}(X) u^p v^q \in \mathbb{Z}[u,v] .
\]
These concepts can be generalised to complex varieties in general (the numbers are then called virtual Hodge numbers), and induce a motivic measure (i.e. a ring morphism) $\HD : K_0\left(\mathrm{Var}_\C\right) \to \mathbb{Z}[u,v]$, as seen in \cite{CLNS18}, chapter 2, proposition 3.2.14. Furthermore $\HD(\mathbf{L}) = uv$, and $\HD$ also provides a motivic measure $\widehat{\mathscr{M}}_\C \to \mathbb{Z}[u,v]\llbracket(uv)^{-1}\rrbracket$. This allows one to retrieve a part of the coefficients of $\HD(X)$ from an approximation of $X$ in motivic rings. In our case, theorem \ref{th_main_intro} provides us with the following partial result, where
\[
\mathscr{P} = \prod_{x\in\CmS} \left( \left(1-\mathbf{L}^{-1}\right) \frac{[V]}{\mathbf{L}^{n-1}} \right)
\]
stands for the motivic Euler product appearing in said theorem. Recall also that
\[
\xi_e - 1 + l = (n-d+1)e - n(g-1) - 1 - (n-1)l .
\]

\begin{corollaryintro} \label{corol_HD_intro}
We have
\[
\HD(\Hom^e(\scrC,V|W)) \equiv (uv)^{\xi_e - 1 + l} \big((u+1)(v+1)\big)^g \left(1-(uv)^{-1}\right)^{\#|\scrS|-1} \cdot \HD(\mathscr{P}) \cdot \HD(W)
\]
modulo monomials $u^p v^q$ with
\[
p+q \leq \dim \Hom^e(\scrC,V|W) + w_{J(\scrC) \times W}\left(\mathcal{E}_e\right) ,
\]
that last weight being the one appearing in theorem \ref{th_main_intro}.
\end{corollaryintro}
See subsection \ref{subsec_HD}. Computing $\HD(\mathscr{P})$ would allow one to obtain the virtual Hodge numbers of these moduli spaces for $p+q$ larger than the above bound, as a function of $\HD(W)$.

\subsection{Structure of the paper}

In section \ref{sect_preli} `Preliminaries' we recall or establish a number of notions and constructions. Among them: properties of the spaces of sections $P_L$ and $P_{|Z}$; construction of vector bundles that represent the latter over a base; presentation of a generalisation of the Grothendieck ring $K_0\left(\mathrm{Var}_S\right)$ and various related properties.

Through section \ref{sect_CM} `Circle method' we study the naive space $\calM_e^W$ (or rather a generalisation of it). After initiating said method, we prove appropriate versions of Weyl's inequality, then progressively rewrite the class of the naive space as a renormalised motivic Euler product, up to a controlled error term. This gives rise to theorem \ref{th_naive_intro} and corollary \ref{corol_naive_intro}.

Finally the goal of section \ref{sect_Hom} is to deduce theorem \ref{th_main_intro} and corollary \ref{corol_main_intro}, which are concerned with the space of main interest, $\Hom^e(\scrC, V | W)$. This is done via motivic inclusion-exclusion, then we detail the bound on the error term and the hypotheses needed for the corollary. We also justify corollary \ref{corol_HD_intro}.

\subsection{Acknowledgements}

I am very grateful to my PhD advisors Margaret Bilu and Antoine Chambert-Loir for suggesting this subject and guiding me through its various steps and convolutions. I would also like to thank Tim Browning, Loïs Faisant and Matthew Hase-Liu for helpful conversations.

No AI was used in the making of this article.

\newpage
\section{Preliminaries}
\label{sect_preli}

\subsection{Integer partitions and configuration spaces} \label{part_conf}

First we fix a handful of notations regarding integer partitions. Recall that a partition of an integer $N\geq0$ is a sequence $\mu = (\mu_i)_{i\geq1}$ with $\mu_i \geq0$ such that $\sum_{i\geq1} \mu_i i = N$; $\mu_i$ is the multiplicity of $i$ in the partition. We write this $\mu \vdash N$, and denote
\[
\Sigma\mu := \sum_{i\geq1} \mu_i i = N .
\]
We also consider the length
\[
|\mu| := \sum_{i\geq1} \mu_i
\]
of $\mu$, as well as
\[
||\mu|| := \#\{ i\geq1 \,|\, \mu_i \geq 1 \} .
\]
The empty partition is $\varnothing = (0)_{i\geq1} \vdash 0$. The sum of two partitions $\mu$ and $\mu'$ is $\mu+\mu' = \nu$ with $\nu_i = \mu_i + \mu'_i$: it can be thought of as their disjoint union. Clearly
\[
\Sigma(\mu+\mu') = \Sigma\mu + \Sigma\mu' \quad\text{and}\quad |\mu+\mu'| = |\mu| + |\mu'| .
\]
Consider two partitions $\mu$ and $\nu$: if there exists $\mu'$ such that $\nu = \mu+\mu'$ we say that $\mu$ is a subpartition of $\nu$ and denote $\mu\leq\nu$. In this case $\mu'$ (which is unique) will be denoted by $\nu-\mu$. We obtain a structure of partially ordered commutative monoid. We denote by $[i]$ the partition $\mu$ of $i$ with $\mu_i = 1$ and other coefficients zero. We also write $j[1] := [1] + \cdots + [1]$ with $j$ terms, that is $\mu$ with $\mu_1 = j$ and other coefficients zero.

Let $g \in k[t]$ be a nonzero polynomial over a field: over a splitting field write $g = c\prod_{j=1}^r (t-x_j)^{n_j}$ with $c\in k^\times$ and pairwise distinct $x_j$. We define the partition of $g$ to be $\pt g = \mu$ where $\mu_i = \#\left\{ 1\leq j\leq r \,\left|\, n_j = i \right.\right\}$. Likewise, the partition of an effective divisor $Z = \sum_{j=1}^r n_j x_j$ (with pairwise distinct $x_j$) on $\scrC$ is $\pt Z = \mu$ with $\mu_i = \#\left\{ 1\leq j\leq r \,\left|\, n_j = i \right.\right\}$ too.

Now let us recall the definition of symmetric powers. Let $S$ be a variety over a field $k$, and $X$ be a quasi-projective $S$-variety. First of all,
\[
\Sym_{/S}^N(X) := X^N / \mathfrak{S}_N
\]
is the quotient of $X^N$ (fibered over $S$) by the symmetric group. For $\mu \vdash N$, $\Sym_{/S}^\mu(X)$ is then the subscheme of $\Sym_{/S}^N(X)$ whose elements have partition $\mu$. It can also be constructed as follows. Take the product $\prod_{i\geq1} X^{\mu_i}$ over $S$, and denote by $\big(\prod_{i\geq1} X^{\mu_i}\big)_*$ the complement of the big diagonal, i.e. keep only the tuples with pairwise distinct coordinates. Then
\[
\Sym_{/S}^\mu(X) = \Bigg( \prod_{i\geq1} X^{\mu_i} \Bigg)_* \ / \mathfrak{S}_\mu
\]
where $\mathfrak{S}_\mu := \prod_{i\geq1} \mathfrak{S}_{\mu_i}$ is the (finite) symmetric group associated to $\mu$. Note that
\[
\dim_S \Sym_{/S}^N(X) = N \dim_S X \quad\text{and}\quad \dim_S \Sym_{/S}^\mu(X) = |\mu| \dim_S X .
\]
Another classical construction is the mixed product
\begin{equation} \label{def_mixed_sym}
\Sym_{/S}^{\mu, \mu'}(X) := \left\{\left. (Z,Z') \in \Sym_{/S}^\mu(X) \times_S \Sym_{/S}^{\mu'}(X) \,\right|\, Z \cap Z' = \varnothing \right\} .
\end{equation}
When $S = \Spec k$ with $k$ clear from the context, we will drop the $/S$ from the notations.

Finally we construct more general \emph{configuration spaces}. Consider a variety $S$ over a field $k$, a quasi-projective variety $X$ over $S$, and a family $(X_i)_{i\geq1}$ of $X$-varieties. For an integer partition $\mu$ we take the product $\prod_{i\geq1} X_i^{\mu_i}$ over $S$, and remove the `big diagonal over $X$', that is we only keep the tuples whose coordinates map to pairwise distinct points of $X$, to obtain $\big(\prod_{i\geq1} X_i^{\mu_i}\big)_{*,X}$. The configuration space we are interested in is then\footnote{This is the $S^\mu\left( (X_i)_{i\geq1} / S \right)$ from \cite{Bil23} (paragraph 3.1.3), the $C_{X/S}^\mu \left((X_i)_{i\geq1}\right)$ from \cite{BH21} (subsection 3.2), or the $\Conf^\mu(X_i)_{i\geq1}$ from \cite{BB25} (subsection 2.4).}
\begin{equation} \label{def_conf}
\Conf_{X/S}^\mu \left((X_i)_{i\geq1}\right) := \Bigg( \prod_{i\geq1} X_i^{\mu_i} \Bigg)_{*,X} \ / \mathfrak{S}_\mu .
\end{equation}
When all the $X_i$ are equal to a common $Y \to X$, we write $\Conf_{X/S}^\mu(Y)$. In particular we retrieve $\Sym_{/S}^\mu(X) = \Conf_{X/S}^\mu(X)$.

\subsection{Sections of certain vector bundles}

\subsubsection{Sections for a given line bundle $L$ or effective divisor $Z$}
\label{global_sect}

This paragraph is a summary of properties of spaces of sections such as our $P_L$ and $P_{|Z}$. Let $Z$ be an effective divisor on $\scrC$, or equivalently a subscheme of $\scrC$ of dimension $\leq 0$. First we justify the isomorphism $H^0(Z, L|_Z) \simeq H^0(Z, \trivlb|_Z) = P_{|Z}$ from remark \ref{sect_loc}, for any line bundle $L$ whose base locus doesn't intersect $Z$.

\begin{proof}[Proof of remark \ref{sect_loc}.]
Write $Z = \sum_{j=1}^r n_j x_j$ with the $x_j \in \scrC$ pairwise distinct and $n_j\geq1$, so that $H^0(Z, L|_Z) = \prod_{j=1}^r H^0(n_j x_j, L|_{n_j x_j})$. Since line bundles are locally trivial and $x_i$ is no basepoint of $L$, one has $L|_{n_j x_j} \simeq \trivlb|_{n_j x_j}$, and as a consequence
\[
H^0(n_j x_j, L|_{n_j x_j}) \simeq H^0(n_j x_j, \trivlb|_{n_j x_j}) :
\]
we deduce $H^0(Z, L|_Z) \simeq P_{|Z}$.
\end{proof}

In relation with this argument, note that $H^0(n_j x_j, \trivlb|_{n_j x_j}) \simeq \C[t]/t^{n_j}$, therefore if $Z = \sum_{j=1}^r n_j x_j$ as above then we have a canonical isomorphism
\begin{equation} \label{PZ_structure}
P_{|Z} \simeq \prod_{j=1}^r \C[t]/t^{n_j} .
\end{equation}
This has three consequences. The first one is that $\dim P_{|Z} = \deg Z$. The second one is the lemma below.
\begin{lemma} \label{PZ_isom}
If $Z$ and $Z'$ are two effective divisors with the same partition, then there exists a (non-canonical) isomorphism of $\C$-algebras $P_{|Z} \simeq P_{|Z'}$.
\end{lemma}
The last consequence of (\ref{PZ_structure}) is the following. Pick arbitrary pairwise distinct $c_j \in \C$: $\C[t]/t^{n_j}$ and $\C[t]/(t-c_j)^{n_j}$ are isomorphic as $\C$-algebras, so $P_{|Z} \simeq \prod_{j=1}^r \C[t]/(t-c_j)^{n_j}$, and the Chinese remainder theorem gives the isomorphism of $\C$-algebras
\begin{equation} \label{PZ_quotient}
P_{|Z} \simeq \C[t]/h
\end{equation}
with $h = \prod_{j=1}^r(t-c_j)^{n_j}$. Note that the polynomial $h$ has the same partition as $Z$ (and in fact can be replaced with any other polynomial with the same partition).

Now let us look at global sections : consider a line bundle $L$ of degree $e$ over $\scrC$. By the Riemann-Roch theorem, as long as $e \geq 2g-1$,
\begin{equation} \label{RR}
\dim P_L = h^0(\scrC, L) = e-g+1 .
\end{equation}
When $Z$ is an effective divisor on $\scrC$ that does not intersect the base locus of $L$ (for instance when $e\geq2g$, which implies that $L$ is basepoint-free), one has the usual cohomology exact sequence (where we use $H^0(Z, L|_Z) \simeq P_{|Z}$)
\begin{equation} \label{exseq_H1}
\xymatrix{
0 \ar[r] & P_{L(-Z)} \ar[r]^-{\subset} & P_L \ar[r]^-{\cdot|_Z} & P_{|Z} \ar[r] & H^1(\scrC, L(-Z))
} .
\end{equation}
If $\deg L(-Z) \geq 2g-1$, i.e. $\deg Z \leq e-2g+1$, then $H^1(\scrC, L(-Z))=0$ and we get the following short exact sequence:
\begin{equation} \label{exseq_ed}
\xymatrix{
0 \ar[r] & P_{L(-Z)} \ar[r]^-{\subset} & P_L \ar[r]^-{\cdot|_Z} & P_{|Z} \ar[r] & 0
}
\end{equation}
(which is split since we are looking at vector spaces).

\begin{example} \label{expl}
Let's consider the base case $\scrC = \projsp^1$, $L = \mathcal{O}(e)$ with $e\geq0$. Take an effective divisor $Z$ of degree $\delta \leq e+1$: it can be written
\[
Z = \sum_{j=1}^r n_j [a_j : b_j]
\]
with $[a_j : b_j] \in \projsp^1$ pairwise distinct and $n_j\geq1$. As we already said, $\C[t]_{\leq e}$ designates the polynomials of degree at most $e$. With the representation (\ref{PZ_quotient}) for $P_{|Z}$, the previous exact sequence (\ref{exseq_ed}) becomes in this case
\[
\xymatrix{
0 \ar[r] & h \cdot \C[t]_{\leq e-\delta} \ar[r]^-{\subset} & \C[t]_{\leq e} \ar[r]^-{\bmod h} & \C[t]/h \ar[r] & 0
}
\]
with a certain $h = \prod_{j=1}^r(t-c_j)^{n_j}$. The choice of $h$ is not canonical in general, but it can be made so here. Naturally an appropriate choice for the roots is $c_j = \frac{b_j}{a_j}$ when $a_j\neq0$, and any distinct $c_j$ for $a_j=0$ (i.e. for the point $\infty \in \projsp^1$). In order to have something completely canonical, we must switch to the homogeneous point of view: $\C[t]_{\leq e}$ is isomorphic to the space $\C[t_0, t_1]_e^\mathrm{hg}$ of homogeneous polynomials of degree $e$ in $t_0, t_1$. The divisor $Z$ is then associated to the following degree $\delta$ homogeneous polynomial, defined up to a multiplicative constant:
\[
h_Z(t_0, t_1) := \prod_{j=1}^r (a_j t_1 - b_j t_0)^{n_j} .
\]
With this the above exact sequence is isomorphic to
\[
\xymatrix{
0 \ar[r] & h_Z \cdot \C[t_0, t_1]_{e-\delta}^\mathrm{hg} \ar[r]^-{\subset} & \C[t_0, t_1]_e^\mathrm{hg} \ar[r] & \C[t_0, t_1]_e^\mathrm{hg}/h_Z \ar[r] & 0
} ,
\]
the second non-trivial arrow being $\bmod\, h_Z$. In any case note that these exact sequences represent the Euclidean division by $h$: the quotient (multiplied by $h$) is on the left and the remainder is on the right. This key example leads us to call the exact sequence (\ref{exseq_ed}) the \emph{Euclidean division exact sequence} for any curve $\scrC$.
\end{example}

Now let $Z_1$ and $Z_2$ be effective divisors: we will consider their intersection $Z_1 \cap Z_2$ and union $Z_1 \cup Z_2$ as effective divisors, whose multiplicities are respectively the minima and the maxima of those of $Z_1$ and $Z_2$ (in light of the previous example, one can think of this as a gcd and an lcm). We define
\[
P_{-Z_2|Z_1} := H^0(Z_1, \trivlb(-Z_2)|_{Z_1}) ,
\]
which is an ideal of $P_{|Z_1}$. Remark that as with $P_{|Z_1}$, one can replace the $\trivlb$ in this definition by any line bundle $L$ whose base locus doesn't intersect $Z_1$, and that for this reason $P_{-Z_2|Z_1} = P_{-(Z_1 \cap Z_2)|Z_1}$. Going back to the polynomial representation (\ref{PZ_quotient}), let $h_1, h_2, h_\cap \in \C[t]$ represent respectively $Z_1, Z_2, Z_1 \cap Z_2$ in a compatible way (i.e. each has same partition as the associated divisor, and $h_\cap = \gcd(h_1, h_2)$). The corresponding isomorphism transforms
\[
P_{-Z_2|Z_1} = P_{-(Z_1 \cap Z_2)|Z_1} \subset P_{|Z_1}
\]
into
\[
h_2 \cdot \C[t]/h_1 = h_\cap \cdot \C[t]/h_1 \subset \C[t]/h_1 .
\]
More precisely on the polynomial side we have the exact sequence
\[
\xymatrix{
0 \ar[r] & h_2 \cdot \C[t]/h_1 \ar[r]^-{\subset} & \C[t]/h_1 \ar[r] & \C[t]/h_\cap \ar[r] & 0
}
\]
(where the second non-trivial arrow is $\bmod\, h_\cap$), which translates to
\begin{equation} \label{exseq_ed_loc}
\xymatrix{
0 \ar[r] & P_{-Z_2|Z_1} \ar[r]^-{\subset} & P_{|Z_1} \ar[r] & P_{|Z_1 \cap Z_2} \ar[r] & 0
}
\end{equation}
(this time the second non-trivial arrow is $\cdot|_{Z_1 \cap Z_2}$). It can be seen as a \emph{local} or restricted version of the Euclidean division exact sequence (\ref{exseq_ed}). A consequence is that
\[
\dim P_{-Z_2|Z_1} = \deg Z_1 - \deg(Z_1 \cap Z_2) .
\]
To complete these preliminaries, let us give a version of exact sequence (\ref{exseq_ed}) with two divisors.
\begin{lemma} \label{exseq_ed_vanish}
Let $Z_1$ and $Z_2$ be two effective divisors such that $\deg(Z_1 \cup Z_2) \leq e-2g+1$. Suppose that $Z_1 \cup Z_2$ doesn't intersect the base locus of $L$. The Euclidean division exact sequence (\ref{exseq_ed}) restricts to
\[
\xymatrix{
0 \ar[r] & P_{L(-(Z_1 \cup Z_2))} \ar[r]^-{\subset} & P_{L(-Z_2)} \ar[r]^-{\cdot|_{Z_1}} & P_{-Z_2|Z_1} \ar[r] & 0
} .
\]
\end{lemma}
In practice we will apply it with $L$ basepoint-free.
\begin{proof}
We need to prove that in this context the restriction to $Z_1$ is surjective and has kernel $P_{L(-(Z_1 \cup Z_2))}$. The second part is direct: if $\mathbf{v} \in P_{L(-Z_2)}$, by definition $\mathbf{v}|_{Z_2} = 0$: then $\mathbf{v}|_{Z_1} = 0$ if and only if $\mathbf{v}|_{Z_1 \cup Z_2} = 0$, if and only if $\mathbf{v} \in P_{L(-(Z_1 \cup Z_2))}$. For the surjectivity it is enough to show that $\mathrm{Im}(\cdot|_{Z_1})$ and $P_{-Z_2|Z_1}$ have same dimension. Since $\deg Z_2 \leq \deg(Z_1 \cup Z_2) \leq e-2g+1$ the Riemann-Roch theorem (\ref{RR}) can be applied to both $Z_2$ and $Z_1 \cup Z_2$, so that $\dim P_{L(-Z_2)} = e-g+1 - \deg Z_2$ and $\dim P_{L(-(Z_1 \cup Z_2))} = e-g+1 - \deg(Z_1 \cup Z_2)$. Hence
\begin{align*}
\dim\mathrm{Im}(\cdot|_{Z_1}) & = \dim P_{L(-Z_2)} - \dim P_{L(-(Z_1 \cup Z_2))} \\
& = \deg(Z_1 \cup Z_2) - \deg Z_2 \\
& = \deg Z_1 - \deg(Z_1 \cap Z_2) \\
& = \dim P_{-Z_2|Z_1} ,
\end{align*}
where the transition between lines 2 and 3 comes from the identity $Z_1 + Z_2 = (Z_1 \cup Z_2) + (Z_1 \cap Z_2)$. As a consequence $\cdot|_{Z_1}$ is surjective.
\end{proof}

\subsubsection{Poincaré bundles}

\label{Poincare}

Let $e \geq 2g$ be a fixed integer: the condition ensures that all $L\in\Pic^e(\scrC)$ are basepoint-free. The circle method will essentially be run using spaces of global sections such as $P_{L(-Z)}$, however it is essential for the application to spaces of morphisms in section \ref{sect_Hom} and the process of \emph{motivic inclusion-exclusion} that will be used, to perform it \emph{over a base,} that is both over $L \in \Pic^e(\scrC)$ and over some space of divisors $Z$.

Recall from \cite{ACGH85} (chapter IV, paragraph 2) the existence of the Poincaré line bundle
\[
\Poinc_e \quad\text{over}\quad \scrC \times \Pic^e(\scrC)
\]
with the property that for any $L \in \Pic^e(\scrC)$, $\Poinc_e|_{\scrC\times\{L\}} = L$ as a line bundle over $\scrC\times\{L\} \simeq \scrC$.
\begin{remark} \label{constr_fibr_Hom}
Before going on let us describe $\Hom^e(\scrC,V)$ and $\Hom^e(\scrC,V|W)$ in a way that is coherent with the description of their fibres (\ref{def_Hom_L}) and (\ref{def_Hom_LW}) over $\Pic^e(\scrC)$ that was given in the introduction. With the projection $\pi : \scrC \times \Pic^e(\scrC) \to \Pic^e(\scrC)$, the vector bundle $\pi_*\Poinc_e$ represents the global sections of all line bundles of degree $e$ in the following way. For $L \in \Pic^e(\scrC)$,
\begin{align*}
(\pi_*\Poinc_e)_L & = \lim_{\mathcal{U}\ni L} H^0(\mathcal{U}, \pi_*\Poinc_e) \\
& = \lim_{\mathcal{U}\ni L} H^0(\scrC\times\mathcal{U}, \Poinc_e) \\
& = H^0\!\left(\scrC\times\{L\}, \Poinc_e|_{\scrC\times\{L\}}\right) \\
& = H^0(\scrC,L) \\
& = P_L .
\end{align*}
As a consequence the points of $\pi_*\Poinc_e$ can be seen as $(L,\mathbf{u})$ with $\mathbf{u} \in P_L$. We may then write
\[
\Hom^e(\scrC,V) = \left\{\left. (L, \bar{\mathbf{u}}) \in (\pi_*\Poinc_e)^{n+1} /\mathbb{G}_m \,\right|\, \bar{\mathbf{u}} \in \Hom^L(\scrC,V) \right\} ,
\]
where the $(\cdot)^{n+1}$ is to be understood as an iterated fibered product over $\Pic^e(\scrC)$, and $\mathbb{G}_m$ acts on $\mathbf{u}$. The same holds for $\Hom^e(\scrC,V|W)$ and $\Hom^L(\scrC,V|W)$.
\end{remark}

In order to have bundles over
\[
\Pic^e(\scrC) \times \Sym^\nu(\scrC)
\]
for some integer partition $\nu$ and not only over $\Pic^e(\scrC)$, we consider the projection
\[
\pr_{1,2} : \scrC \times \Pic^e(\scrC) \times \Sym^\nu(\scrC) \longrightarrow \scrC \times \Pic^e(\scrC)
\]
and the line bundle
\[
\calL_e := \pr_{1,2}^* \Poinc_e .
\]
It satisfies
\begin{equation}
\calL_e|_{\scrC\times\{(L,Z)\}} = \Poinc_e|_{\scrC\times\{L\}} = L
\end{equation}
for any $L \in \Pic^e(\scrC)$ and $Z\in\Sym^\nu(\scrC)$. Next we want a second bundle that provides $L(-Z)$ under the same restriction. In \cite{ACGH85}, in order to define the Poincaré line bundles, they first construct the \emph{universal effective divisor of degree $\delta$ on $\scrC$}, which is a divisor $\Delta_\delta$ of $\scrC \times \Sym^\delta(\scrC)$ whose restriction to $\scrC \times \{Z\} \simeq \scrC$ is precisely $Z$. With $\delta = \Sigma\nu$ we can restrict it to $\Sym^\nu(\scrC) \subset \Sym^\delta(\scrC)$, leading to a second divisor $\Delta_\nu$ ; then we pull back the latter through $\pr_{1,3}$ to get a third divisor $\Delta_\nu'$ over $\scrC \times \Pic^e(\scrC) \times \Sym^\nu(\scrC)$, whose restriction to $\scrC \times \{(L,Z)\}$ is again $Z$. Now put
\[
\calL_{e(-\nu)} := \calL_e(-\Delta_\nu') :
\]
by construction we have
\begin{equation}
\calL_{e(-\nu)}|_{\scrC\times\{(L,Z)\}} = L(-Z)
\end{equation}
for all $(L,Z)$. Finally we consider
\[
\calL_{|\nu} := \calL_e / \calL_{e(-\nu)} :
\]
we have
\begin{align}
\calL_{|\nu}|_{\scrC\times\{(L,Z)\}} & = \calL_e|_{\scrC\times\{(L,Z)\}} \left/ \calL_{e(-\nu)}|_{\scrC\times\{(L,Z)\}} \right. \nonumber \\
& = L/L(-Z) \nonumber \\
& = (i_Z)_* L|_Z
\end{align}
where $i_Z : Z \hookrightarrow \scrC$. As seen in the proof of remark \ref{sect_loc}, $L|_Z \simeq \mathcal{O}_Z$, hence the absence of $e$ in the notation. We have the exact sequence
\begin{equation} \label{exseq_calL}
\xymatrix{
0 \ar[r] & \calL_{e(-\nu)} \ar[r] & \calL_e \ar[r] & \calL_{|\nu} \ar[r] & 0
} ,
\end{equation}
which represents the family of exact sequences
\[
\xymatrix{
0 \ar[r] & L(-Z) \ar[r] & L \ar[r] & (i_Z)_* L|_Z \ar[r] & 0
} .
\]

Since we are interested in the global sections of the latter, the analog over our base is as follows. Let
\[
\pi : \scrC \times \Pic^e(\scrC) \times \Sym^\nu(\scrC) \longrightarrow \Pic^e(\scrC) \times \Sym^\nu(\scrC)
\]
be the projection that forgets the point of $\scrC$: we take the notations
\begin{equation} \label{def_calP}
\calP_{e(-\nu)} := \pi_*\calL_{e(-\nu)} , \qquad \calP_e := \pi_*\calL_e , \qquad \calP_{|\nu} := \pi_*\calL_{|\nu} .
\end{equation}
Just like in remark \ref{constr_fibr_Hom}, the fibres over a closed point $(L,Z)$ of $\Pic^e(\scrC) \times \Sym^\nu(\scrC)$ are
\[
(\calP_e)_{L,Z} = P_L, \qquad (\calP_{e(-\nu)})_{L,Z} = P_{L(-Z)} ,
\]
and
\[
(\calP_{|\nu})_{L,Z} = H^0(\scrC, (i_Z)_* L|_Z) = H^0(Z, L|_Z) \simeq P_{|Z} .
\]
\begin{remark} \label{def_gen_fibres}
As for the generic fibres, we can adopt the converse point of view: the spaces of sections are not defined yet so we use the above as a definition. More precisely, let $\pr_1$ and $\pr_2$ be the projections on $\Pic^e(\scrC)$ and $\Sym^\nu(\scrC)$ respectively:
\begin{itemize}
\item for all $L \in \Pic^e(\scrC)$ non closed we put $P_L := (\pr_{1*}\calP_e)_L$;
\item for all $(L,Z) \in \Pic^e(\scrC) \times \Sym^\nu(\scrC)$ non closed
\footnote{Not all generic point is a couple $(L,Z)$, but for readability we will still write them this way.}
we put $P_{L(-Z)} := (\calP_{e(-\nu)})_{L,Z}$;
\item for all $Z \in \Sym^\nu(\scrC)$ non closed we put $P_{|Z} := (\pr_{2*}\calP_{|\nu})_Z$.
\end{itemize}
Note that with these conventions, $(\calP_e)_{L,Z} = P_L \otimes_{\kappa(L)} \kappa(L,Z)$ and $(\calP_{|\nu})_{L,Z} = P_{|Z} \otimes_{\kappa(Z)} \kappa(L,Z)$.
\end{remark}

\begin{lemma}[Euclidean division over a base] \label{exseq_fam_nu}
If $\Sigma\nu \leq e-2g+1$ then we have the exact sequence of vector bundles
\[
\xymatrix{
0 \ar[r] & \calP_{e(-\nu)} \ar[r] & \calP_e \ar[r] & \calP_{|\nu} \ar[r] &0
} .
\]
\end{lemma}
\begin{proof}
The functor $\pi_*$ is left-exact so exact sequence (\ref{exseq_calL}) gives rise to
\[
\xymatrix{
0 \ar[r] & \calP_{e(-\nu)} \ar[r] & \calP_e \ar[r] & \calP_{|\nu}
} .
\]
For every closed point $(L,Z)$ the fibre is
\[
\xymatrix{
0 \ar[r] & P_{L(-Z)} \ar[r] & P_L \ar[r] & P_{|Z}
} ,
\]
which admits a $0$ on its right by (\ref{exseq_ed}), because $\deg Z = \Sigma\nu \leq e-2g+1$. Since this is true for all closed points, for dimension reasons the exact sequences corresponding to the generic fibres can be completed with a $0$ as well. An application of \cite{Har77}, chapter II, exercise 1.2 concludes.
\end{proof}

If $(L,Z)$ is a closed point, we know how $P_L$ and $P_{|Z}$ behave with respect to the product: remark \ref{sect_loc} yields a commutative square
\begin{equation} \label{sq_prod}
\xymatrix{
P_L^2 \ar[r]^{\cdot|_Z} \ar[d]_{\times} & P_{|Z}^2 \ar[d]^{\times} \\
P_{L^{\otimes 2}} \ar[r]_{\cdot|_Z} & P_{|Z}
} .
\end{equation}
We now look at the generic case. Let $(L',Z')$ be a non closed point of residue field $\kappa'$: by construction and by remark \ref{def_gen_fibres}, we have a map $P_{L'}^2 \to P_{|Z'}^2$. By setting
\[
\tau_2 : \left\{\begin{array}{rcl}
\Pic^e(\scrC) \times \Sym^\nu(\scrC) & \longrightarrow & \Pic^{2e}(\scrC) \times \Sym^\nu(\scrC) \\
(L, \scrT) & \longmapsto & (L^{\otimes 2}, \scrT)
\end{array}\right.
\]
and considering $\tau_2^* \calP_{2e} \to \tau_2^* \calP_{|\nu} = \calP_{|\nu}$, we also get $P_{L'{}^{\otimes 2}} \to P_{|Z'}$. For dimensional reasons, these two maps of $\kappa'$-vector spaces are respectively isomorphic to the first and second line of (\ref{sq_prod}) tensored with $\kappa'$, which allows us to simply transport the two product arrows and obtain a second commutative square
\[
\xymatrix{
P_{L'}^2 \ar[r]^{\cdot|_{Z'}} \ar[d]_{\times} & P_{|Z'}^2 \ar[d]^{\times} \\
P_{L'{}^{\otimes 2}} \ar[r]_{\cdot|_{Z'}} & P_{|Z}
} ,
\]
isomorphic to (\ref{sq_prod}) tensored with $\kappa'$. In particular, via the projection that forgets $L \in \Pic^e(\scrC)$, this endows $P_{|Z'}$ with a structure of $\kappa(Z')$-algebra that satisfies the property below. Note that by lemma \ref{PZ_isom}, having taken a different closed $Z \in \Sym^\nu(\scrC)$ would have resulted in the same algebra structure.
\begin{lemma} \label{PZ_gen_isom}
If $Z,Z'$ are points of $\Sym^\nu(\scrC)$, the former being closed, then we have an isomorphism of $\kappa(Z')$-algebras
\[
P_{|Z'} \simeq P_{|Z} \otimes \kappa(Z') .
\]
\end{lemma}

Let us state a second lemma that we will use later.
\begin{lemma} \label{PZ_gen_decomp}
For all $(Z_1, Z_2) \in \Sym^{\nu_1, \nu_2}(\scrC)$, with $\kappa = \kappa(Z_1, Z_2)$,
\[
(\cdot|_{Z_1}, \cdot|_{Z_2}) : P_{|Z_1+Z_2} \overset{\sim}{\longrightarrow} \left(P_{|Z_1} \otimes_{\kappa(Z_1)} \kappa\right) \times \left(P_{|Z_2} \otimes_{\kappa(Z_2)} \kappa\right)
\]
is an isomorphism of $\kappa$-algebras.
\end{lemma}
\begin{proof}
This map is well-defined and is an isomorphism for all closed points $(Z_1, Z_2)$, so it has to be an isomorphism of $\kappa$-vector spaces for the non closed points as well. Furthermore it commutes with the product from the previous construction.
\end{proof}

Now that we have the necessary properties, let us add something to the base we are working over. In fact we will need to add vanishing conditions to the sections we will consider, in the following sense. Recall our $\scrS$ fixed in (\ref{def_scrS}). If $\scrT$ is an effective divisor on $\scrC$, we denote its support by $|\scrT|$. We fix an integer partition $\mu$: to every $\scrT\in\Sym^\mu(\CmS)$ we associate a condition locus
\begin{equation} \label{def_ST}
\scrST := \scrS+|\scrT| .
\end{equation}
By slightly adapting the construction resulting in lemma \ref{exseq_fam_nu}, if $l+|\mu| \leq e - 2g + 1$ we can construct an exact sequence
\begin{equation} \label{exseq_fam_Smu}
\xymatrix{
0 \ar[r] & \calP_{e(-\scrS_\mu)} \ar[r] & \calP_e \ar[r] & \calP_{|\scrS_\mu} \ar[r] &0
}
\end{equation}
of vector bundles over $\Pic^e(\scrC) \times \Sym^\mu(\CmS)$ whose fibre at $(L,\scrT)$ is
\[
\xymatrix{
0 \ar[r] & P_{L(-\scrST)} \ar[r] & P_L \ar[r] & P_{|\scrST} \ar[r] &0
} .
\]
Actually under the assumption
\[
l + |\mu| + \Sigma\nu \leq e - 2g + 1
\]
we can put both version of the exact sequence in interaction by pulling back everything over
\[
\Pic^e(\scrC) \times \Sym^\mu(\CmS) \times \Sym^\nu(\scrC) ,
\]
which provides a diagram
\footnote{The `$e(-\scrS_\mu)$', `$\scrS_\mu \cup \nu$', etc. are only notations meant to reflect what happens at the fibers, and have no meaning by themselves.}
\footnote{For simplicity we keep the notations $\calP_e$ etc. regardless of which base we are over, as opposed to adding some `$\pr^*$'.}
\begin{equation} \label{exseq_9}
\xymatrix{
& 0 \ar[d] & 0 \ar[d] & 0 \ar[d] &\\
0 \ar[r] & \calP_{e(-(\scrS_\mu \cup \nu))} \ar[r]\ar[d] & \calP_{e(-\nu)} \ar[r]\ar[d] & \calP_{-\nu|\scrS_\mu} \ar[r]\ar[d] & 0 \\
0 \ar[r] & \calP_{e(-\scrS_\mu)} \ar[r]\ar[d] & \calP_e \ar[r]\ar[d] & \calP_{|\scrS_\mu} \ar[r]\ar[d] & 0 \\
0 \ar[r] & \calP_{-\scrS_\mu|\nu} \ar[r]\ar[d] & \calP_{|\nu} \ar[r]\ar[d] & \calP_{|\scrS_\mu \cap \nu} \ar[r]\ar[d] & 0 \\
& 0 & 0 & 0 &
}
\end{equation}
whose fibre over $(L,\scrT,Z)$ is
\footnote{In the fashion of remark \ref{def_gen_fibres}, this lets us define $P_{-\scrST|Z}$ and $P_{-Z|\scrST}$ for generic points.}
\begin{equation} \label{exseq_9_fibre}
\xymatrix{
& 0 \ar[d] & 0 \ar[d] & 0 \ar[d] &\\
0 \ar[r] & P_{L(-(\scrST \cup Z))} \ar[r]\ar[d] & P_{L(-Z)} \ar[r]\ar[d] & P_{-Z|\scrST} \ar[r]\ar[d] & 0 \\
0 \ar[r] & P_{L(-\scrST)} \ar[r]\ar[d] & P_L \ar[r]\ar[d] & P_{|\scrST} \ar[r]\ar[d] & 0 \\
0 \ar[r] & P_{-\scrST|Z} \ar[r]\ar[d] & P_{|Z} \ar[r]\ar[d] & P_{|\scrST \cap Z} \ar[r]\ar[d] & 0 \\
& 0 & 0 & 0 &
}
\end{equation}
tensored with $\kappa(L,\scrT,Z)$. Both diagrams are commutative and have exact rows and columns.
This is a consequence of the construction for \ref{exseq_fam_nu} and (\ref{exseq_fam_Smu}), as well as similar constructions based on exact sequences \ref{exseq_ed_vanish} and (\ref{exseq_ed_loc}).

\subsection{Jet schemes}
\label{jet_sch}

Let us say a few words about jet schemes. Let $k$ be a field, $X$ a scheme of finite type over $k$, and $m \in \mathbb{N}$. Recall, for instance from section 2 of \cite{Ish07}, that the $m$-th jet scheme of $X$ is the finite type $k$-scheme $\mathcal{J}_mX$ such that
\[
\mathcal{J}_mX(A) = X\left(A[t]/t^{m+1}\right)
\]
for any finite type $k$-algebra $A$. Generally speaking, one expects $\mathcal{J}_mX$ to have dimension $(m+1) \dim X$.

We mention them because they appear in two ways in this work (over $k = \C$). One is via the $\Lambda_N(x)$ that appear in theorem \ref{th_naive_intro} and that we will work with in paragraph \ref{locdens}. They are defined as the spaces of solutions of $f(\mathbf{u}) = 0$ over $Nx$ for any closed $x \in \scrC$:
\begin{equation} \label{def_Lambda}
\Lambda_N(x) := \left\{\left. \mathbf{u} \in P_{|Nx}^{n+1} \,\right|\, f(\mathbf{u}) = 0 \right\} .
\end{equation}
Let $\tilde V \subset \mathbb{A}_\C^{n+1}$ be the affine hypersurface defined by $f$ (i.e. the universal torsor above $V \subset \projsp_\C^n$). Because $P_{|Nx} \simeq \C[t]/t^N$ as $\C$-algebras, the $A$-points of $\Lambda_N(x)$ are $\tilde V \left(A[t]/t^N\right)$ for any $\C$-algebra $A$, and
\[
\Lambda_N(x) \simeq \mathcal{J}_{N-1} \tilde V
\]
as long as $N\geq1$. We know that this is irreducible of dimension $nN$ from corollary 6.15 of \cite{BB25}, because $n\geq d$ (beware that once again, our $n$ is their $n-1$).

The other jet scheme we consider is naturally $\Hom((m+1)p,V)$ with $m\geq0$ and $p \in \scrC$, that we described in (\ref{Hom_lp}). In fact $(m+1)p = \Spec\left( \C[t]/t^{m+1} \right)$ and
\[
\Hom((m+1)p,V) = \mathcal{J}_m V .
\]
Let us justify the description with $P_{|(m+1)p}$. For any $\C$-algebra of finite type $A$ we have
\begin{align*}
V\left(A[t]/t^{m+1}\right) & = \left\{\left. [\mathbf{g}] \in \projsp^n\!\left(A[t]/t^{m+1}\right) \,\right|\, f(\mathbf{g}) = 0 \right\} \\
& = \left\{ \mathbf{g} \in \left(A[t]/t^{m+1}\right)^{n+1} \,\left|\,
\begin{array}{l}
\exists i,\, g_i \in \left(A[t]/t^{m+1}\right)^\times \\
f(\mathbf{g}) = 0
\end{array} \right.\right\} \left/ \left(A[t]/t^{m+1}\right)^\times \right. .
\end{align*}
The condition $g_i \in \left(A[t]/t^{m+1}\right)^\times$ is equivalent to $g_i(0) \in A^\times$. Since $P_{|(m+1)p} \simeq \C[t]/t^{m+1}$ as $\C$-algebras (which sends $w_i(p)$ to $g_i(0)$), the above consists of the $A$-points of
\[
\left\{ \mathbf{w} \in P_{|(m+1)p}^{n+1} \,\left|\,
\begin{array}{l}
\exists i,\, w_i(p) \in \mathbb{G}_m \\
f(\mathbf{w}) = 0
\end{array} \right.\right\} \left/ \mathbb{G}_{m, P_{|(m+1)p}} \right. ;
\]
focusing on the geometric points, the $w_i(p) \in \mathbb{G}_m$ can be replaced with $w_i(p)\neq 0$, and we obtain (\ref{Hom_lp}) with $m+1 = l_p$.

\begin{remark} \label{dim_jet}
Note that $\Hom((m+1)p,V) = \mathcal{J}_m V$ is the quotient of a dense open of $\mathcal{J}_m \tilde V$ by $\mathbb{G}_{m, P_{|(m+1)p}}$. In particular we mentioned that $\dim \mathcal{J}_m \tilde V = n(m+1)$, while $\dim \mathbb{G}_{m, P_{|(m+1)p}} = m+1$, thus
\[
\dim \mathcal{J}_m V = (n-1)(m+1) .
\]
\end{remark}

\subsection{Grothendieck rings}

\subsubsection{Grothendieck ring of varieties with exponentials}
\label{parag_GR}

Let $S$ be a noetherian scheme. In the introduction we defined the usual Grothendieck ring of varieties $K_0\left(\mathrm{Var}_S\right)$ and its localisation $\mathscr{M}_S$ at $\mathbf{L} := \left[\mathbb{A}_S^1\right]$ (we also mentioned a weight topology, which will be back in a further paragraph). Now we construct a more general ring, that is the Grothendieck ring of varieties with exponentials, which we will use for computations throughout a large part of the paper.

To begin with, a \emph{variety with exponential} over $S$ is a couple $(X,u)$ where $X$ is an $S$-variety and $u : X \to \mathbb{A}_{\mathbb{Z}}^1$ is a map. A morphism from $(X,u)$ to $(Y,v)$ is an $S$-morphism $\varphi : X \to Y$ such that $v \circ \varphi = u$. We define the \emph{Grothendieck ring of varieties with exponentials over $S$}, denoted $K_0\left(\mathrm{ExpVar}_S\right)$, first as an abelian group by generators and relations. The generators are the isomorphism classes of $S$-varieties with exponentials, and the relations are of two types: the \emph{cut-and-paste} relations
\[
(X,u) - (Y,u|_Y) - (X \smallsetminus Y, u|_{X \smallsetminus Y})
\]
for any $S$-variety with exponential $(X,u)$ and $Y$ closed subscheme of $X$; and
\begin{equation} \label{relation2}
(X \times_{\mathbb{Z}} \mathbb{A}_{\mathbb{Z}}^1, \pr_2)
\end{equation}
for any $S$-variety $X$, where $\pr_2$ is the projection on $\mathbb{A}_{\mathbb{Z}}^1$. The class of $(X,u)$ in $K_0\left(\mathrm{ExpVar}_S\right)$ will be written $[X,u]$, or $[X,u]_S$ if we want to keep track of the base $S$. Such an element is said to be \emph{effective} (as opposed to, for instance, $-[S,0]$).
\begin{remark} \label{rk_not_class_subscript}
For practicality, if $X$ is not itself a $S$-variety but both $X$ and $S$ are $T$-varieties in a way that is clear from the context, we will allow ourselves to still write
\[
[X,u]_S = [X \times_T S,\, u \circ \pr_1] \in K_0\left(\mathrm{ExpVar}_S\right) .
\]
In particular if $X$ and $S$ are not related apart from being over a common field $k$, we write $[X,u]_S$ as above with $T = \Spec k$.
\end{remark}
The ring structure is given by the product
\[
[ X,u] \cdot [Y,v] := \left[ X \times_S Y,\, u\circ\pr_1 + v\circ\pr_2 \right] ;
\]
we have $0 = [\varnothing, 0]$ and $1 = [S, 0]$. We take again the notation $\mathbf{L} = \mathbf{L}_S := \left[\mathbb{A}_S^1, 0\right]$, and invert this class to obtain the ring $\mathscr{E}xp\mathscr{M}_S := K_0\left(\mathrm{ExpVar}_S\right)\left[\mathbf{L}^{-1}\right]$. After lemma 1.1.3 of \cite{CLL16}, $[X] \mapsto [X,0]$ defines injective ring morphisms
\[
K_0\left(\mathrm{Var}_S\right) \longhookrightarrow K_0\left(\mathrm{ExpVar}_S\right) \quad\text{and}\quad \mathscr{M}_S \longhookrightarrow \mathscr{E}xp\mathscr{M}_S .
\]
It should be noted that \cite{Bor18} proved $\mathbf{L}$ to be a zero divisor in characteristic zero. In particular the localisation is not injective in this case.

\begin{remark}[Interpretation as exponential sums]
Let us look at the finite field case for an instant: let $\psi : \mathbb{F}_q \to \C^\times$ be a non trivial additive character. Then we have a morphism of rings
\[
\left\{\begin{array}{rcl}
K_0\big(\mathrm{ExpVar}_{\mathbb{F}_q}\big) & \longrightarrow & \C \\
{[X,u]} & \longmapsto & \sum_{x\in X(\mathbb{F}_q)} \psi(u(x))
\end{array}\right. .
\]
In particular the second type (\ref{relation2}) of relations defining the Grothendieck group is justified because it reflects the fact that $\sum_{a\in \mathbb{F}_q} \psi(a) = 0$. Note that this map extends the map
\begin{equation} \label{count_measure}
\left\{\begin{array}{rcl}
K_0\left(\mathrm{Var}_{\mathbb{F}_q}\right) & \longrightarrow & \mathbb{Z} \\
{[X]} & \longmapsto & \#X(\mathbb{F}_q)
\end{array}\right.
\end{equation}
(which is called the counting measure). For any base $S$, we will think of the elements of $K_0\left(\mathrm{ExpVar}_S\right)$ as a \emph{motivic} version of exponential sums.
\end{remark}

If we have a map $S \to T$ of noetherian schemes, ring morphisms
\begin{equation} \label{GR_functor}
K_0\left(\mathrm{ExpVar}_T\right) \longrightarrow K_0\left(\mathrm{ExpVar}_S\right) \quad\text{and}\quad \mathscr{E}xp\mathscr{M}_T \longrightarrow \mathscr{E}xp\mathscr{M}_S
\end{equation}
are induced by this map. In particular we can write expressions combining elements over $T$ and $S$, with the result living over $S$. For instance if $T$ and $T'$ are noetherian schemes, we can consider expressions involving elements over both of them, the total lying over $T \times T'$. Let us now describe a general \emph{functional interpretation} of our classes. Any $C \in K_0\left(\mathrm{ExpVar}_S\right)$ can be seen as a motivic function with source $S$ by denoting, for any $s \in S$,
\[
C(s) \in K_0\big(\mathrm{ExpVar}_{\kappa(s)}\big)
\]
the image of $C$ via the map $K_0\left(\mathrm{ExpVar}_S\right) \longrightarrow K_0\big(\mathrm{ExpVar}_{\kappa(s)}\big)$ from above (and same with $\mathscr{E}xp\mathscr{M}_\bullet$ instead). The following lemma guarantees that a motivic function is determined by its values.
\begin{lemma}[\cite{CLL16}, lemma 1.1.8] \label{lift_equal}
Let $C \in K_0\left(\mathrm{ExpVar}_S\right)$ (respectively $\mathscr{E}xp\mathscr{M}_S$). If for all $s\in S$, $C(s) = 0$, then $C = 0$.
\end{lemma}

Finally, here is another lemma from the same article, that will be useful to us.
\begin{lemma}[\cite{CLL16}, lemma 1.1.11] \label{lin_form_zero}
Let $V$ be a finite-dimensional $k$-vector space and $\phi : V \to k$ a linear form. Then
\[
[V,\phi] =
\left\{ \begin{array}{ll}
\mathbf{L}^{\dim V} & \text{if } \phi=0 \\
0 & \text{otherwise}
\end{array} \right. .
\]
\end{lemma}
\begin{remark} \label{affine_form_zero}
With the same notations, if $\phi\neq0$ and $c\in k$ is a constant, then
\[
[V, \phi+c] = [V, \phi] \cdot [\Spec k, c] = 0 .
\]
\end{remark}

\subsubsection{Motivic Euler products and zeta functions} \label{MEP}

Motivic Euler products were constructed by Bilu as a formal version of Euler products with sequences of elements of a Grothendieck ring instead of numbers, with the goal of playing their role in a geometric setting, like here. Let us recall their construction from section 3 of \cite{Bil23}. It is based on the configuration spaces (\ref{def_conf}), except we are going to look at elements of $K_0\left(\mathrm{ExpVar}_X\right)$ (respectively $\mathscr{E}xp\mathscr{M}_X$): let us define the configuration space of a family of varieties with exponentials. Suppose that $S$ is a $k$-variety for some field $k$, $X$ is a quasi-projective $S$-variety, and $(X_i, u_i)_{i\geq1}$ is a family of $X$-varieties with exponentials. We already have the
\[
\Conf_{X/S}^\mu \left((X_i)_{i\geq1}\right) = \Bigg( \prod_{i\geq1} X_i^{\mu_i} \Bigg)_{*,X} \ / \mathfrak{S}_\mu 
\]
from (\ref{def_conf}), to which we add the map
\[
u : \left\{\begin{array}{rcl}
\Conf_{X/S}^\mu \left((X_i)_{i\geq1}\right) & \longrightarrow & \mathbb{A}_\mathbb{Z}^1 \\
(x_{i,j})_{i,j} & \longmapsto & \sum_{i,j} u_i(x_{i,j})
\end{array}\right.
\]
to obtain a configuration space with exponential
\begin{equation} \label{def_conf_exp}
\Conf_{X/S}^\mu \left((X_i,u_i)_{i\geq1}\right) := \left( \Conf_{X/S}^\mu \left((X_i)_{i\geq1}\right) , u \right) .
\end{equation}
As for elements of the Grothendieck rings, one can show that if $[X_i,u_i] = [Y_i,v_i]$ for all $i$, then $\left[\Conf_{X/S}^\mu \left((X_i,u_i)_{i\geq1}\right)\right] = \left[\Conf_{X/S}^\mu \left((Y_i,v_i)_{i\geq1}\right)\right]$, which allows us to consider the configuration space of effective classes. Furthermore, in its subsection 3.5, \cite{Bil23} extends the definition to
\[
\Conf_{X/S}^\mu \left((C_i)_{i\geq1}\right)
\]
for any sequence of $C_i \in K_0\left(\mathrm{ExpVar}_X\right)$ (respectively $\mathscr{E}xp\mathscr{M}_X$). Thanks to this notion we can now define motivic Euler products.

\begin{definition}[Motivic Euler products -- \cite{Bil23}, subsection 3.8] \label{def_MEP}
Let $k$ be a field, $S$ a $k$-variety, $X$ a quasi-projective $S$-variety, and $(C_i)_{i\geq1}$ a family of classes in $K_0\left(\mathrm{ExpVar}_X\right)$ (respectively $\mathscr{E}xp\mathscr{M}_X$). We take the notation
\[
\prod_{x \in X/S} \left(1 + \sum_{i\geq1} C_{i,x} T^i\right) := \sum_{j\geq0} \left(\sum_{\mu\vdash j} \Conf_{X/S}^\mu \left((C_i)_{i\geq1}\right)\right) T^j \in K_0\left(\mathrm{ExpVar}_S\right) \llbracket T \rrbracket
\]
(respectively $\in \mathscr{E}xp\mathscr{M}_S \llbracket T \rrbracket$).
\end{definition}
Let us stress that the above is purely a \emph{notation for a formal power series}, and is not an actual product. It however behaves like a product in a number of way as explained in \cite{Bil23}, which justifies the notation. To mention only one of these product properties (which we will use later), associativity holds in the sense that if $Y$ is a closed subscheme of $X$, the product over $X$ is equal to the one over $Y$ times the one over $X \smallsetminus Y$ (see \cite{Bil23}, (3.14)).

Just like symmetric powers in subsection \ref{part_conf}, in the case where $S = \Spec k$ with $k$ clear from the context, we write $\prod_{x \in X}$ in place of $\prod_{x \in X/k}$.
\begin{remark} \label{MEP_short_index}
More generally, if $X = X_0 \times_k S$ with $S$ clear from the context, we will allow ourselves to write $\prod_{x \in X_0}$ for $\prod_{x \in (X_0 \times_k S) / S}$. Keep in mind that this is only a shortening of the proper notation, with the goal of keeping formulas readable while preserving intuition.
\end{remark}

We next define a related object.
\begin{definition}[Kapranov's zeta function -- \cite{Kap00}, (1.3.1)] \label{def_Kapranov}
The Kapranov zeta function of a quasi-projective $k$-variety $X$ is defined as
\[
Z_X(T) := \sum_{j\geq0} \left[\Sym^j(X)\right] T^j \in K_0\left(\mathrm{Var}_k\right) \llbracket T \rrbracket .
\]
\end{definition}
It can easily be written as a motivic Euler product:
\begin{equation} \label{zeta_prod}
Z_X(T) = \prod_{x \in X} \sum_{i\geq0} T^i = \prod_{x \in X} (1-T)^{-1} .
\end{equation}
Remark that over a finite field $k = \mathbb{F}_q$, the counting measure (\ref{count_measure}) sends $Z_X(T)$ to the Weil zeta function $\zeta_X(s)$ with $q^{-s} = T$.

Let us also reproduce a piece of notation and a short computation from \cite{VW15}. In their section 2 they denote by $\mathcal{Q}$ the set of integer partitions $\mu = (\mu_i)_{i\geq1}$ for which there exists an integer $r\geq0$ such that $\mu_i \geq 1$ for $i\leq r$ and $\mu_i = 0$ for $i>r$. Note that $||\mu|| = r$ then. Using this $\mathcal{Q}$, just after their proposition 3.5 they rewrite the inverse of the zeta function as
\begin{equation} \label{inv_zeta_sum}
Z_X^{-1}(T) = \sum_{r\geq0}\left( -\sum_{j\geq1} \left[\Sym^j(X)\right] T^j \right)^r = \sum_{\mu \in \mathcal{Q}} (-1)^{||\mu||} \left[\Sym^\mu(X)\right] T^{|\mu|} .
\end{equation}

\subsubsection{Weight function}
Let $S$ be a complex variety. In order to carry out estimations, we need a topology on the ring $\mathscr{E}xp\mathscr{M}_S$. Like we mentioned in the introduction, one thinks first of the dimensional one: however it is unfit to perform the circle method, and more particularly Weyl differencing (proposition \ref{Weyl_dif}). To resolve this issue, \cite{BB25} constructed in their section 4 a weight function
\[
w_S : \mathscr{E}xp\mathscr{M}_S \to \mathbb{Z} \cup \{-\infty\}
\]
based on Hodge theory, namely mixed Hodge modules. It satisfies
\[
w_S(C_1 \cdot C_2) = w_S(C_1) + w_S(C_2) .
\]
for any $C_1, C_2 \in \mathscr{E}xp\mathscr{M}_S$. We will heavily use this weight function and the associated topology. The completion with respect to it will be denoted $\widehat{\mathscr{E}xp\mathscr{M}}_S$. The only class of weight $-\infty$ is zero, which means that a sequence of classes converges to some limit when the weight of the difference with this limit tends to $-\infty$. Just like we write $\dim$ for $\dim_\C$ we denote simply $w$ for $w_\C$. Let us state some properties that we are going to use.

\begin{proposition}[\cite{BB25}, proposition 4.4] \label{prop_weight}
~
\begin{enumerate}

\item If $C_1, C_2 \in \mathscr{E}xp\mathscr{M}_S$ then
\[
w_S(C_1+C_2) \leq \max(w_S(C_1), w_S(C_2)) .
\]

\item Let $p : T \to S$ be a dominant morphism with fibres of dimension $\leq d$, and $C \in \mathscr{E}xp\mathscr{M}_S$. Then there exists $S_0 \subset S$ open and non-empty such that, denoting by $C'$ the pullback of $C$ to $p^{-1}(S_0)$, we have
\[
w_{S_0}(p_!C') \leq w_{p^{-1}(S_0)}(C') + 2d .
\]

\item Let $i : T \hookrightarrow S$ be a locally closed subvariety: for all $C \in \mathscr{E}xp\mathscr{M}_T$ we have
\[
w_S(i_!C) \leq w_T(C) .
\]

\end{enumerate}
\end{proposition}

\begin{lemma}[\cite{BB25}, lemma 4.6] \label{weight_var}
For $X$ a $S$-variety we have
\[
w_S(X) = 2\dim_S X .
\]
\end{lemma}
In particular $w_S(\mathbf{L}^r) = 2r$ for all $r \in \mathbb{Z}$.

\begin{lemma}[Triangular inequality -- \cite{BB25}, lemma 4.7] \label{triang_ineq}
For any variety with exponential $(X,u)$ over $S$,
\[
w_S([X,u]) \leq w_S(X) = 2\dim_S X .
\]
\end{lemma}

Another important property of this weight is related to irreducible components.
\begin{lemma}[Cancellation of maximal weights -- \cite{BB25}, lemma 4.8] \label{cancel_max}
If $X$ and $Y$ are irreducible complex varieties of same dimension $d$, then
\[
w([X]-[Y]) \leq 2d-1 .
\]
\end{lemma}

We deduce the fact below.
\begin{lemma} \label{nb_compo_max}
Let $X$ be a complex variety. There is a unique integer $\gamma$ such that
\[
w\left( \mathbf{L}^{-\dim X} [X] - \gamma \right) < 0 ,
\]
which is the number of irreducible components of $X$ of maximal dimension.
\end{lemma}

\begin{proof}
Let $d$ be the dimension of $\dim X$ and $\gamma$ its number of irreducible components of dimension $d$. First $[X]$ is the sum of the classes of these components on the one hand, and of other classes of dimension $<d$ (possibly with negative coefficients) on the other hand. After lemma \ref{cancel_max} above, any irreducible component of dimension $d$ writes as $\mathbf{L}^d$ plus some class of weight $<2d$. As a result,
\[
w\left([X] - \gamma \mathbf{L}^d\right) < 2d ,
\]
which is equivalent to the claimed formula. As for the uniqueness, if $\gamma'$ is also an integer satisfying the condition, point 1 of proposition \ref{prop_weight} yields $w(\gamma'-\gamma)<0$, which forces $\gamma' = \gamma$.
\end{proof}

Finally one can study the convergence of power series in $\widehat{\mathscr{E}xp\mathscr{M}}_S$. The \emph{radius of convergence} of a formal power series
\[
F(T) = \sum_{i\geq0} C_i T^i \in \widehat{\mathscr{E}xp\mathscr{M}}_S\llbracket T \rrbracket
\]
is defined as
\begin{equation} \label{def_rad_cv}
\limsup_{i\to\infty} \frac{w_S(C_i)}{2i} .
\end{equation}

\begin{lemma}[\cite{BB25}, subsection 4.4] \label{rad_cv}
If a formal power series $F(T)$ has $r$ for radius of convergence, then for any $N>r$, the evaluation $F\left(\mathbf{L}^{-N}\right) = \sum_{i\geq0} C_i \mathbf{L}^{-Ni}$ exists in $\widehat{\mathscr{E}xp\mathscr{M}}_S$.
\end{lemma}

An important application for us is to motivic Euler products, since they are power series. These products do not behave perfectly with respect to evaluation, so we will a priori write values as
\[
\prod_{x \in X/S} \left(1 + \sum_{i\geq1} C_x T^i\right)_{|T=\mathbf{L}^{-N}} ,
\]
as opposed to directly replacing the $T$ with $\mathbf{L}^{-N}$ in the factors. An exception, already present in other works such as \cite{Fai25_2} and \cite{BB25}, is the notation
\begin{equation} \label{MEP_noT}
\prod_{x \in X/S} C_x := \prod_{x \in X/S} (1 + (C_x-1)T)_{|T=1}
\end{equation}
for any $C \in K_0\left(\mathrm{ExpVar}_X\right)$ (respectively $\mathscr{E}xp\mathscr{M}_X$) such that the right-hand side converges. In particular the products present in theorems \ref{th_naive_intro} and \ref{th_main_intro} are of this sort.

\begin{remark} \label{MEP_noT_sum}
The sum form of such products is
\[
\prod_{x \in X/S} C_x = \sum_{j\geq0} \Conf_{X/S}^{j[1]}(C-1) ,
\]
with $j[1] = [1] + \cdots + [1]$ as defined in subsection \ref{part_conf}. Let us justify this in the case $C = 1 + [Y]$ over $X$. By construction the left-hand side is
\[
\prod_{x \in X/S} (1 + [Y_x]T)_{|T=1} = \sum_{j\geq0} \sum_{\mu\vdash j} \Conf_{X/S}^\mu([Y], 0, 0, \dots) .
\]
These configuration spaces are made up of certain $j$-tuples (up to permutation) of elements of the varieties in argument: here they all have to live in $Y$ so the configuration spaces for partitions not of the form $j[1]$ are empty. Then the $0$ classes do not play a role anymore in the argument, so they can be replaced with $[Y]$ for instance, which results in $\Conf_{X/S}^{j[1]}([Y])$.
\end{remark}

\newpage
\section{Circle method}
\label{sect_CM}

In this section we run the circle method for the naive spaces. During the first three subsections we introduce some more notation, start the machinery, and partition a certain space of linear forms into \emph{major and minor arcs} (the former will give rise to the main term of the estimation while the latter will end up in the error term). Then in subsection \ref{subsect_Weyl} we prove two versions of Weyl's inequality: a `global' one concerned with sections of $P_L$, following \cite{HL24}'s strategy with slope theory; and a `local' one with sections of $P_{|Z}$, where we use the more classical approach through geometry of numbers.

Next we estimate the contribution of the major arcs in the long subsection \ref{subsect_maj}. We first write it as a finite sum of certain classes of varieties with exponentials (proposition \ref{maj_1}), and show that this sum can be replaced, up to an error term (thanks to the local Weyl inequality), with the corresponding converging series, called the \emph{singular series}. In paragraph \ref{subsubsect_sing_MEP}, we rewrite the singular series as a motivic Euler product, leading to corollary \ref{sing_mep_2}. The two following paragraphs consist in reexpressing the factors respectively as local densities as we mentioned in (\ref{intro_mep_local_dens}), and as factors closer to the ones appearing in our main theorem (\ref{equidistrib_C}).

Finally in subsection \ref{subsect_min} we tackle the minor arcs, applying the global Weyl inequality, before concluding with subsection \ref{subsect_concl_naive}.

\subsection{Framework}
\label{CM_framework}

During the whole section we suppose that our hypersurface $V \subset \projsp_\C^n$ has degree $d\geq2$, and that $n \geq 2^d(d-1)$. Fix an integer $e \geq 2g$, which ensures that all $L\in\Pic^e(\scrC)$ are basepoint-free by (\ref{cond_BPF}). The goal is to prove theorem \ref{th_naive_intro}, or rather a generalisation of it where we impose some vanishing conditions on the sections in the naive space. Fix an integer partition $\mu$ as well and recall that in (\ref{def_ST}) we took the notation
\[
\scrST = \scrS+|\scrT|
\]
for any $\scrT\in\Sym^\mu(\CmS)$. All these new conditions divisors have degree
\[
l' := l + |\mu| .
\]
We put
\[
0_{|\scrT|} := \{0\} \subset P_{|\,|\scrT|}^{n+1} \big/ \mathbb{G}_{m, P_{|\, |\scrT|}} :
\]
this is a set of conditions just like $W \subset P_{|\scrS}^{n+1} \big/ \mathbb{G}_{m, P_{|\scrS}}$, except it is not in $\Hom(|\scrT|,V)$. Combining the two we can consider total conditions
\[
W \times 0_{|\scrT|} \subset P_{|\scrST}^{n+1} \big/ \mathbb{G}_{m, P_{|\scrST}} ,
\]
the lift
\[
\widetilde{W \times 0_{|\scrT|}} = \tilde W \times \{0\} \subset P_{|\scrST}^{n+1}
\]
and the naive spaces
\[
M_L^{W \times 0_{|\scrT|}} := \left\{ \mathbf{u} \in P_L^{n+1} \,\left|\,
\begin{array}{l}
f(\mathbf{u}) = 0 \\
\mathbf{u}|_{\scrST} \in \tilde W \times \{0\}
\end{array} \right.\right\}
\]
that generalise the $M_L^W$ from (\ref{def_ML}), with $L \in\Pic^e(\scrC)$.

We make use of the vector bundles we constructed in paragraph \ref{Poincare}, more precisely the ones over $\Pic^e(\scrC) \times \Sym^\mu(\CmS)$: It is going to be useful to take the short notation
\[
X \,\diamond\, Y:= X \times_{\Pic^e(\scrC) \times \Sym^\mu(\CmS)} Y
\]
for the fibre product of two varieties over that base. Like in remark \ref{constr_fibr_Hom}, the points of $\calP_e^{\diamond(n+1)}$ (for instance) are the $(L, \scrT, \mathbf{u})$ with $\mathbf{u} \in P_L^{n+1}$. Now, for reasons that will become clear in section \ref{sect_Hom}, the object we will be studying in this section is
\begin{equation} \label{def_calM}
\calM_e^{W,\mu} := \left\{ (L, \scrT, \mathbf{u}) \in \calP_e^{\diamond(n+1)} \,\left|\,
\begin{array}{l}
f(\mathbf{u}) = 0 \\
\mathbf{u}|_{\scrST} \in \tilde W \times \{0\}
\end{array} \right.\right\} .
\end{equation}
Naturally its fibre over a $(L,\scrT) \in \Pic^e(\scrC) \times \Sym^\mu(\CmS)$ is precisely $M_L^{W \times 0_{|\scrT|}}$. If $\mu=\varnothing$ then $\Sym^\varnothing(\CmS) = \{\varnothing\}$ and $\calM_e^{W,\varnothing}$ is the $\calM_e^W$ we announced in (\ref{fibr_Me}). The main base over which we will view $\calM_e^{W,\mu}$ is
\begin{equation} \label{def_Bmu}
\calB_\mu := J(\scrC) \times \Sym^\mu(\CmS) \times \tilde W
\end{equation}
(recall that $\Pic^e(\scrC) \simeq J(\scrC)$).

To the condition $e \geq 2g$ we add the second requirement
\begin{equation} \label{cond_CM}
e \geq 2g-1+l' ,
\end{equation}
which is the condition needed to apply the Euclidean division exact sequence \ref{exseq_fam_Smu}:
\[
\xymatrix{
0 \ar[r] & \calP_{e(-\scrS_\mu)} \ar[r] & \calP_e \ar[r] & \calP_{|\scrS_\mu} \ar[r] &0
} .
\]
Let
\[
\sigma : \calP_{|\scrS_\mu} \longrightarrow \calP_e
\]
be a section of this short exact sequence; for a given $(L,\scrT)$ we denote its fibre by
\begin{equation} \label{def_sect_sigma_1}
\sigma_L^\scrT : P_{|\scrST} \to P_L .
\end{equation}
It provides the isomorphism
\[
\left\{\begin{array}{rcl}
\calP_{e(-\scrS_\mu)} \,\diamond\, \calP_{|\scrS_\mu} & \overset{\sim}{\longrightarrow} & \calP_e \\
(L, \scrT, v, w') & \longmapsto & (L, \scrT, v + \sigma_L^\scrT(w'))
\end{array}\right.
\]
(here $v \in P_{L(-\scrST)}$, $w' \in P_{|\scrST}$ and $v + \sigma_L^\scrT(w') \in P_L$).
We are interested in $\mathbf{w} \in \tilde W \subset P_{|\scrS}^{n+1} \hookrightarrow P_{|\scrST}^{n+1}$: we extend $\sigma$ to $\calP_{|\scrS_\mu}^{\diamond(n+1)} \to \calP_e^{\diamond(n+1)}$, then write $\sigma_L^\scrT(\mathbf{w}) := \sigma_L^\scrT(\mathbf{w},0)$. Via this isomorphism our first naive spaces become
\begin{equation} \label{reformulation_M}
\calM_e^{W,\mu} \simeq \left\{\left. (L, \scrT, \mathbf{v}, \mathbf{w}) \in \calP_{e(-\scrS_\mu)}^{\diamond(n+1)} \times \tilde W \,\right|\, f\left(\mathbf{v} + \sigma_L^\scrT(\mathbf{w})\right)=0 \right\} .
\end{equation}

\begin{remark} \label{expect_dim_naive}
We can give a heuristic for the expected dimension of the above, as announced in (\ref{def_xi}). Its fibres are
\[
\left(\calM_e^{W,\mu}\right)_{L,\scrT,w} \simeq \left\{\left. \mathbf{v} \in P_{L(-\scrST)}^{n+1} \,\right|\, f\left(\mathbf{v} + \sigma_L^\scrT(\mathbf{w})\right)=0 \right\} ;
\]
write $\kappa = \kappa(L,\scrT,\mathbf{w})$. By Riemann-Roch (\ref{RR}), under the right conditions,
\[
\dim_\kappa P_{L(-\scrST)}^{n+1} = (n+1)(e-l-|\mu|-g+1) .
\]
Furthermore $\mathbf{v} + \sigma_L^\scrT(\mathbf{w}) \in P_{L(-|\scrT|)}^{n+1}$ so $f\left(\mathbf{v} + \sigma_L^\scrT(\mathbf{w})\right) \in P_{L(-|\scrT|)^{\otimes d}}$, and the latter has $\scrS$-restriction $f(\mathbf{w})$, which fixes $l$ of its coefficients already; as a result, it lives in an affine space of dimension
\[
\dim P_{L(-|\scrT|)^{\otimes d}} - l = d(e-|\mu|)-g+1-l .
\]
The expected $\kappa$-dimension of the fibre is then the difference of these two expressions, that is
\begin{equation} \label{def_xi_mu}
\xi_{e,|\mu|} := (n-d+1)(e-|\mu|) - n(g-1) - nl .
\end{equation}
This is also the expected $\calB_\mu$-dimension of $\calM_e^{W,\mu}$. When $\mu=\varnothing$ we retrieve the $\xi_e$ from (\ref{def_xi}).
\end{remark}

\subsection{Initialisation}

As seen earlier, for any $L \in \Pic^e(\scrC)$, $f$ defines a map $P_L \to P_{L^{\otimes d}}$. To put this second space over the appropriate base, we pull the vector buncle $\calP_{de}$ back through
\begin{equation} \label{def_taud}
\tau_d : \left\{\begin{array}{rcl}
\Pic^e(\scrC) \times \Sym^\mu(\CmS) & \longrightarrow & \Pic^{de}(\scrC) \times \Sym^\mu(\CmS) \\
(L, \scrT) & \longmapsto & (L^{\otimes d}, \scrT)
\end{array}\right. ,
\end{equation}
obtaining $\tau_d^* \calP_{de}$ over $\Pic^e(\scrC) \times \Sym^\mu(\CmS)$, with fibres $(\tau_d^* \calP_{de})_{L,\scrT} = (\calP_{de})_{L^{\otimes d}, \scrT} = P_{L^{\otimes d}}$. Then we consider the dual
\[
\tau_d^* \calP_{de}^\vee ,
\]
whose fibres are the spaces $P_{L^{\otimes d}}^\vee$ of linear forms over $P_{L^{\otimes d}}$. First for a fixed $L$ we give the following version of the motivic orthogonality lemma, which is analogous to the classical integral identity (\ref{orth_int}).

\begin{lemma}[orthogonality] \label{lem_orth}
For any fixed $u\in P_{L^{\otimes d}}$ we have, over $\kappa(u)$,
\[
\left[ P_{L^{\otimes d}}^\vee ,\, \alpha(u) \right]_{\kappa(u)} =
\left\{ \begin{array}{ll}
\mathbf{L}^{de-g+1} & \text{if } u=0 \\
0 & \text{otherwise}
\end{array} \right. .
\]
\end{lemma}

\begin{proof}
It is a direct application of lemma \ref{lin_form_zero}, using that $\alpha \mapsto \alpha(u)$ is a linear form on $P_{L^{\otimes d}}^\vee$, and is zero if and only if $u=0$. We also use that under our assumptions, Riemann-Roch (\ref{RR}) provides $\dim P_{L^{\otimes d}} = de-g+1$.
\end{proof}

We deduce the result below, starting point of the circle method.

\begin{proposition} \label{init}
In the Grothendieck ring $\mathscr{E}xp\mathscr{M}_{\calP_{e(-\scrS_\mu)}^{\diamond(n+1)} \times \tilde W}$ we have the relation
\[
\left[ \calM_e^{W,\mu} \right] = \mathbf{L}^{-de+g-1} \left[ \tau_d^* \calP_{de}^\vee \,\diamond\, \calP_{e(-\scrS_\mu)}^{\diamond(n+1)} \times \tilde W ,\, \alpha\left(f\left( \mathbf{v} + \sigma_L^\scrT(\mathbf{w}) \right)\right) \right] .
\]
\end{proposition}
On the right hand side, $\alpha \in P_{L^{\otimes d}}^\vee$, $\mathbf{v} \in P_{L(-\scrST)}^{n+1}$ and $\mathbf{w} \in \tilde W$.

\begin{proof}
Fix $(L, \scrT, \mathbf{v}, \mathbf{w}) \in \calP_{e(-\scrS_\mu)}^{\diamond(n+1)} \times \tilde W$, and consider the associated residue field $\kappa = \kappa(L,\scrT,\mathbf{v},\mathbf{w})$. We evaluate both sides of the claimed equality at $(L,\scrT,\mathbf{v},\mathbf{w})$, in the motivic sense. First the right-hand side simply becomes
\[
\mathbf{L}^{-de+g-1} \left[ P_{L^{\otimes d}}^\vee ,\, \alpha\left(f\left( \mathbf{v} + \sigma_L^\scrT(\mathbf{w}) \right)\right) \right]_\kappa ,
\]
which by the previous lemma is $1$ if $f\left( \mathbf{v} + \sigma_L^\scrT(\mathbf{w}) \right) = 0$ and $0$ otherwise. As for the left-hand side, upon evaluation it becomes $\{(L, \scrT, \mathbf{v}, \mathbf{w})\} = 1$ if this point belongs to $\calM_e^{W,\mu}$, and $0$ otherwise. By description (\ref{reformulation_M}) the vanishing of $f$ caracterises $\calM_e^{W,\mu}$, so in all cases the two values coincide over $\kappa$. By lemma \ref{lift_equal}, the equality is therefore verified over our base.
\end{proof}

\begin{remark} \label{why_Pvee}
The use we make of lemma \ref{lem_orth} can be thought of as a motivic integration over $P_{L^{\otimes d}}^\vee$ (in simple cases like this one, this amounts to considering the concerned class over the right base and renormalising it). Let us give an insight on why we make this dual space appear.\footnote{We mention several articles but try to keep one coherent set of notations.}

In genus zero and over finite fields, such as in \cite{Lee11} (section 2), we are starting with $f(\mathbf{u}) \in k[t]_{\leq de}$: let $k((t^{-1}))$ be the space of formal Laurent series in $t^{-1}$, $\mathbb{T} := t^{-1} k\llbracket t^{-1}\rrbracket$ be the subspace of those with no polynomial part, and $\res$ be the residue of a given Laurent series, i.e. its coefficient of degree $-1$. In this context one integrates an exponential of $\res(\alpha f(\mathbf{u}))$ over $\alpha \in \mathbb{T}$, with a result analogous to the above proposition. In the motivic version of \cite{BB25} (lemma 3.6), with $k=\C$, it is essentially the same but motivic integration requires a finite-dimensional domain of integration, so $\mathbb{T}$ is replaced with a subspace where only the coefficients of relevant degree of $\alpha$ are considered.

When $g\geq1$, $f(\mathbf{u})$ is no longer a polynomial but lives in $P_{L^{\otimes d}}$. To replace the previous residue, note that the $\res(\alpha\cdot)$ were linear forms over $k[t]_{\leq de}$ (see the proof of lemma \ref{N'_poly} for an avatar of this reasoning): one can directly replace the space of the $\alpha$ with $k[t]_{\leq de}^\vee$, or more generally $P_{L^{\otimes d}}^\vee$. This is what \cite{HL24} does over finite fields (above lemma 16; although the integral is simply a finite sum). For a more thorough investigation of this correspondance, see his section 6. Coming back to us, we use that same $P_{L^{\otimes d}}^\vee$, but over $\C$ and with motivic integration.
\end{remark}

\subsection{Factorisation of linear forms and partition into major and minor arcs}
\label{def_arcs}

In order to carry out the circle method, we need to cut $\tau_d^* \calP_{de}^\vee$ into two: one part, called the \emph{major arcs}, on which we will refine and transform the previous expression, and another, called the \emph{minor arcs}, whose contribution we will bounded. To continue on remark \ref{why_Pvee}, the major arcs are classically defined as the $\alpha \in \mathbb{T}$ that are well-approximated by a rational fraction with small denominator. Let $L\in\Pic^e(\scrC)$: for linear forms in $P_{L^{\otimes d}}^\vee$, this potential denominator (or rather its roots) is analogous to an effective divisor $Z$ of $\scrC$, and the good approximation amounts to the linear form factoring through $Z$, as explained below.

\begin{definition}
We say an $\alpha \in P_{L^{\otimes d}}^\vee$ factors through $Z$ if it factors through $P_{|Z}$, that is if there exists $\beta \in P_{|Z}^\vee$ such that for all $u \in P_{L^{\otimes d}}$, $\alpha(u) = \beta(u|_Z)$. We also say the factorisation is minimal if $\alpha$ factors through no strict subscheme of $Z$. We respectively write $\alpha \leadsto Z$ and $\alpha \underset{\min}{\leadsto} Z$.
\end{definition}
Consider the dual exact sequence of (\ref{exseq_H1}) applied to $L^{\otimes d}$ and $Z$:
\begin{equation} \label{exseq_dual}
\xymatrix{
H^1(\scrC, L^{\otimes d}(-Z))^\vee \ar[r] & P_{|Z}^\vee \ar[r] & P_{L^{\otimes d}}^\vee \ar[r] & P_{L^{\otimes d}(-Z)}^\vee \ar[r] & 0
} ,
\end{equation}
where the last non-trivial arrow is simply the restriction to $P_{L^{\otimes d}(-Z)}$. Note that $\alpha \leadsto Z$ is equivalent to $\alpha \in \mathrm{Im} \left[ P_{|Z}^\vee \to P_{L^{\otimes d}}^\vee \right]$, itself equivalent to $\alpha|_{P_{L^{\otimes d}(-Z)}} = 0$. Moreover under the condition $\deg Z \leq de-2g+1$ the $H^1$ above disappears, providing an isomorphism
\[
P_{|Z}^\vee \overset{\sim}{\longrightarrow} \left\{\left. \alpha \in P_{L^{\otimes d}}^\vee \,\right|\, \alpha \leadsto Z \right\} .
\]
We define $P_{|Z}^{\vee*} \subset P_{|Z}^\vee$ as the subvariety of the $\beta$ that factor through no strict subscheme of $Z$: the above restricts to
\begin{equation} \label{Pvee*}
P_{|Z}^{\vee*} \overset{\sim}{\longrightarrow} \left\{ \alpha \in P_{L^{\otimes d}}^\vee \,\left|\, \alpha \underset{\min}{\leadsto} Z \right.\right\} .
\end{equation}
Now let us state two results taken from \cite{HL24}. Again in this article the computations take place over a finite field, but these two proofs do not depend on the field.
\begin{lemma}[\cite{HL24}, lemma 12] \label{max_deg}
Suppose $e-2g+1 \geq0$. Then any $\alpha \in P_{L^{\otimes d}}^\vee$ factors through a $Z$ of degree at most $\frac{d}{2}e+1$.
\end{lemma}

Under this hypothesis $e-2g+1 \geq0$, for $\alpha \in P_{L^{\otimes d}}^\vee$ we define
\[
\deg\alpha := \min\{ \deg Z \,;\, \alpha \leadsto Z \}
\]
and therefore have $0 \leq \deg\alpha \leq \frac{d}{2}e+1$.

\begin{lemma}[\cite{HL24}, remark 14] \label{factor_min_unique}
If $\deg\alpha \leq e - 2g + 1$, there exists a unique effective divisor $Z$ such that $\alpha \underset{\min}{\leadsto} Z$.
\end{lemma}
Of course in that case, $\deg Z = \deg\alpha \leq e - 2g + 1$.

Now is the time to define the major arcs. We want to work with $Z$ satisfying several properties, therefore imposing conditions on the degree:
\begin{itemize}
\item $\deg Z \leq de-2g+1$ for isomorphism (\ref{Pvee*});
\item $\deg Z \leq e-2g+1$ for the above lemma \ref{factor_min_unique};
\item $l' + \deg Z \leq e-2g+1$ for Euclidean division vanishing over $\scrST$ as in lemma \ref{exseq_ed_vanish}.
\end{itemize}
The latter is the strongest, which leads us to set\footnote{We don't make the dependency on $l'$ explicit in the notation because it will only appear in the current context of the circle method, where $l' = l+|\mu|$ is fixed.}
\begin{equation} \label{def_eta}
\eta_e := e - l' - 2g + 1
\end{equation}
(which is $\geq0$ under assumption (\ref{cond_CM})) and define the major arcs as
\[
\mathfrak{M} := \left\{\left. (L, \scrT, \alpha) \in \tau_d^* \calP_{de}^\vee \,\right|\, \deg\alpha \leq \eta_e \right\} ;
\]
naturally the minor arcs will be $\mathfrak{m} := \tau_d^* \calP_{de}^\vee \smallsetminus \mathfrak{M}$. Going back to proposition \ref{init}, the respective contributions will be written
\[
\Gamma_\mathfrak{M} := \mathbf{L}^{-de+g-1} \left[ \mathfrak{M} \,\diamond\, \calP_{e(-\scrS_\mu)}^{\diamond(n+1)} \times \tilde W ,\, \alpha\left(f\left( \mathbf{v} + \sigma_L^\scrT(\mathbf{w}) \right)\right) \right]
\]
and the same for $\Gamma_\mathfrak{m}$ with $\mathfrak{m}$ in place of $\mathfrak{M}$, such that
\begin{equation} \label{part_maj_min}
\left[ \calM_e^{W,\mu} \right] = \Gamma_\mathfrak{M} + \Gamma_\mathfrak{m} .
\end{equation}

\subsection{Weyl's inequality}
\label{subsect_Weyl}

Before going on the computations with the $\Gamma_\mathfrak{M}$ and $\Gamma_\mathfrak{m}$ we just defined, in this subsection we make a detour to establish the main tool we will use for bounding the error terms that will arise: Weyl's inequality (or rather, two separate avatars of it).

\subsubsection{Weyl differencing}
\label{Weyl_gen}

If $\mathbf{x} = (x_0, \dots, x_n)$ is a $(n+1)$-tuple of indeterminates, for any ring $R$ we will denote the corresponding polynomial ring $R[\mathbf{x}] = R[x_0, \dots, x_n]$.

\begin{definition}[\cite{Sch82}, section 2] \label{def_Gd}
Let $G : A \to A'$ be a map between two abelian groups (we do not request $G$ to be a group morphism): we define the map
\[
G^{(d)} :
\left\{\begin{array}{rcl}
A^d & \longrightarrow & A' \\
(y^{(1)}, \dots, y^{(d)}) & \longmapsto & \sum_{\varepsilon_1, \dots, \varepsilon_d \in \{0,1\}} (-1)^{\varepsilon_1 + \cdots + \varepsilon_d} G\left( \varepsilon_1 y^{(1)} + \cdots + \varepsilon_d y^{(d)} \right)
\end{array}\right. .
\]
In particular if $G \in R[\mathbf{y}]$ is a polynomial over a ring $R$, the same formula defines a polynomial $G^{(d)} \in R\!\left[ \mathbf{y}^{(1)}, \dots, \mathbf{y}^{(d)} \right]$.
\end{definition}

Let $\calV$ be a vector bundle over a base $B$, endowed with a map $G : \calV \to \mathbb{A}_B^1$: it induces maps $G_b : \calV_b \to \mathbb{A}_{\kappa(b)}^1$. Write $\calV^{d-1}$ for the $(d-1)$-th power of $\calV$ over $B$. We denote by
\begin{equation} \label{def_VG}
\calV^{d-1}(G) \longhookrightarrow \calV^{d-1}
\end{equation}
the subvariety made up of the $\left(b, \mathbf{y}^{(1)}, \dots, \mathbf{y}^{(d-1)}\right)$ for which
\[
\mathbf{x} \longmapsto G_b^{(d)}\left( \mathbf{y}^{(1)}, \dots, \mathbf{y}^{(d-1)}, \mathbf{x} \right)
\]
is constant over $\calV_b$. This is a generalisation of the `$V(G)$' appearing in the proposition 5.3 of \cite{BB25}, with $\calV$ any vector bundle instead of the trivial one $\mathbb{A}_B^r$.
\begin{remark}
In \cite{BB25}, the considered function is actually
\[
\sum_{\varepsilon_1, \dots, \varepsilon_{d-1} \in \{0,1\}} (-1)^{\varepsilon_1 + \cdots + \varepsilon_{d-1}} G_b\left( \varepsilon_1 \mathbf{y}^{(1)} + \cdots + \varepsilon_{d-1} \mathbf{y}^{(d-1)} + \mathbf{x} \right) .
\]
The difference is the terms for which the $\varepsilon_d$ in front of $\mathbf{x}$ is zero. Since these terms do not depend on $\mathbf{x}$, the resulting variety is identical.
\end{remark}

Suppose that $B$ is a complex variety. The two following results are generalisations of lemma 5.2 and proposition 5.3 from \cite{BB25}.

\begin{lemma} \label{key_Weyl}
If $(X,u)$ is a $B$-variety with exponential, then
\[
w_B([X,u] \cdot [X,-u]) = 2w_B([X,u]) .
\]
\end{lemma}
The proof from \cite{BB25} holds identically for any variety with exponential. In fact the first step is to link $[X,u]$ to a linear combination of simple pure polarisable variations of Hodge structures, which requires no hypothesis on $[X,u]$. This combination is then put together with the analogous one for $[X,-u]$, resulting in something twice the weight of $[X,u]$.

\begin{proposition}[Weyl differencing] \label{Weyl_dif}
Suppose that $\calV \to B$ is a vector bundle of rank $r$, coming with $G : \calV \to \mathbb{A}_B^1$ such that for all $b\in B$, $G_b : \calV_b \to \mathbb{A}_{\kappa(b)}^1$ is polynomial of degree $\leq d$ in the $r$ coefficients. Then the weight function $w_B$ satisfies
\[
w_B([\calV,G]) \leq 2^{-(d-2)} \left( \dim_B \calV^{d-1}(G) + \left(2^{d-1}-d+1\right)r \right) .
\]
\end{proposition}

\begin{proof}
The proof is the same as in \cite{BB25}, which was itself adapted from classical arguments as found in \cite{BS20} (proposition 5.5). As opposed to the previous proof, it relies on the fact we are working with vector bundles, so this time we fully write it. It is an induction on $d\geq1$.

Assume first that $d=1$. Then $\calV^0(G)$ is the subvariety of $B$ given by the $b$ for which $G_b^{(1)}(\mathbf{x}) = G_b(0) - G_b(\mathbf{x})$ is constant in $\mathbf{x}$, or equivalently $G_b(\mathbf{x})$ is. In the case $\calV^0(G) = \varnothing$, for all $b$ the polynomial $G_b$ is of degree $1$, so $[\calV_b, G_b] = 0$ by remark \ref{affine_form_zero}. By lemma \ref{lift_equal} this entails $[\calV, G] = 0$ over $B$, and both sides of the desired inequality are equal to $-\infty$. Now suppose $\calV^0(G) \neq \varnothing$: by the triangular inequality \ref{triang_ineq}, $w_B([\calV,G]) \leq w_B(\calV) = 2r$, which is again the expected inequality.

Now let $d\geq2$ be fixed, and assume the result to be true for $d' \leq d-1$. The key is starting with lemma \ref{key_Weyl}:
\begin{align} \label{Weyl_dif_proof_key}
2 w_B([\calV,G]) & = w_B([\calV,G] \cdot [\calV,-G]) \nonumber \\
& = w_B([\calV \times_B \calV, G(\mathbf{x})-G(\mathbf{y})]) \nonumber \\
& = w_B([\calV \times_B \calV, G(\mathbf{x})-G(\mathbf{x}+\mathbf{y})])
\end{align}
via a change of variables. For the rest of the proof we write
\[
H_{\mathbf{y}}(\mathbf{x}) := G(\mathbf{x})-G(\mathbf{x}+\mathbf{y}) :
\]
for a fixed $\mathbf{y}$, this has degree $\leq d-1$ in the $r$ coefficients of $\mathbf{x}$. Let $p : \calV^{d-1}(G) \to \calV$ be the projection on the first coordinate. Its fibre over $\mathbf{y}^{(1)}$ is $\{\mathbf{y}^{(1)}\} \times \calV^{d-2}(H_{\mathbf{y}^{(1)}})$, which follows from the definition (\ref{def_VG}) of these varieties and from
\begin{align*}
G_b^{(d)} \left(\mathbf{y}^{(1)}, \dots, \mathbf{y}^{(d-1)}, \mathbf{x}\right) & = \sum_{\varepsilon_2, \dots, \varepsilon_d \in \{0,1\}} (-1)^{\varepsilon_2 + \cdots + \varepsilon_d} G_b\left({\textstyle \sum_{i=2}^{d-1} \varepsilon_i \mathbf{y}^{(i)} + \varepsilon_d \mathbf{x} }\right) \\
& \qquad - \sum_{\varepsilon_2, \dots, \varepsilon_d \in \{0,1\}} (-1)^{\varepsilon_2 + \cdots + \varepsilon_d} G_b\left({\textstyle \mathbf{y}^{(1)} + \sum_{i=2}^{d-1} \varepsilon_i \mathbf{y}^{(i)} + \varepsilon_d \mathbf{x} }\right) \\
& =  \sum_{\varepsilon_2, \dots, \varepsilon_d \in \{0,1\}} (-1)^{\varepsilon_2 + \cdots + \varepsilon_d} H_{\mathbf{y}^{(1)}}\left({\textstyle \sum_{i=2}^{d-1} \varepsilon_i \mathbf{y}^{(i)} + \varepsilon_d \mathbf{x} }\right) \\
& = H_{\mathbf{y}^{(1)}}^{(d-1)} \left(\mathbf{y}^{(2)}, \dots, \mathbf{y}^{(d-1)}, \mathbf{x}\right) .
\end{align*}
Let $(\calW_j)_j$ be a finite stratification of $\calV$ such that the dimension of $\calV^{d-2}(H_{\mathbf{y}^{(1)}})$ is constant over each $\calW_j$. Denote $q_j : \calW_j \to B$ the structure maps, and $B_j := q_j(\calW_j)$: up to refining the stratification we can assume that the $q_j : \calW_j \to B_j$ have fibres of constant dimension, and that the $B_j$ are locally closed in $B$. By proposition \ref{prop_weight} (first the third statement with $B_j \hookrightarrow B$, then the second one with the dominant map $q_j$) we have
\begin{align} \label{Weyl_dif_proof_appliprop}
w_B([\calV \times_B \calW_j, H]) & \leq w_{B_j}([\calV \times_{B_j} \calW_j, H]) \nonumber \\
& \leq 2\dim_{B_j} \calW_j + w_{\calW_j}([\calV \times_{B_j} \calW_j, H]) \nonumber \\
& = 2\dim_B \calW_j + w_{\calW_j}([\calV \times_{B_j} \calW_j, H]) .
\end{align}
Now we view $\calV \times_{B_j} \calW_j$ as a vector bundle over $\calW_j$: it has rank $r$ and comes with the map $H$ to $\mathbb{A}_{\calW_j}^1$ such that $H_{\mathbf{y}}$ is polynomial of degree $\leq d-1$ for all $\mathbf{y}$, so the induction hypothesis yields
\begin{equation} \label{Weyl_dif_proof_indhyp}
w_{\calW_j}([\calV \times_{B_j} \calW_j, H]) \leq 2^{-(d-3)} \left( \dim_{\calW_j} \calV^{d-2}(H) + \left(2^{d-2}-d+2\right)r \right) .
\end{equation}
Note that $\calV^{d-2}(H) = p^{-1}(\calW_j) \subset \calV^{d-1}(G)$, so that
\begin{align*}
\dim_{\calW_j} \calV^{d-2}(H) & = \dim_B \calV^{d-2}(H) - \dim_B \calW_j \\
& \leq \dim_B \calV^{d-1}(G) - \dim_B \calW_j .
\end{align*}
We plug this into (\ref{Weyl_dif_proof_indhyp}), that we in turn plug into (\ref{Weyl_dif_proof_appliprop}) to obtain
\[
w_B([\calV \times_B \calW_j, H]) \leq 2^{-(d-3)} \left( \dim_B \calV^{d-1}(G) + (2^{d-2}-1)\dim_B \calW_j + \left(2^{d-2}-d+2\right)r \right) .
\]
Since $\dim_B \calW_j \leq \dim_B \calV = r$ the bound is smaller or equal to
\[
2^{-(d-3)} \left( \dim_B \calV^{d-1}(G) + \left(2^{d-1}-d+1\right)r \right) .
\]
Going back to (\ref{Weyl_dif_proof_key}) and using the first statement of proposition \ref{prop_weight},
\[
2 w_B([\calV,G]) = w_B([\calV \times_B \calV, H]) \leq \max_j w_B([\calV \times_B \calW_j, H]):
\]
by replacing the right-hand side with the bound we found we get the desired inequality.
\end{proof}

The above will be related to what follows. Our homogeneous polynomial $f \in \C[\mathbf{x}]$ can be written
\[
f(\mathbf{x}) = \sum_{j_1, \dots, j_d = 0}^n c_\mathbf{j} \ x_{j_1} \cdots x_{j_d}
\]
with symmetric coefficients (i.e. $c_\mathbf{j} = c_{\sigma \cdot \mathbf{j}}$ for all $\sigma \in \mathfrak{S}_d$). The multilinear forms associated to $f$ are the
\begin{equation} \label{def_Psi}
\Psi_j\left( \mathbf{y}^{(1)}, \dots, \mathbf{y}^{(d-1)} \right) := d! \sum_{j_1, \dots, j_{d-1} = 0}^n c_{j_1, \dots, j_{d-1}, j} \ y_{j_1}^{(1)} \cdots y_{j_{d-1}}^{(d-1)}
\end{equation}
for $0\leq j\leq n$.

\begin{lemma} \label{link_multilin}
Let $\mathbf{a} = (a_0, \dots, a_n)$ be an additional tuple of indeterminates. We can apply definition \ref{def_Gd} to $f(\mathbf{x}+\mathbf{a})$ where $\mathbf{x}$ is treated as the entry and $\mathbf{a}$ as constants, and consider the result $f(\cdot+\mathbf{a})^{(d)}$. Then in $\C[\mathbf{a}]\!\left[ \mathbf{y}^{(1)}, \dots, \mathbf{y}^{(d-1)}, \mathbf{x} \right]$ the following identity holds:
\[
f(\cdot+\mathbf{a})^{(d)}\left( \mathbf{y}^{(1)}, \dots, \mathbf{y}^{(d-1)}, \mathbf{x} \right) = \sum_{j=0}^n \Psi_j\left( \mathbf{y}^{(1)}, \dots, \mathbf{y}^{(d-1)} \right) x_j .
\]
In particular $\mathbf{a}$ does not actually appear in $f(\cdot+\mathbf{a})^{(d)}$.
\end{lemma}
\begin{proof}
$f(\mathbf{x}+\mathbf{a})$ is a polynomial with coefficients in $\C[\mathbf{a}]$, of degree $d$ (in $\mathbf{x}$) and with homogeneous part $f$: it is classical\footnote{See section 11 of \cite{Sch85} for related arguments.} that both sides of the desired identity are equal to $f^{(d)}\left( \mathbf{y}^{(1)}, \dots, \mathbf{y}^{(d-1)}, \mathbf{x} \right)$.
\end{proof}

With these different ingredients we can properly work towards the two versions of Weyl's inequality we need.

\subsubsection{Global inequality}
\label{Weyl_shrink_global}

First we focus on what we call the `global' version, that is the one that involves global sections in $P_{L(-\scrST)}$, inspired by \cite{HL24}. The following object is standard in analytical number theory and in particular appears in the proofs of the various iterations of Weyl's inequality (usually as a cardinality, but in our motivic context we consider the variety itself instead).

\begin{definition} \label{def_N_global}
For $(L,\scrT,\alpha) \in \tau_d^* \calP_{de}^\vee$, where $\tau_d$ was defined in (\ref{def_taud}) (i.e. $\alpha \in P_{L^{\otimes d}}^\vee \times_{\kappa(L)} \kappa(L,\scrT)$), we define the variety
\[
N(L,\scrT;\alpha) := \left\{\left. \underline{\mathbf{v}} \in \left(P_{L(-\scrST)}^{n+1}\right)^{d-1} \,\right|\, \forall j, \forall u \in P_{L(-\scrST)},\,  \alpha\left( \Psi_j(\underline{\mathbf{v}}) \cdot u \right) = 0 \right\} ,
\]
where we write $\underline{\mathbf{v}} = \left(\mathbf{v}^{(1)}, \dots, \mathbf{v}^{(d-1)}\right) \in \left(P_{L(-\scrST)}^{n+1}\right)^{d-1}$.
\end{definition}

The following lemma lets us bound the weight of varieties with exponential of the shape we work with (and that we started meeting with proposition \ref{init}) using this variety $N(L,\scrT;\alpha)$.

\begin{lemma} \label{Weyl_global}
If $A$ is a subvariety of $\tau_d^* \calP_{de}^\vee$, the following bound holds:
\begin{align*}
& w_{A \times \tilde W}\left(\left[ A \,\diamond\, \calP_{e(-\scrS_\mu)}^{\diamond(n+1)} \times \tilde W ,\, \alpha\left(f\left( \mathbf{v} + \sigma_L^\scrT(\mathbf{w}) \right)\right) \right]\right) \\
& \leq 2^{-(d-2)} \left(\sup_{(L,\scrT,\alpha) \in A} \dim_{\kappa(L,\scrT,\alpha)} N(L,\scrT;\alpha) + \left(2^{d-1}-d+1\right)(n+1)(e-l'-g+1)\right) .
\end{align*}
\end{lemma}
\begin{proof}
We are going to apply Weyl differencing \ref{Weyl_dif} with the following data. Let $B = A \times \tilde W$ and
\[
\calV = A \,\diamond\, \calP_{e(-\scrS_\mu)}^{\diamond(n+1)} \times \tilde W .
\]
We assumed $e \geq 2g-1+l'$ in (\ref{cond_CM}), so by construction (\ref{exseq_fam_Smu}) of $\calP_{e(-\scrS_\mu)}$ and by Riemann-Roch, $\calV$ is a vector bundle of rank $r = (n+1)(e-l'-g+1)$ over $B$. We define $G$ by
\[
G_b(\mathbf{v}) = \alpha\left(f\left( \mathbf{v} + \sigma_L^\scrT(\mathbf{w}) \right)\right)
\]
for all $b = (L,\scrT,\alpha,\mathbf{w}) \in B$ and $\mathbf{v} \in \calV_b = P_{L(-\scrST)}^{n+1} \times_{\kappa(L,\scrT)} \kappa(b)$. This $G_b$ is polynomial in the coefficients of $\mathbf{v}$, of degree $\leq d$. As a consequence proposition \ref{Weyl_dif} can in fact be applied, and we only need to relate $\calV^{d-1}(G)$ to $N(L,\scrT;\alpha)$. Fix a $b = (L,\scrT,\alpha,\mathbf{w}) \in B$ and consider $\underline{\mathbf{v}} = \left( \mathbf{v}^{(1)}, \dots, \mathbf{v}^{(d-1)} \right) \in \calV_b^{d-1}$ and $\mathbf{u} \in \calV_b$. Lemma \ref{link_multilin} provides
\[
f\left( \cdot + \sigma_L^\scrT(\mathbf{w}) \right)^{(d)}(\underline{\mathbf{v}}, \mathbf{u}) = \sum_{j=0}^n \Psi_j(\underline{\mathbf{v}}) u_j ,
\]
from which we deduce, by linearity of $\alpha$,
\[
G_b^{(d)}(\underline{\mathbf{v}}, \mathbf{u}) = \sum_{j=0}^n \alpha\left(\Psi_j(\underline{\mathbf{v}}) u_j\right) .
\]
By definition (\ref{def_VG}), $\underline{\mathbf{v}}$ is in $\calV^{d-1}(G)_b$ when the above is constant in $\mathbf{u} \in \calV_b$. Suppose this is the case: taking $\mathbf{u} = 0$ shows that the constant value is zero, then taking $u_{j'} = 0$ for all $j'$ other than a given $j$ gives $\alpha\left(\Psi_j(\underline{\mathbf{v}}) u_j\right) = 0$; clearly the converse is true. By definition \ref{def_N_global}, and recalling that $b = (L,\scrT,\alpha,\mathbf{w})$, this shows that
\begin{equation}
\calV^{d-1}(G)_b = N(L,\scrT;\alpha) \times_{\kappa(L,\scrT)} \kappa(b) .
\end{equation}
In particular
\begin{align*}
\dim_B \calV^{d-1}(G) & \leq \sup_{b \in B} \dim_{\kappa(b)} \calV^{d-1}(G)_b \\
& = \sup_{(L,\scrT,\alpha) \in A} \dim_{\kappa(L,\scrT)} N(L,\scrT;\alpha) ,
\end{align*}
which completes the proof.
\end{proof}

Our general goal being to bound the weight in the above statement, it suffices now to bound the dimension of each $N(L,\scrT;\alpha)$ (we can restrict to closed $\alpha \in P_{L^{\otimes d}}^\vee \times_{\kappa(L)} \kappa(L,\scrT)$): we fix a $(L,\scrT)$ and write $\kappa = \kappa(L,\scrT)$ and
\[
N(\alpha) = N(L,\scrT;\alpha) .
\]
Let $\infty$ be a fixed arbitrary closed point of $\scrC$: for any integers $s\geq0$ and $0\leq c\leq d-1$ and any $\alpha \in P_{L^{\otimes d}}^\vee \times_{\kappa(L)} \kappa$ we set
\begin{align} \label{def_Nsc}
N_s^c(\alpha) := & \bigg\{\left. \left(\mathbf{v}^{(1)}, \dots, \mathbf{v}^{(c)}, \mathbf{v}^{(c+1)}, \dots, \mathbf{v}^{(d-1)}\right) \in \left(P_{L(-\scrST-s\infty)}^{n+1}\right)^c \times \left(P_{L(-\scrST)}^{n+1}\right)^{d-1-c} \,\right|\, \nonumber \\
& \quad \forall j, \forall u \in P_{L(-\scrST+cs\infty)} ,\, \alpha\left( \Psi_j(\underline{\mathbf{v}}) \cdot u \right) = 0 \bigg\} .
\end{align}
This corresponds to the $N(\alpha)$ from definition \ref{def_N_global}, where the size of the boxes in which the first $c$ components live has been shrunk. In particular note that
\[
N(\alpha) = N_s^0(\alpha) ,
\]
and denote
\begin{equation}
N_s(\alpha) := N_s^{d-1}(\alpha) .
\end{equation}
The strategy is now to relate the dimensions of $N_s^c(\alpha)$ and $N_s^{c+1}(\alpha)$ in order to obtain a relation between the two extremal ones, then to bound $\dim_\kappa N_s(\alpha)$ for an appropriate $s$. This is the standard process when showing some version of Weyl's inequality. Here we adapt sections 7 and 8 from \cite{HL24}. His work is presented in positive characteristic, but everything works in characteristic zero as well, provided that one uses dimensions instead of cardinalities.

Let us start with section 7. We fix a $0 \leq c \leq d-2$, then some
\[
\mathbf{v}^{(1)}, \dots, \mathbf{v}^{(c)} \in P_{L(-\scrST-s\infty)}^{n+1} \quad\text{and}\quad \mathbf{v}^{(c+2)}, \dots, \mathbf{v}^{(d-1)} \in P_{L(-\scrST)}^{n+1}
\]
(that is, all of them except for $\mathbf{v}^{(c+1)}$). We denote their $(d-2)$-tuple as $\widehat{\underline{\mathbf{v}}}$. When a given $\mathbf{v}^{(c+1)}$ is involved, we will write $\underline{\mathbf{v}}$ for the combined data of the fixed $\widehat{\underline{\mathbf{v}}}$ and that $\mathbf{v}^{(c+1)}$. Let us set
\begin{equation} \label{def_Nv}
N_s^c(\alpha)_{\widehat{\underline{\mathbf{v}}}} := \left\{\left. \mathbf{v}^{(c+1)} \in P_{L(-\scrST)}^{n+1} \,\right|\, \forall j, \forall u \in P_{L(-\scrST + cs\infty)},\,  \alpha\left( \Psi_j(\underline{\mathbf{v}}) \cdot u \right) = 0 \right\} 
\end{equation}
and
\[
N_s^{c+1}(\alpha)_{\widehat{\underline{\mathbf{v}}}}' := \left\{\left. \mathbf{v}^{(c+1)} \in P_{L(-\scrST-s\infty)}^{n+1} \,\right|\, \forall j, \forall u \in P_{L(-\scrST + (c+1)s\infty)},\,  \alpha\left( \Psi_j(\underline{\mathbf{v}}) \cdot u \right) = 0 \right\} :
\]
they are fibres over $\widehat{\underline{\mathbf{v}}}$ of the projection from respectively $N_s^c(\alpha)$ and $N_s^{c+1}(\alpha)$ on all entries but the $(c+1)$-th one. Next is a construction of a vector bundle on $\scrC_\kappa = \scrC \times \kappa$ that realises these fibres. Classically, this step relies on geometry of numbers instead, but in \cite{HL24} since we are working with global sections over a curve, it is necessary to translate the usual process into the theory of vector bundles and their slopes.

\begin{remark}
Consider the field
\[
K_\infty := \Frac\left(\hatOCinf\right) .
\]
Thanks to the smoothness of $\scrC_\kappa$, it is isomorphic to $\kappa((t))$ with a uniformiser $t$. As explained in section 6 of \cite{HL24}, our linear form $\alpha$ can be lifted to an element of $K_\infty$, that we still \break denote $\alpha$. Furthermore, by pulling back the involved spaces of sections to $\Spec K_\infty$ we may see $\mathbf{v}^{(c+1)} \mapsto \Psi_j(\underline{\mathbf{v}})$ as a linear form with coefficients in $K_\infty$.
\end{remark}
Let $K_\scrC$ be the canonical divisor on $\scrC_\kappa$ (not to be confused with $K_\infty$). Consider the vector bundles
\[
\left.\left( L(-\scrST)^{n+1} \oplus L^\vee(\scrST - cs\infty + K_\scrC)^{n+1} \right)\right|_{\scrC_\kappa\smallsetminus\{\infty\}} \longrightarrow \scrC_\kappa\smallsetminus\{\infty\}
\]
and
\[
\hatOCinf^{2n+2} \longrightarrow \hatOCinf ,
\]
as well as the invertible matrix
\[
\begin{bmatrix}
\Id_{n+1} & 0 \\
t^{-cs}\alpha\gamma & \Id_{n+1}
\end{bmatrix}
\in K_\infty^{(2n+2)\times(2n+2)} ,
\]
where $\gamma$ is the $K_\infty$-matrix whose rows correspond to $\mathbf{v}^{(c+1)} \mapsto \Psi_j(\underline{\mathbf{v}})$, as in the previous remark, for $0\leq j \leq n$. By the Beauville-Laszlo theorem (\cite{SP}, tag 0BP2), these two bundles can be glued via this transition matrix, resulting in a vector bundle\footnote{$E^c$ also depends on $s$, $\alpha$ and $\widehat{\underline{\mathbf{v}}}$ but we omit them from the notation.}
\[
E^c \longrightarrow \scrC_\kappa .
\]
That is the construction from \cite{HL24} with $\mathcal{O}(D)$ replaced with $L(-\scrST)$.
\begin{lemma}[\cite{HL24}, lemmas 31 and 32] \label{link_E_N}
One has
\[
H^0(\scrC,E^c) \simeq N_s^c(\alpha)_{\widehat{\underline{\mathbf{v}}}} \quad\text{and}\quad H^0(\scrC,E^c(-s\infty)) \simeq N_s^{c+1}(\alpha)_{\widehat{\underline{\mathbf{v}}}}' .
\]
\end{lemma}
\begin{proof}
See Hase-Liu. Here is a short summary of the justification for the first isomorphism, the second one being similar. By construction a global section of $E^c$ is equivalent to a $\mathbf{v} \in P_{L(-\scrST)}^{n+1}$ and a
\[
\mathbf{v}' \in H^0\left(\scrC_\kappa\smallsetminus\{\infty\}, L^\vee(\scrST - cs\infty + K_\scrC)\right)^{n+1}
\]
such that $v_\infty(\alpha\gamma\mathbf{v} + \mathbf{v}') \geq 0$. For a given $\mathbf{v}$, such a $\mathbf{v}'$ is unique, so the data of a global section equates that of $\mathbf{v}$ alone with some condition. This condition reformulates as $\alpha\gamma\mathbf{v} = 0$ in $P_{L(-\scrST+cs\infty)}^\vee$, which is the definition (\ref{def_Nv}) of $\mathbf{v} \in N_s^c(\alpha)_{\widehat{\underline{\mathbf{v}}}}$.
\end{proof}

As announced, this can be used to show the following.
\begin{lemma}[Shrinking lemma, global] \label{shrink_global}
As soon as $s \geq \max(2g-1, 2)$ we have
\[
\dim_\kappa N(\alpha) \leq (d-1)(n+1)s + (n+1)\frac{g+1}{2}\mathbbm{1}_{g\geq2} + \dim_\kappa N_s(\alpha) .
\]
\end{lemma}
\begin{proof}
Just like in \cite{HL24}, one shows that $(E^c)^\vee \otimes K_\scrC(-cs\infty) \simeq E^c$. This identity, together with slope theory, is the key to obtain inequalities between some $h^0$ and $h^1$ of the Harder–Narasimhan filtration of $E^c$. Eventually it is shown that if $s \geq \max(2g-1,2)$ then
\begin{equation} \label{bound_h0E}
h^0(\scrC, E^c) - h^0(\scrC, E^c(-s\infty)) \leq
\left\{ \begin{array}{ll}
(n+1)s & \text{if } g\leq1 \text{ or } c\geq1 \\
(n+1)\left(s + \frac{g+1}{2}\right) & \text{if } g\geq2 \text{ and } c=0
\end{array} \right.
\end{equation}
(\cite{HL24}, lemma 38). The only thing to note here compared to Hase-Liu's version is that the case $g=0$ can be included along with $g=1$.

Then we can apply lemma \ref{link_E_N} above : the difference $\dim_\kappa N_s^c(\alpha)_{\widehat{\underline{\mathbf{v}}}} - \dim_\kappa N_s^{c+1}(\alpha)_{\widehat{\underline{\mathbf{v}}}}'$ admits the same bound as (\ref{bound_h0E}), which is independent of $\widehat{\underline{\mathbf{v}}}$. As a consequence and by definition (\ref{def_Nv}) of these fibres, this inequality can be lifted to
\[
\dim_\kappa N_s^c(\alpha) - \dim_\kappa N_s^{c+1}(\alpha) \leq
\left\{ \begin{array}{ll}
(n+1)s & \text{if } g\leq1 \text{ or } c\geq1 \\
(n+1)\left(s + \frac{g+1}{2}\right) & \text{if } g\geq2 \text{ and } c=0
\end{array} \right. .
\]
By chaining this relation for all $0 \leq c \leq d-1$ we obtain the claimed bound.
\end{proof}

To finish, we establish a dimensional bound for $N_s(\alpha)$ when $s$ is chosen large enough.
\begin{lemma} \label{bound_Ns_global}
If
\[
s > \max\left( \frac{\deg\alpha-e+l'+2g-2}{d-1},\, e-l'-\frac{\deg\alpha}{d-1} \right)
\]
then
\[
\dim_\kappa N_s(\alpha) \leq (d-2)(n+1) h^0(\scrC, L(-\scrST-s\infty)) .
\]
\end{lemma}
\begin{proof}
This is the object of section 8 of \cite{HL24}. The proof of his lemma 40 shows (the only difference being the replacement of $e = \deg L$ with $e-l' = \deg(L(-\scrST))$) that under the assumption on $s$, the condition
\[
\forall u \in P_{L(-\scrST+(d-1)s\infty)} ,\, \alpha\left( \Psi_j(\underline{\mathbf{v}}) \cdot u \right) = 0
\]
from the definition (\ref{def_Nsc}) of $N_s(\alpha)$ equates to $\Psi_j(\underline{\mathbf{v}}) = 0$. From there, since $\mathbf{v}^{(1)}, \dots, \mathbf{v}^{(d-1)} \in P_{L(-\scrST-s\infty)}$, the final lines of section 8 justify the claimed bound.
\end{proof}

For $\delta \in \mathbb{N}$ we set
\begin{equation} \label{def_s}
s(\delta) := \max_{1\leq i\leq3} s_i(\delta)
\end{equation}
with
\[
\left\{\begin{array}{rl}
\bigg. s_1(\delta) := & \left\lfloor \frac{\delta-e+l'+2g-2}{d-1} \right\rfloor + 1 \\
\bigg. s_2(\delta) := & e-l'+1 - \left\lceil \frac{\delta}{d-1} \right\rceil \\
\bigg. s_3(\delta) = s_3 := & \max(2g-1, 2)
\end{array}\right. .
\]

The previous results combine to provide the following one.
\begin{proposition}[Weyl's inequality, global] \label{concl_Weyl_global}
Suppose that $A$ is a subvariety of $\tau_d^* \calP_{de}^\vee$ such that for all $(L,\scrT,\alpha) \in A$, $s(\deg\alpha) \leq e-l'-2g+1$: then
\begin{align*}
& w_{A \times \tilde W}\left(\left[ A \,\diamond\, \calP_{e(-\scrS_\mu)}^{\diamond(n+1)} \times \tilde W ,\, \alpha\left(f\left( \mathbf{v} + \sigma_L^\scrT(\mathbf{w}) \right)\right) \right]\right) \\
& \leq 2^{-(d-2)}(n+1) \left(\sup_{(L,\scrT,\alpha) \in A} s(\deg\alpha) + \left(2^{d-1}-1\right)(e-l'-g+1) + \frac{g+1}{2}\mathbbm{1}_{g\geq2} \right) .
\end{align*}
\end{proposition}
\begin{proof}
The hypothesis on $s = s(\deg\alpha)$ ensures that we can apply Riemann-Roch to get
\[
h^0(\scrC, L(-\scrST-s\infty)) = (e-l'-s)-g+1 .
\]
Putting this together with lemmas \ref{bound_Ns_global} and \ref{shrink_global} provides
\[
\dim_\kappa N(L,\scrT;\alpha) \leq (n+1)\left( s(\deg\alpha) + (d-2)(e-l'-g+1) + \frac{g+1}{2}\mathbbm{1}_{g\geq2} \right) ,
\]
which we plug in lemma \ref{Weyl_global} to conclude.
\end{proof}

\subsubsection{Local inequality}
\label{Weyl_shrink_local}

We now establish a `local' version of Weyl's inequality, with sections living in $P_{-\scrST|Z}$. This version is more classical than the previous one because such sections are essentially polynomials. Fix closed points $\scrT \in \Sym^\mu(\CmS)$ and $Z \in \Sym^\nu(\scrC)$. For $\beta \in P_{|Z}^\vee$ we define an alternative local version of the $N(\alpha)$ from \ref{def_N_global}:
\begin{equation} \label{def_N'}
N'(\beta) := \left\{\left. \underline{\mathbf{r}} \in \left(P_{-\scrST|Z}^{n+1}\right)^{d-1} \,\right|\, \forall j, \forall u \in P_{-\scrST|Z},\, \beta\left( \Psi_j(\underline{\mathbf{r}}) \cdot u \right) = 0 \right\} .
\end{equation}
We also denote $\delta := \Sigma\nu = \deg Z$ and $\delta' := \deg(\scrST \cap Z)$, and take a section $\sigma' : \calP_{|\scrS_\mu \cap \nu} \to \calP_{|\nu}$ of the exact sequence
\[
\xymatrix{
0 \ar[r] & \calP_{-\scrS_\mu|\nu} \ar[r] & \calP_{|\nu} \ar[r] & \calP_{|\scrS_\mu \cap \nu} \ar[r] & 0
}
\]
from the bottom of (\ref{exseq_9}). Let $\sigma_{|Z}^\scrT$ be the composition of the restriction $P_{|\scrST} \twoheadrightarrow P_{|\scrST\cap Z}$ with $\sigma'_{\scrT,Z} : P_{|\scrST\cap Z} \to P_{|Z}$, so that
\begin{equation} \label{def_sect_sigma_2}
\xymatrix @C=0.5pc {
P_{|\scrST} \ar@{->>}[rd] \ar[rr]^{\sigma_{|Z}^\scrT} && P_{|Z} \ar@{->>}[ld] \\
& P_{|\scrST\cap Z} &
} .
\end{equation}
Like with the $\sigma_L^\scrT$ from (\ref{def_sect_sigma_1}), for $\mathbf{w} \in \tilde W \subset P_{|\scrS}^{n+1}$ we will write $\sigma_{|Z}^\scrT(\mathbf{w}) := \sigma_{|Z}^\scrT(\mathbf{w},0) \in P_{|Z}^{n+1}$. As with the global case we start with an application of Weyl differencing.

\begin{proposition}
For $A$ a subvariety of $P_{|Z}^\vee$ we have
\begin{align*}
& w_{A \times \tilde W}\left(\left[ A \times P_{-\scrST|Z}^{n+1} \times \tilde W ,\, \beta\left(f\left( \mathbf{r} + \sigma_{|Z}^{\scrT}(\mathbf{w}) \right)\right) \right]\right) \\
& \leq 2^{-(d-2)} \left(\sup_{\beta \in B} \dim N'(\beta) + \left(2^{d-1}-d+1\right)(n+1)(\delta-\delta')\right) .
\end{align*}
\end{proposition}
\begin{proof}
Apply proposition \ref{Weyl_dif} with $B = A \times \tilde W$, $\calV = A \times P_{-\scrST|Z}^{n+1} \times \tilde W$, $r = (n+1)(\delta-\delta')$, and
\[
G_{\beta,\mathbf{w}}(\mathbf{r}) = \beta\left(f\left( \mathbf{r} + \sigma_{|Z}^{\scrT}(\mathbf{w}) \right)\right) .
\]
Just like in the global case one then shows, using lemma \ref{link_multilin}, that $\calV^{d-1}(G)_{\beta,\mathbf{w}} = N'(\beta)$ .
\end{proof}

Like we did in the previous paragraph with $N(\alpha)$, the next goal is to bound the dimension of $N'(\beta)$, mainly using a shrinking argument. This $N'(\beta)$ is however of a more number-theoretic nature than the previous $N(\alpha)$, and the most direct way to deal with it is to compare it with the version studied in \cite{BB25}. That is the object of what follows.

Let us switch to a more polynomial point of view (reminiscent of remark \ref{why_Pvee}): we start by introducing some notation. If $\alpha = \sum_{i=-\infty}^r a_i t^i \in \C((t^{-1}))$ is a formal Laurent series in $t^{-1}$ with $a_r \neq 0$, its order is $\ord\alpha := r$ (or $-\infty$ for $\alpha=0$; this extends the degree on $\C[t]$), its fractional part is $\{\alpha\} := \sum_{i=-\infty}^{-1} a_i t^i$, and its residue is $\res\alpha := a_{-1}$. Note that $\C(t) \subset \C((t^{-1}))$.
\begin{lemma} \label{N'_poly}
Let $h_\cap | h$ be polynomials with the same partitions as $\scrST \cap Z \subset Z$. If $\beta \in P_{|Z}^{\vee*}$, there exists a $g \in \C[t]_{<\delta}$ coprime with $h$ such that
\[
N'(\beta) \simeq N'\left(\frac{g}{h}\right) := \left\{\left. \underline{\mathbf{r}} \in \left(h_\cap \C[t]_{<\delta-\delta'}^{n+1}\right)^{d-1} \,\right|\, \forall j,\, \ord\left\{ \frac{g}{h} h_\cap \Psi_j(\underline{\mathbf{r}}) \right\} < -(\delta-\delta') \right\} .
\]
\end{lemma}
\begin{remark}
In the case $\scrST=\varnothing$ we retrieve exactly the $N\left(\frac{g}{h}\right)$ from \cite{BB25}, proposition 5.4:
\[
\left\{\left. \underline{\mathbf{r}} \in \left(\C[t]_{<\delta}^{n+1}\right)^{d-1} \,\right|\, \forall j,\, \ord\left\{ \frac{g}{h} \Psi_j(\underline{\mathbf{r}}) \right\} < -\delta \right\} .
\]
\end{remark}

\begin{proof}[Proof of lemma \ref{N'_poly}]
After (\ref{exseq_ed_loc}), consider the isomorphism $P_{|Z} \overset{\sim}{\to} \C[t]/h$ that restricts to $P_{-\scrST|Z} \overset{\sim}{\to} h_\cap \C[t]/h$. It induces a dual isomorphism: let $\beta' \in (\C[t]/h)^\vee$ be the image of $\beta$ via the latter. Then, with the definition (\ref{def_N'}) of $N'(\beta)$ in mind,
\begin{align} \label{N'_intermed}
N'(\beta) & \simeq \left\{\left. \underline{\mathbf{r}}' \in \left(\left(h_\cap \C[t]/h\right)^{n+1}\right)^{d-1} \,\right|\, \forall j, \forall u' \in h_\cap \C[t]/h,\, \beta'\!\left( \Psi_j(\underline{\mathbf{r}}') \cdot u' \right) = 0 \right\} \nonumber \\
& \simeq \left\{\left. \underline{\mathbf{r}}'' \in \left(h_\cap \C[t]_{<\delta-\delta'}^{n+1}\right)^{d-1} \,\right|\, \forall j, \forall i \in \llbracket0;\delta-\delta'-1\rrbracket ,\, \beta'\!\left( \Psi_j(\underline{\mathbf{r}}'') \cdot t^i h_\cap \bmod h \right) = 0 \right\} .
\end{align}

We claim that we have $\beta'(\cdot \bmod h) = \res\left(\frac{g}{h} \cdot\right)$ for some $g \in \C[t]_{<\delta}$ coprime with $h$. In fact let us endow $\C[t]$ with its usual hermitian product
\[
{\textstyle \left\langle \sum_i a_i t^i, \sum_i b_i t^i \right\rangle = \sum_i \overline{a_i} b_i },
\]
where $\overline{\cdot}$ is the complex conjugation. The linear form
\[
\beta'(\cdot \bmod h) : \C[t] \longrightarrow \C[t]/h \longrightarrow \C
\]
is continuous, so by the Riesz representation theorem it can be written $\langle g', \cdot \rangle$ for some $g' \in \C[t]$. Note that for any $v \in \C[t]$, $\langle g', v \rangle$ is the constant coefficient of $\overline{g'}(t^{-1}) v$, that is $\res\left( t^{-1}\overline{g'}(t^{-1}) v \right)$. Since the kernel of our linear form contains $(h)$, for all $i \geq 0$ we have $\res\left( t^{-1} \overline{g'}(t^{-1}) \cdot t^i h \right) = 0$; but this is the $(-i-1)$-th coefficient of $t^{-1} \overline{g'}(t^{-1}) h$: the latter is therefore a polynomial, that we denote $g$. Its degree is $< \deg h = \delta$. Finally $t^{-1} \overline{g'}(t^{-1}) = \frac{g}{h}$ and
\[
\beta'(\cdot \bmod h) = \res\left(\frac{g}{h} \cdot\right) .
\]
We still need to check that $g$ is coprime with $h$. Suppose that $g=h'v$ with $h'$ a nontrivial divisor of $h$: our linear form $\res\left(\frac{g}{h} \cdot\right) = \res\left(\frac{v}{h/h'} \cdot\right)$ factors through $\C[t]/(h/h')$, and $\beta'$ as well. By construction, our original $\beta$ then factors through $P_{|Z'}$ for some $Z' \subsetneq Z$, which contradicts the hypothesis $\beta \in P_{|Z}^{\vee*}$. This ends the proof of the claim.

The condition in (\ref{N'_intermed}) becomes, for all $0 \leq j \leq n$ and $0 \leq i \leq \delta-\delta'-1$,
\[
\res\left( \frac{g}{h} \Psi_j(\underline{\mathbf{r}}) \cdot t^i h_\cap \right) = 0
\]
(we go back to writing $\underline{\mathbf{r}}$ instead of $\underline{\mathbf{r}}''$). This residue is the $(-i-1)$-th coefficient of $\frac{g}{h} h_\cap \Psi_j(\underline{\mathbf{r}})$: altogether we demand the coefficients $-(\delta-\delta')$ to $-1$ to vanish. This can be reformulated as
\[
\ord\left\{ \frac{g}{h} h_\cap \Psi_j(\underline{\mathbf{r}}) \right\} < -(\delta-\delta') .
\]
\end{proof}

Now that the variety we are interested in is very similar to the one from \cite{BB25}, we follow their proof. The strategy remains the same: first shrinking the size of the `box' we are taking our tuples in, and then bounding the dimension of the shrunk variety. Here the tool is classical: geometry of numbers (and not slope theory as in \cite{HL24} and the previous paragraph). This process is the object of lemma 5.7 of \cite{BB25}, in particular the first lines of the corresponding proof deal with the shrinking step. The slight differences in the studied variety ($h_\cap$ divides $\underline{\mathbf{r}}$ and appears in the condition) do not impact this part of the reasoning, and in our case the result is as follows.

\begin{lemma}[Shrinking lemma, local]
With $h, h_\cap, g$ as above, for any integer $s\geq0$ we have
\[
\dim N'\left(\frac{g}{h}\right) \leq (d-1)(n+1)s + \dim N'_s\left(\frac{g}{h}\right)
\]
with
\[
N'_s\left(\frac{g}{h}\right) := \left\{\left. \underline{\mathbf{r}} \in \left(h_\cap \C[t]_{<\delta-\delta'-s}^{n+1}\right)^{d-1} \,\right|\, \forall j,\, \ord\left\{ \frac{g}{h} h_\cap \Psi_j(\underline{\mathbf{r}}) \right\} < -(\delta-\delta')-(d-1)s \right\} .
\]
\end{lemma}

Then comes the time to pick an $s$ large enough for $\dim N_s'(\frac{g}{h})$ to be bounded easily.

\begin{lemma}
If
\[
s > \frac{(d-2)\delta+\delta'}{d-1} - 1
\]
then
\[
\dim N'_s\left(\frac{g}{h}\right) \leq (d-2)(n+1)(\delta-\delta'-s) .
\]
\end{lemma}
\begin{proof}
We adapt the proof of lemma 5.7 from \cite{BB25}, more exactly the case of a fraction (we do not need the case of a small formal Laurent series $\theta$). The first step for us is to remark that the conditions in the definition of $N'_s\left(\frac{g}{h}\right)$ are equivalent to $h/h_\cap$ dividing $\Psi_j(\underline{\mathbf{r}})$ (in particular $s$ does not actually appear in it, and only plays a role in the degree of the tuples of polynomials we consider)\footnote{We could have directly remarked this in lemma \ref{N'_poly}, but this would not be convenient for the shrinking lemma.}
. In fact the smallest order that appears in the expansion of
\[
\frac{g}{h} h_\cap \Psi_j(\underline{\mathbf{r}}) = \frac{g \Psi_j(\underline{\mathbf{r}})}{h/h_\cap}
\]
is $\geq -\deg(h/h_\cap) = -(\delta-\delta')$, so:
\begin{align*}
& \ord\left\{ \frac{g \Psi_j(\underline{\mathbf{r}})}{h/h_\cap} \right\} < -(\delta-\delta')-(d-1)s \\
& \Longleftrightarrow\quad \left\{ \frac{g \Psi_j(\underline{\mathbf{r}})}{h/h_\cap} \right\} = 0 \quad\Longleftrightarrow\quad h/h_\cap \big| g \Psi_j(\underline{\mathbf{r}}) \quad\Longleftrightarrow\quad h/h_\cap \big| \Psi_j(\underline{\mathbf{r}})
\end{align*}
(the last equivalence is because $g$ and $h$ are coprime).

Then we want $s$ to be large enough to force the conditions to become $\Psi_j(\underline{\mathbf{r}}) = 0$. After what we just showed, this can be achieved by having $\deg\Psi_j(\underline{\mathbf{r}}) < \delta-\delta'$. By definition of $\Psi_j$ (\ref{def_Psi}) and $N'_s\left(\frac{g}{h}\right)$, one always has $\deg\Psi_j(\underline{\mathbf{r}}) \leq (d-1)(\delta-s-1)$ : it suffices to impose the condition $(d-1)(\delta-s-1) < \delta-\delta'$, which rewrites as
\[
s > \frac{(d-2)\delta+\delta'}{d-1} - 1 ,
\]
hence the hypothesis of the lemma.

Now for such a choice of $s$ we have
\[
N'_s\left(\frac{g}{h}\right) = \left\{\left. \underline{\mathbf{r}}' \in \left(\C[t]_{<\delta-\delta'-s}^{n+1}\right)^{d-1} \,\right|\, \forall j,\, \Psi_j(\underline{\mathbf{r}}') = 0 \right\} .
\]
Lemma 5.6 of \cite{BB25} directly provides
\[
\dim N'_s\left(\frac{g}{h}\right) \leq (d-2)(\delta-\delta'-s)(n+1)
\]
(again, mind that their $n$ is our $n+1$), which is what we announced.
\end{proof}

The results of this paragraph combine to provide the following. Recall that $\delta = \deg Z$ and $\delta' = \deg(\scrST \cap Z)$.
\begin{proposition}[Weyl's inequality, local] \label{concl_Weyl_dif_local}
If $A$ is a subvariety of $P_{|Z}^{\vee*}$ and if we take
\[
s = \left\lfloor \frac{(d-2)\delta+\delta'}{d-1} \right\rfloor ,
\]
then
\begin{align*}
& w_{A \times \tilde W}\left(\left[ A \times P_{-\scrST|Z}^{n+1} \times \tilde W ,\, \beta\left(f\left( \mathbf{r} + \sigma_{|Z}^{\scrT}(\mathbf{w}) \right)\right) \right]\right) \\
& \leq 2^{-(d-2)}(n+1) \left(s + \left(2^{d-1}-1\right)(\delta-\delta')\right) .
\end{align*}
\end{proposition}

\subsection{Major arcs}
\label{subsect_maj}

With appropriate instances of Weyl's inequality in our pocket, we can now work on the major and minor arcs; we start with the major arcs $\mathfrak{M}$. This gives rise to the present subsection, which is longer than the other ones. The goal of the first paragraph is proposition \ref{maj_1}, which expresses $\Gamma_\mathfrak{M}$ as a renormalised sum of varieties with exponentials, simpler than the one we start with and with only local sections. In the next paragraph this sum is completed into the `singular series' $\mathfrak{S}^\mu$ (lemma \ref{sing_cv}). Then begins a succession of rewritings of $\mathfrak{S}^\mu$, the first one of them being as a motivic Euler product (corollary \ref{sing_mep_2}), thanks to the factorisation lemma \ref{lem_fact}. The factors of this product are reformulated as local densities of solutions (more precisely, as limits of normalised jet schemes) as summed up in corollary \ref{sing_locdens}. Finally we deduce an even more explicit form for these factors in proposition \ref{sing_prod_final}.

This general trajectory is inspired by \cite{BB25} (which, again, does the case $\scrC = \projsp^1$, $\scrS=\varnothing$), with of course differences, among which: the varieties with exponentials are more complicated (mainly because of the conditions over $\scrST$), we need to stratify things by partition and not only by degree, and we have no singular integral (see remark \ref{no_intsing}).

\subsubsection{First manipulations}
\label{maj_start}

Going back to the definition of the major arcs under (\ref{def_eta}), let $\mathfrak{M}_\nu$ be the stratum for which $\alpha$ factors minimally through a $Z \in \Sym^\nu(\scrC)$, so that
\[
\mathfrak{M} = \bigsqcup_{\delta=0}^{\eta_e} \bigsqcup_{\nu\vdash\delta} \mathfrak{M}_\nu .
\]
Write
\[
\Gamma_\nu := \mathbf{L}^{-de+g-1} \left[ \mathfrak{M}_\nu \,\diamond\, \calP_{e(-\scrS_\mu)}^{\diamond(n+1)} \times \tilde W ,\, \alpha\left(f\left( \mathbf{v} + \sigma_L^\scrT(\mathbf{w}) \right)\right) \right] ,
\]
in order to have
\[
\Gamma_\mathfrak{M} = \sum_{\delta=0}^{\eta_e} \sum_{\nu\vdash\delta} \Gamma_\nu .
\]
Let's fix a $0 \leq \delta \leq \eta_e$ and a $\nu\vdash\delta$: we are going to reformulate $\Gamma_\nu$. The definition of $\mathfrak{M}_\nu$ places it over $\Sym^\nu(\scrC)$, so we can see $\Gamma_\nu$ as over
\[
\calB_\mu \times \Sym^\nu(\scrC) = J(\scrC) \times \Sym^\mu(\CmS) \times \Sym^\nu(\scrC) \times \tilde W
\]
(recall definition (\ref{def_Bmu})).

For now fix $L \in \Pic^e(\scrC)$, $\scrT \in \Sym^\mu(\CmS)$ and $Z \in \Sym^\nu(\scrC)$: we can consider the residue field $\kappa = \kappa(L, \scrT, Z)$ and the motivic value
\[
\Gamma_\nu(L,\scrT,Z) \in \mathscr{E}xp\mathscr{M}_{\kappa \times \tilde W} .
\]
The fibre $(\mathfrak{M}_\nu)_{L,\scrT,Z}$ is the space of $\alpha \in P_{L^{\otimes d}}^\vee$ that factor minimally through $Z$, which is isomorphic to $P_{|Z}^{\vee*}$ by (\ref{Pvee*}). Thus
\begin{align*}
\Gamma_\nu(L,\scrT,Z) & = \mathbf{L}^{-de+g-1} \left[ P_{|Z}^{\vee*} \times P_{L(-\scrST)}^{n+1} \times \tilde W ,\, \beta\left(\left.f\left( \mathbf{v} + \sigma_L^\scrT(\mathbf{w}) \right)\right|_Z\right) \right]_{\kappa \times \tilde W} \\
& = \mathbf{L}^{-de+g-1} \left[ P_{|Z}^{\vee*} \times P_{L(-\scrST)}^{n+1} \times \tilde W ,\, \beta\left(f\left( \mathbf{v}|_Z + \sigma_L^\scrT(\mathbf{w})|_Z \right)\right) \right]_{\kappa \times \tilde W} .
\end{align*}
Here the exponential doesn't depend on the whole $\mathbf{v}$ but only on $\mathbf{v}|_Z$: since lemma \ref{exseq_ed_vanish} gives us $P_{L(-\scrST)} \simeq P_{L(-(\scrST \cup Z))} \times P_{-\scrST|Z}$, the above rewrites as
\begin{equation} \label{maj_EucliDiv}
\mathbf{L}^{-de+g-1} \left[P_{L(-(\scrST \cup Z))}^{n+1}\right]_{\kappa \times \tilde W} \cdot \left[ P_{|Z}^{\vee*} \times P_{-\scrST|Z}^{n+1} \times \tilde W ,\, \beta\left(f\left( \mathbf{r} + \sigma_L^\scrT(\mathbf{w})|_Z \right)\right) \right]_{\kappa \times \tilde W} .
\end{equation}
Note that $P_{L(-(\scrST \cup Z))}$ is a vector space of dimension
\begin{align*}
e - \deg(\scrST \cup Z) - g + 1 & = e-g+1 - (\deg\scrST + \deg Z - \deg(\scrST \cap Z)) \\
& = e - g + 1 - l' - \delta + \deg(\scrST \cap Z) ,
\end{align*}
so we have
\begin{equation} \label{maj_intermed_1}
\Gamma_\nu(L,\scrT,Z) = \mathbf{L}^{\zeta_{e,|\mu|} - (n+1)(\delta-\deg(\scrST \cap Z))} \left[ P_{|Z}^{\vee*} \times P_{-\scrST|Z}^{n+1} \times \tilde W ,\, \beta\left(f\left( \mathbf{r} + \sigma_L^\scrT(\mathbf{w})|_Z \right)\right) \right]_{\kappa \times \tilde W}
\end{equation}
where
\begin{align} \label{def_zeta}
\zeta_{e,|\mu|} := &\, -de+g-1 + (n+1)(e - g + 1 - l') \nonumber \\
= &\, (n-d+1)e - n(g-1) - (n+1)l'
\end{align}
depends on $|\mu|$ via $l' = l+|\mu|$.

We are now going to show that $\Gamma_\nu(L,\scrT,Z)$ actually does not depend on $L$, and only depends on $(\scrT,Z)$ through the part of $Z$ over $|\scrS|$ and the partition of its part over $|\scrT|$. If $T$ is a subset of $\scrC$ (for instance the support of an effective divisor) we denote by $Z|_T$ the part of $Z$ that is over $T$ (that is, $\big|Z|_T\big| \subset T$ while $\big|Z-Z|_T\big| \cap T = \varnothing$)
\footnote{$Z_1|_{|Z_2|}$ should not be mistaken for $Z_1 \cap Z_2$: it is only true that as effective divisors, $Z_1 \cap Z_2 \leq Z_1|_{|Z_2|}$.}.
As $(\scrT,Z)$ varies, the possible $Z|_{|\scrS|}$ are the effective divisors $X$ such that $|X| \subset |\scrS|$ and $\pt X \leq \nu$, and then the possibilities for $\pt\left(Z|_{|\scrT|}\right)$ are the subpartitions $\nu^0$ of $\nu - \pt X$ with $|\nu^0| \leq |\mu|$. Let $I_\nu$ be the set of such couples $(X, \nu^0)$: it is finite. If $(X, \nu^0) \in I_\nu$, write
\begin{equation} \label{def_Unu}
U_\nu(X, \nu^0) := \left\{ (\scrT,Z) \in \Sym^\mu(\CmS) \times \Sym^\nu(\scrC) \,\left|\,
\begin{array}{l}
Z|_{|\scrS|} = X \\
\pt\left(Z|_{|\scrT|}\right) = \nu^0
\end{array} \right.\right\} :
\end{equation}
these form a finite stratification of $\Sym^\mu(\CmS) \times \Sym^\nu(\scrC)$. It is well defined as a variety, in fact it is isomorphic to
\[
\left\{\left. (\scrT,Y,Y^0) \in \Sym^{\mu, \nu - \pt X - \nu^0}(\CmS) \times \Sym^{\nu^0}(\CmS) \,\right|\, |Y^0| \subset |\scrT| \right\}
\]
where $Y$ and $Y^0$ respectively stand for $Z|_{\scrC \smallsetminus |\scrST|}$ and $Z|_{|\scrT|}$.

Recall the $\sigma_{|Z}^\scrT$ constructed in (\ref{def_sect_sigma_2}). We can take the temporary notation, for every $(\scrT,Z)$ in $\Sym^\mu(\CmS) \times \Sym^\nu(\scrC)$,
\begin{equation}
C(\scrT,Z) := \left[ P_{|Z}^{\vee*} \times P_{-\scrST|Z}^{n+1} \times \tilde W ,\, \beta\left(f\left( \mathbf{r} + \sigma_{|Z}^{\scrT}(\mathbf{w}) \right)\right) \right]_{\kappa(\scrT,Z) \times \tilde W} ,
\end{equation}
and state the following lemma.
\begin{lemma} \label{maj_intermed_2}
For all $L \in \Pic^e(\scrC)$ and $(\scrT,Z) \in U_\nu(X,\nu^0)$, if we denote $\kappa = \kappa(L,\scrT,Z)$ we have
\[
\Gamma_\nu(L,\scrT,Z) = \mathbf{L}^{\zeta_{e,|\mu|} - (n+1)\left(\delta - \deg(\scrS \cap X) - |\nu^0|\right)} \cdot C(\scrT,Z)_{\kappa \times \tilde W} .
\]
\end{lemma}

\begin{proof}
We previously obtained (\ref{maj_intermed_1}):
\[
\Gamma_\nu(L,\scrT,Z) = \mathbf{L}^{\zeta_{e,|\mu|} - (n+1)(\delta-\deg(\scrST \cap Z))} \left[ P_{|Z}^{\vee*} \times P_{-\scrST|Z}^{n+1} \times \tilde W ,\, \beta\left(f\left( \mathbf{r} + \sigma_L^\scrT(\mathbf{w})|_Z \right)\right) \right]_{\kappa \times \tilde W} .
\]
The first thing to see is that in the exponent of the $\mathbf{L}$,
\begin{align*}
\deg(\scrST \cap Z) & = \deg(\scrS \cap Z) + \deg(|\scrT| \cap Z) \\
& = \deg(\scrS \cap X) + |\nu^0|
\end{align*}
because $\deg(|\scrT| \cap Z)$ is the number of points of $Z|_{|\scrT|}$, that is $|\nu^0|$. Next we wish to replace $\sigma_L^\scrT(\mathbf{w})|_Z$ with $\sigma_{|Z}^\scrT(\mathbf{w})$. While they have no reason to be equal, we claim that the translation
\[
\mathbf{r} \longmapsto \mathbf{r} + \sigma_L^\scrT(\mathbf{w})|_Z - \sigma_{|Z}^\scrT(\mathbf{w})
\]
is an automorphism of $P_{-\scrST|Z}^{n+1}$. The only thing to check is that $\sigma_L^\scrT(\mathbf{w})|_Z - \sigma_{|Z}^\scrT(\mathbf{w})$ is itself an element of $P_{-\scrST|Z}^{n+1}$, in other words that its restriction to $\scrST \cap Z$ is zero. By construction, $\sigma_L^\scrT(\mathbf{w})$ and $\sigma_{|Z}^\scrT(\mathbf{w})$ both give $(\mathbf{w}|_{\scrS \cap Z},0) \in P_{|\scrST \cap Z}^{n+1}$ via this restriction, which justifies the claim. Now if $\mathbf{r}'$ is the image of $\mathbf{r}$ by this translation, $\mathbf{r} + \sigma_L^\scrT(\mathbf{w})|_Z = \mathbf{r}' + \sigma_{|Z}^\scrT(\mathbf{w})$, and applying this to the class at the start of the proof we get $C(\scrT,Z)_{\kappa \times \tilde W}$.
\end{proof}

Now we show that $C(\scrT,Z)$ is constant over $U_\nu(X,\nu^0)$.
\begin{lemma} \label{maj_C_cst}
Take $(\scrT,Z), (\scrT',Z') \in U_\nu(X,\nu^0)$ and assume the former is a closed point. Then in $K_0\left(\mathrm{ExpVar}_{\kappa(\scrT',Z') \times \tilde W}\right)$,
\[
C(\scrT',Z') = [\kappa(\scrT',Z')] \cdot C(\scrT,Z) .
\]
\end{lemma}

\begin{proof}
One can write $Y^0 = Z|_{|\scrT|}$ and $Y = Z|_{\scrC \smallsetminus |\scrST|}$, as well as $Y^0{}' = Z'|_{|\scrT'|}$ and $Y' = Z'|_{\scrC \smallsetminus |\scrS_{\scrT'}|}$, so that
\[
Z = X + Y^0 + Y \qquad\text{and}\qquad Z' = X + Y^0{}' + Y' .
\]
By hypothesis $\pt(Y^0) = \nu^0 = \pt(Y^0{}')$, which forces $\pt(Y) = \pt(Y')$. Using lemma \ref{PZ_gen_decomp} we get
\begin{equation} \label{maj_C_cst_intermed1}
P_{|Z'} \simeq \left(P_{|X} \otimes \kappa(Z')\right) \times \left(P_{|Y^0{}'} \otimes_{\kappa(Y^0{}')} \kappa(Z')\right) \times \left(P_{|Y'} \otimes_{\kappa(Y')} \kappa(Z')\right)
\end{equation}
(note that $\kappa(X)=\C$). We can then use lemma \ref{PZ_gen_isom} to write $P_{|Y^0{}'} \simeq P_{|Y^0} \otimes \kappa(Y^0{}')$, which provides
\[
P_{|Y^0{}'} \otimes_{\kappa(Y^0{}')} \kappa(Z') \simeq P_{|Y^0} \otimes \kappa(Z') .
\]
Likewise the third factor of (\ref{maj_C_cst_intermed1}) is isomorphic to $P_{|Y} \otimes \kappa(Z')$, so
\begin{align*}
P_{|Z'} & \simeq \left(P_{|X} \otimes \kappa(Z')\right) \times \left(P_{|Y^0} \otimes \kappa(Z')\right) \times \left(P_{|Y} \otimes \kappa(Z')\right) \\
& = \left( P_{|X} \times P_{|Y^0} \times P_{|Y} \right) \otimes \kappa(Z') \\
& = P_{|Z} \otimes \kappa(Z')
\end{align*}
We call the isomorphism of $\kappa(Z')$-algebras we constructed
\begin{equation} \label{maj_C_cst_intermed2}
\varphi : P_{|Z'} \overset{\sim}{\longrightarrow} P_{|Z} \otimes \kappa(Z') .
\end{equation}
Note that via
\[
P_{-\scrS_{\scrT'}|Z'} \overset{\sim}{\longrightarrow} \left(P_{-\scrS|X} \otimes \kappa(Z')\right) \times \{0\} \times \left(P_{|Y'} \otimes_{\kappa(Y')} \kappa(Z')\right)
\]
and
\[
P_{-\scrST|Z} \overset{\sim}{\longrightarrow} P_{-\scrS|X} \times \{0\} \times P_{|Y} ,
\]
$\varphi$ restricts to
\[
P_{-\scrS_{\scrT'}|Z'} \overset{\sim}{\longrightarrow} P_{-\scrST|Z} \otimes \kappa(Z')
\]
(remark that directly applying lemma \ref{PZ_gen_isom} to get an isomorphism (\ref{maj_C_cst_intermed2}) would not have implied that it restricts in this way, hence the prior decomposition we performed on $P_{|Z}$ and $P_{|Z'}$).
We also consider the isomorphism dual to $\varphi$,
\[
\psi : P_{|Z'}^\vee \overset{\sim}{\longrightarrow} P_{|Z}^\vee \otimes \kappa(Z') ,
\]
which clearly restricts to $P_{|Z'}^{\vee*} \overset{\sim}{\to} P_{|Z}^{\vee*} \otimes \kappa(Z')$. Let $\beta' \in P_{|Z'}^{\vee*}$, $\mathbf{r}' \in P_{-\scrS_{\scrT'}|Z'}^{n+1}$, $\mathbf{w} \in \tilde W$: extending $\varphi$ to $(n+1)$-tuples, we have
\begin{align*}
\beta'\left(f\left( \mathbf{r}' + \sigma_{|Z'}^{\scrT'}(\mathbf{w}) \right)\right) & = \psi(\beta')\left( \varphi\left(f\left( \mathbf{r}' + \sigma_{|Z'}^{\scrT'}(\mathbf{w}) \right)\right) \right) \\
& = \psi(\beta')\left( f\left( \varphi(\mathbf{r}') + \varphi\big(\sigma_{|Z'}^{\scrT'}(\mathbf{w})\big) \right) \right) .
\end{align*}
As a consequence, by respectively applying isomorphisms $\psi$ and $\varphi$ to the first and second factors of
\[
C(\scrT',Z') = \left[ P_{|Z'}^{\vee*} \times P_{-\scrS_{\scrT'}|Z'}^{n+1} \times \tilde W ,\, \beta'\left(f\left( \mathbf{r}' + \sigma_{|Z'}^{\scrT'}(\mathbf{w}) \right)\right) \right]_{\kappa(\scrT',Z') \times \tilde W} ,
\]
we obtain
\begin{equation} \label{maj_intermed_3}
\left[ P_{|Z}^{\vee*} \times P_{-\scrST|Z}^{n+1} \times \kappa(Z') \times \tilde W ,\, \beta\left(f\left( \mathbf{r} + \varphi\big(\sigma_{|Z'}^{\scrT'}(\mathbf{w})\big) \right)\right) \right]_{\kappa(\scrT',Z') \times \tilde W}
\footnote{Inside the symbol $[\cdot]$, in order to emphasize that we see the spaces of sections as varieties, we write a fibre product instead of a tensor product to extend the scalars (which has the same effect in this context).} .
\end{equation}
To retrieve $C(\scrT,Z)$, we want to substitute $\varphi\big(\sigma_{|Z'}^{\scrT'}(\mathbf{w})\big)$ for $\sigma_{|Z}^\scrT(\mathbf{w})$. Like in the proof of lemma \ref{maj_intermed_2} these two have no reason to be equal, but the translation
\[
\mathbf{r} \longmapsto \mathbf{r} + \varphi\big(\sigma_{|Z'}^{\scrT'}(\mathbf{w})\big) - \sigma_{|Z}^\scrT(\mathbf{w})
\]
is an automorphism of $P_{-\scrST|Z}^{n+1} \otimes \kappa(Z')$. In fact $\varphi\big(\sigma_{|Z'}^{\scrT'}(\mathbf{w})\big)$ and $\sigma_{|Z}^\scrT(\mathbf{w})$ have the same restriction to $\scrST \cap Z = (\scrS \cap X) + (|\scrT| \cap Y^0)$: by construction of $\varphi$, $\sigma_{|Z}^\scrT$ and $\sigma_{|Z'}^{\scrT'}$, on the one hand we have
\[
\varphi\big(\sigma_{|Z'}^{\scrT'}(\mathbf{w})\big)|_{\scrS \cap X} = \sigma_{|Z'}^{\scrT'}(\mathbf{w})|_{\scrS \cap X} = \mathbf{w}|_{\scrS \cap X} = \sigma_{|Z}^\scrT(\mathbf{w})|_{\scrS \cap X} ,
\]
while one the other hand,
\[
\varphi\big(\sigma_{|Z'}^{\scrT'}(\mathbf{w})\big)|_{|\scrT| \cap Y^0} = \varphi\left(\sigma_{|Z'}^{\scrT'}(\mathbf{w})|_{|\scrT'| \cap Y^0{}'}\right) = \varphi(0) = 0 = \sigma_{|Z}^\scrT(\mathbf{w})|_{|\scrT| \cap Y^0} .
\]
Applying this translation to (\ref{maj_intermed_3}) gives
\begin{align*}
& \left[ P_{|Z}^{\vee*} \times P_{-\scrST|Z}^{n+1} \times \kappa(Z') \times \tilde W ,\, \beta\left(f\left( \mathbf{r} + \sigma_{|Z}^\scrT(\mathbf{w}) \right)\right) \right]_{\kappa(\scrT',Z') \times \tilde W} \\
& = \left[ P_{|Z}^{\vee*} \times P_{-\scrST|Z}^{n+1} \times \kappa(\scrT',Z') \times \tilde W ,\, \beta\left(f\left( \mathbf{r} + \sigma_{|Z}^\scrT(\mathbf{w}) \right)\right) \right] \\
& = [\kappa(\scrT',Z')] \cdot C(\scrT,Z) .
\end{align*}
\end{proof}

In particular $C(\scrT,Z)$ is the same for any closed point $(\scrT,Z) \in U_\nu(X,\nu^0)$: we denote this class
\begin{equation} \label{def_S'}
S'_{\nu;X,\nu^0} := \left[ P_{|Z}^{\vee*} \times P_{-\scrST|Z}^{n+1} \times \tilde W ,\, \beta\left(f\left( \mathbf{r} + \sigma_{|Z}^{\scrT}(\mathbf{w}) \right)\right) \right] \in K_0\left(\mathrm{ExpVar}_{\tilde W}\right) .
\end{equation}
Like we briefly mentioned at the start of this subsection \ref{subsect_maj}, classically, and in equation (6.7) of \cite{BB25} in particular, at this point one would have stratified the polynomials by degree. The analog for us is the space of sections $P_{-\scrST|Z}$ with $\pt Z = \nu$: we need to keep track of the partition because the isomorphisms used in the proof above rely on it; this gives rise to a more refined stratification.

\begin{lemma} \label{Gamma_S'}
The following identity holds over  $\calB_\mu \times \Sym^\nu(\scrC)$:
\[
\Gamma_\nu = \sum_{(X, \nu^0) \in I_\nu} \mathbf{L}^{\zeta_{e,|\mu|} - (n+1)\left(\delta - \deg(\scrS \cap X) - |\nu^0|\right)} \cdot \left[J(\scrC) \times U_\nu(X,\nu^0)\right] \cdot S'_{\nu;X,\nu^0} .
\]
\end{lemma}
\begin{proof}
By lemma \ref{maj_C_cst}, the $C(\scrT,Z)_{\kappa \times \tilde W}$ from lemma \ref{maj_intermed_2}, where $\kappa = \kappa(L,\scrT,Z)$, can be rewritten
\[
C(\scrT,Z)_{\kappa \times \tilde W} = \left( [\kappa(\scrT,Z)] \cdot S'_{\nu;X,\nu^0} \right)_{\kappa \times \tilde W} = [\kappa] \cdot S'_{\nu;X,\nu^0} ;
\]
as a consequence 
\[
\Gamma_\nu(L,\scrT,Z) = \mathbf{L}^{\zeta_{e,|\mu|} - (n+1)\left(\delta - \deg(\scrS \cap X) - |\nu^0|\right)} \cdot [\kappa] \cdot S'_{\nu;X,\nu^0} .
\]
This shows that the two sides of the seeked identity have the same motivic values for all $(L, \scrT, Z) \in \Pic^e(\scrC) \times U_\nu(X, \nu^0)$ (we use the identification $\Pic^e(\scrC) \simeq J(\scrC)$). Putting together all the $(X,\nu^0) \in I_\nu$ concludes.
\end{proof}

With this result in mind, it will be useful to take the notation
\begin{align} \label{def_S}
S_\nu & := \sum_{(X, \nu^0) \in I_\nu} \mathbf{L}^{(n+1)\left(\deg(\scrS \cap X) + |\nu^0|\right)} \cdot \left[J(\scrC) \times U_\nu(X,\nu^0)\right] \cdot S'_{\nu;X,\nu^0} \\
& \in K_0\left(\mathrm{ExpVar}_{\calB_\mu \times \Sym^\nu(\scrC)}\right) , \nonumber
\end{align}
such that
\[
\Gamma_\nu = \mathbf{L}^{\zeta_{e,|\mu|} - (n+1)\delta} S_\nu .
\]
A concrete example of $S_\nu$ will be computed in lemma \ref{S_i}. To express the contribution of all major arcs, remains only to sum over all partitions $\nu$ of integers at most $\eta_e$.
\begin{proposition} \label{maj_1}
In $\mathscr{E}xp\mathscr{M}_{\calB_\mu}$ one has
\[
\Gamma_\mathfrak{M} = \mathbf{L}^{\zeta_{e,|\mu|}} \sum_{\delta=0}^{\eta_e} \mathbf{L}^{-(n+1)\delta} \sum_{\nu\vdash\delta} S_\nu .
\]
\end{proposition}

\begin{remark} \label{no_intsing}
The above expression is the product of a normalisation and a finite sum (that will become the singular series); classically at this stage one also has a (motivic) integral, called the \emph{singular integral}: see for instance \cite{BB25}, proposition 6.3. Ultimately, this integral usually plays the role of the missing factor of the singular series at $\infty \in \projsp^1$; while we will find out in proposition \ref{sing_mep_1} that the motivic Euler product development of our singular series does not miss any point and is indexed by the whole\footnote{At first. Independently of this, later the factors at $p\in|\scrS|$ will be separated from the rest.} of $\scrC$. This is a consequence of our computations preserving the symmetry without singling out a point of the curve (which was a natural thing to do in genus zero, since one usually works with polynomials): the data that would go in the singular integral is already contained in our singular series.

More concretely, if this singular integral were to appear, that would be in step (\ref{maj_EucliDiv}), when we apply the Euclidean division exact sequence to the tuple of sections. The factor corresponding to the quotient has no exponential and is just absorbed into the power of $\mathbf{L}$. However in \cite{BB25} (start of subsection 6.2) this is not the case: the factor with the quotient still has an exponential, that depends on the `$\theta$', i.e. the quality of the approximation of the Laurent series $\alpha$ by a rational fraction; and this factor becomes the singular integral. With the first few lines of our subsection \ref{def_arcs} in mind, this suggests that the preservation of symmetry between all the points of the curve in our framework is responsible for the fact that the `approximation of $\alpha$' is simply replaced with the factorisation of $\alpha$ by an effective divisor, and that we have no direct analog to the `$\theta$'.
\end{remark}

\subsubsection{Convergence of the singular series}

In the previous proposition, the sum of exponential sums goes up to $\eta_e = e-l'-2g+1$. This quantity goes to $\infty$ as $e\to\infty$, so our next goal is to replace the sum by the full series, up to an error term (namely, the remainder of the series).

In the technical computations of this paragraph and further on, instead of the dimension $n$ of our projective space, it is rather the quantity
\begin{equation} \label{def_q}
q := \frac{n+1}{2^{d-2}(d-1)} > 4
\end{equation}
that will naturally arise and play a key role in the bounds. The inequality $q>4$ is equivalent to the assumption $n \geq 2^d(d-1)$ we made.

\begin{lemma} \label{sing_cv}
The formal series $\sum_{\delta=0}^\infty \left( \sum_{\nu\vdash\delta} S_\nu \right) T^\delta$ converges at $T=\mathbf{L}^{-(n+1)}$. We call its value the \emph{singular series} and denote it
\[
\mathfrak{S}^\mu := \sum_{\delta=0}^\infty \mathbf{L}^{-(n+1)\delta} \sum_{\nu\vdash\delta} S_\nu \in \widehat{\mathscr{E}xp\mathscr{M}}_{\calB_\mu} .
\]
Furthermore the remainder is bounded in weight as follows:
\[
w_{\calB_\mu}\left( \mathfrak{S}^\mu - \sum_{\delta=0}^{\eta_e} \mathbf{L}^{-(n+1)\delta} \sum_{\nu\vdash\delta} S_\nu \right) \leq -(q-4)e + 2(q-4)(g-1) + ((d+1)q-4)l' .
\]
\end{lemma}
Note that $\mathfrak{S}^\mu$ depends on $W$ and $\mu$ but not on $e$.

\begin{proof}
We want to bound the weight of the degree $\delta$ term, $\sum_{\nu\vdash\delta} S_\nu$. By definition (\ref{def_S}) this rewrites as a finite sum of $\mathbf{L}^{(n+1)\left(\deg(\scrS \cap X) + |\nu^0|\right)} \cdot \left[J(\scrC) \times U_\nu(X,\nu^0)\right] \cdot S'_{\nu;X,\nu^0}$: it suffices to bound each of those terms. By the second point of proposition \ref{prop_weight},
\begin{align*}
& w_{\calB_\mu}\left(\left[J(\scrC) \times U_\nu(X,\nu^0)\right] \cdot S'_{\nu;X,\nu^0}\right) \\
& \leq 2\dim_{\calB_\mu}\left( J(\scrC) \times U_\nu(X,\nu^0) \times P_{|Z}^{\vee*} \times \tilde W \right) \\
& \quad + w_{J(\scrC) \times U_\nu(X,\nu^0) \times P_{|Z}^{\vee*} \times \tilde W}\left(\left[J(\scrC) \times U_\nu(X,\nu^0)\right] \cdot S'_{\nu;X,\nu^0}\right) .
\end{align*}
Going back to the definition (\ref{def_Bmu}) of $\calB_\mu$, the first term is
\begin{align*}
2 \dim_{\Sym^\mu(\CmS)}\left( U_\nu(X,\nu^0) \times P_{|Z}^{\vee*} \right) & = 2\left( \dim_{\Sym^\mu(\CmS)} U_\nu(X,\nu^0) + \dim P_{|Z}^{\vee*} \right) \\
& \leq 2\left( \dim \Sym^\nu(\scrC) + \dim P_{|Z}^\vee \right) \\
& = 2\left( |\nu| + \delta \right) \\
& \leq 4\delta .
\end{align*}
For the second term, with the definition (\ref{def_S'}) of $S'_{\nu;X,\nu^0}$ in mind, one can apply our local version of Weyl's inequality: proposition \ref{concl_Weyl_dif_local} provides the bound
\[
2^{-(d-2)}(n+1) \left(s + \left(2^{d-1}-1\right)(\delta-\delta')\right)
\]
with $\delta' = \deg(\scrS \cap X)+|\nu^0|$ and $s = \left\lfloor \frac{(d-2)\delta+\delta'}{d-1} \right\rfloor$. When putting everything together, including the normalisation, we get
\begin{align*}
& w_{\calB_\mu}\left(\mathbf{L}^{(n+1)\delta'} \cdot \left[J(\scrC) \times U_\nu(X,\nu^0)\right] \cdot S'_{\nu;X,\nu^0}\right) \\
& \leq 2(n+1)\delta' + 4\delta + 2^{-(d-2)}(n+1) \left(s + \left(2^{d-1}-1\right)(\delta-\delta')\right) \\
& = 2(n+1)\delta + 4\delta + 2^{-(d-2)}(n+1) \left(s - (\delta-\delta')\right) \\
& \leq 2(n+1)\delta + 4\delta + 2^{-(d-2)}(n+1) \left(\frac{(d-2)\delta+\delta'}{d-1} - \delta + \delta'\right) \\
& = 2(n+1)\delta + 4\delta + \frac{n+1}{2^{d-2}(d-1)}(-\delta + d\delta') \\
& = 2(n+1)\delta - (q-4)\delta + dq\delta' .
\end{align*}
Remark that $\delta' \leq \deg\scrS + |\mu| = l'$: this shows that
\begin{equation} \label{weight_S}
w_{\calB_\mu}\left( \sum_{\nu\vdash\delta} S_\nu \right) \leq 2(n+1)\delta - (q-4)\delta + dql' .
\end{equation}

By definition (\ref{def_rad_cv}) the formal series we are considering has radius of convergence $r = \limsup_{\delta\to\infty} \frac{1}{2\delta}w_{\calB_\mu}\left( \sum_{\nu\vdash\delta} S_\nu \right)$. The previous computations imply $r \leq n+1 - \frac{q-4}{2}$, in particular $r < n+1$ since $q>4$. After lemma \ref{rad_cv}, the series converges at $T=\mathbf{L}^{-(n+1)}$, which lets us define $\mathfrak{S}^\mu$.

As for the remainder of index $\eta_e$, it rewrites as
\[
\sum_{\delta=\eta_e+1}^\infty \mathbf{L}^{-(n+1)\delta} \sum_{\nu\vdash\delta} S_\nu .
\]
By (\ref{weight_S}), the term of index $\delta$ is bounded in weight by $-(q-4)\delta + dql'$. This is decreasing in $\delta$ so the weight of the whole remainder is at most this bound applied at $\delta = \eta_e+1$:
\begin{align*}
-(q-4)(\eta_e+1) + dql' & = -(q-4)(e-l'-2g+2) + dql' \\
& = -(q-4)e + 2(q-4)(g-1) + ((d+1)q-4)l' .
\end{align*}
\end{proof}

\subsubsection{Writing the singular series as a motivic Euler product}
\label{subsubsect_sing_MEP}

Our goal is then to express the singular series
\[
\mathfrak{S}^\mu = \sum_{\delta=0}^\infty \mathbf{L}^{-(n+1)\delta} \sum_{\nu\vdash\delta} S_\nu
\]
as a motivic Euler product indexed by $\scrC$. Considering first its formal version $\sum_{\delta=0}^\infty \left( \sum_{\nu\vdash\delta} S_\nu \right) T^\delta$, by definition \ref{def_MEP} we seek to establish that for all $\nu$, $S_\nu$ is the configuration space of a certain family of elements of $K_0\left(\mathrm{ExpVar}_{\calB_\mu \times \scrC}\right)$. Not delving directly into these considerations, we start with the following lemma, analogous to \cite{BB25}'s proposition 6.7.

\begin{lemma}[factorisation] \label{lem_fact}
Take $\nu_1, \nu_2$ two integer partitions and consider the mixed product $\Sym^{\nu_1,\nu_2}(\scrC)$ defined in (\ref{def_mixed_sym}). It comes with the two projections $\pi_1 : \Sym^{\nu_1,\nu_2}(\scrC) \to \Sym^{\nu_1}(\scrC)$ and $\pi_2 : \Sym^{\nu_1,\nu_2}(\scrC) \to \Sym^{\nu_2}(\scrC)$ and the sum morphism
\[
\tau :
\left\{\begin{array}{rcl}
\Sym^{\nu_1,\nu_2}(\scrC) & \longrightarrow & \Sym^{\nu_1+\nu_2}(\scrC) \\
(Z_1, Z_2) & \longmapsto & Z_1+Z_2
\end{array}\right. .
\]
Then in $K_0\left(\mathrm{ExpVar}_{\calB_\mu \times \Sym^{\nu_1,\nu_2}(\scrC)}\right)$ we have
\[
\tau^* S_{\nu_1+\nu_2} = \pi_1^* S_{\nu_1} \cdot \pi_2^* S_{\nu_2} .
\]
Moreover in the case $\nu_1 = \nu_2$, this identity commutes with the automorphism of $\Sym^{\nu_1,\nu_2}(\scrC)$ that exchanges $Z_1$ and $Z_2$.
\end{lemma}

\begin{proof}
We start on the right-hand side. Recall definition (\ref{def_S}):
\[
S_{\nu_1} = \sum_{(X_1, \nu_1^0) \in I_{\nu_1}} \mathbf{L}^{(n+1)\left(\deg(\scrS \cap X_1) + |\nu_1^0|\right)} \cdot \left[J(\scrC) \times U_{\nu_1}(X_1,\nu_1^0)\right] \cdot S'_{\nu_1 ; X_1,\nu_1^0}
\]
with $S'_{\nu_1;X_1,\nu_1^0} \in K_0\left(\mathrm{ExpVar}_{\tilde W}\right) $. We have
\begin{align*}
& \pi_1^*\left(\left[J(\scrC) \times U_{\nu_1}(X_1,\nu_1^0)\right] \cdot S'_{\nu_1;X_1,\nu_1^0}\right) \\
& = \left[J(\scrC) \times \left( U_{\nu_1}(X_1, \nu_1^0) \times_{\Sym^{\nu_1}(\scrC)} \Sym^{\nu_1,\nu_2}(\scrC) \right)\right] \cdot S'_{\nu_1;X_1,\nu_1^0} ,
\end{align*}
therefore $\pi_1^* S_{\nu_1}$ is the sum of these, weighted by the powers of $\mathbf{L}$ above. We write things analogously for $\nu_2$: putting them together we have
\begin{align} \label{lem_factor_sum1}
\pi_1^* S_{\nu_1} \cdot \pi_2^* S_{\nu_2} = \sum_{\substack{(X_1, \nu_1^0) \in I_{\nu_1} \\ (X_2, \nu_2^0) \in I_{\nu_2}}} \mathbf{L}^{(n+1)\delta'} \cdot \left[ J(\scrC) \times U_{\nu_1,\nu_2}(X_1,\nu_1^0 ; X_2,\nu_2^0) \right] \cdot S'_{\nu_1;X_1,\nu_1^0} \cdot S'_{\nu_2;X_2,\nu_2^0}
\end{align}
where $\delta' = \deg(\scrS \cap X_1) + \deg(\scrS \cap X_2) + |\nu_1^0| + |\nu_2^0|$ and
\begin{align*}
& U_{\nu_1,\nu_2}(X_1,\nu_1^0 ; X_2,\nu_2^0) \\
& := \left( U_{\nu_1}(X_1, \nu_1^0) \times_{\Sym^{\nu_1}(\scrC)} \Sym^{\nu_1,\nu_2}(\scrC) \right) \\
& \qquad \times_{\Sym^\mu(\CmS) \times \Sym^{\nu_1,\nu_2}(\scrC)} \left( U_{\nu_2}(X_2, \nu_2^0) \times_{\Sym^{\nu_2}(\scrC)} \Sym^{\nu_1,\nu_2}(\scrC) \right) \\
& = \left\{ (\scrT, (Z_1,Z_2)) \in \Sym^\mu(\CmS) \times \Sym^{\nu_1,\nu_2}(\scrC) \,\left|\,
\begin{array}{l}
(\scrT,Z_1) \in U_{\nu_1}(X_1, \nu_1^0) \\
(\scrT,Z_2) \in U_{\nu_2}(X_2, \nu_2^0)
\end{array} \right.\right\} .
\end{align*}
If $(\scrT, (Z_1,Z_2))$ is in the set above, remark that the condition $Z_1 \cap Z_2 = \varnothing$ imposes $X_1 \cap X_2 = \varnothing$. Similarly, $Z_1|_{|\scrT|}$ and $Z_2|_{|\scrT|}$ are disjoint and have respective partitions $\nu_1^0$ and $\nu_2^0$, so $|\scrT|$ contains at least $|\nu_1^0| + |\nu_2^0|$ points: $|\nu_1^0| + |\nu_2^0| \leq |\mu|$. This shows that if one of these two conditions is not met, $U_{\nu_1,\nu_2}(X_1,\nu_1^0 ; X_2,\nu_2^0)$ is empty, and the corresponding term in the sum (\ref{lem_factor_sum1}) is zero. As a consequence that sum can be indexed by
\begin{equation} \label{def_I2}
I_{\nu_1,\nu_2} := \left\{ ((X_1,\nu_1^0), (X_2,\nu_2^0)) \in I_{\nu_1} \times I_{\nu_2} \,\left|\,
\begin{array}{l}
X_1 \cap X_2 = \varnothing \\
|\nu_1^0| + |\nu_2^0| \leq |\mu|
\end{array} \right.\right\} .
\end{equation}
In particular in this case, $\deg(\scrS \cap X_1) + \deg(\scrS \cap X_2) = \deg(\scrS \cap (X_1+X_2))$, so
\begin{align} \label{lem_factor_sum2}
&\pi_1^* S_{\nu_1} \cdot \pi_2^* S_{\nu_2} \nonumber \\
& = \sum_{((X_1,\nu_1^0), (X_2,\nu_2^0)) \in I_{\nu_1,\nu_2}} \mathbf{L}^{(n+1)\delta'} \cdot \left[ J(\scrC) \times U_{\nu_1,\nu_2}(X_1,\nu_1^0 ; X_2,\nu_2^0) \right] \cdot S'_{\nu_1;X_1,\nu_1^0} \cdot S'_{\nu_2;X_2,\nu_2^0}
\end{align}
with $\delta' = \deg(\scrS \cap (X_1+X_2)) + |\nu_1^0 + \nu_2^0|$.

Now let us look at the side with $\tau$. By definition (\ref{def_S}) again, we have
\[
\tau^* S_{\nu_1+\nu_2} = \sum_{(X, \nu^0) \in I_{\nu_1+\nu_2}} \mathbf{L}^{(n+1)\left(\deg(\scrS \cap X) + |\nu^0|\right)} \cdot \left[J(\scrC) \times U' \right] \cdot S'_{\nu_1+\nu_2 ; X,\nu^0}
\]
where
\begin{align*}
U' & = \left( U_{\nu_1+\nu_2}(X,\nu^0) \times_{\Sym^{\nu_1+\nu_2}(\scrC)} \Sym^{\nu_1,\nu_2}(\scrC) \right) \\
& \simeq \bigsqcup_{\substack{((X_1,\nu_1^0), (X_2,\nu_2^0)) \in I_{\nu_1,\nu_2} \\ X_1+X_2 = X \\ \nu_1^0+\nu_2^0 = \nu^0}} U_{\nu_1,\nu_2}(X_1,\nu_1^0 ; X_2,\nu_2^0) .
\end{align*}
The justification for this stratification is that in the first line we have $(\scrT,Z) \in U_{\nu_1+\nu_2}(X,\nu^0)$ with $Z=Z_1+Z_2$, which amounts to taking $(\scrT, (Z_1,Z_2)) \in U_{\nu_1,\nu_2}(X_1,\nu_1^0 ; X_2,\nu_2^0)$ for some decompositions $X=X_1+X_2$ and $\nu^0=\nu_1^0+\nu_2^0$ (recall the definitions of $I_\nu$ above (\ref{def_Unu}) and of $I_{\nu_1,\nu_2}$ in (\ref{def_I2})). As a consequence, $\tau^* S_{\nu_1+\nu_2}$ rewrites as
\begin{align*}
& \sum_{(X, \nu^0) \in I_{\nu_1+\nu_2}} \ \sum_{\substack{((X_1,\nu_1^0), (X_2,\nu_2^0)) \in I_{\nu_1,\nu_2} \\ X_1+X_2 = X \\ \nu_1^0+\nu_2^0 = \nu^0}} \mathbf{L}^{(n+1)\delta'} \cdot \left[J(\scrC) \times U_{\nu_1,\nu_2}(X_1,\nu_1^0 ; X_2,\nu_2^0) \right] \cdot S'_{\nu_1+\nu_2 ; X_1+X_2 ,\nu_1^0+\nu_2^0} \\
& = \sum_{((X_1,\nu_1^0), (X_2,\nu_2^0)) \in I_{\nu_1,\nu_2}} \mathbf{L}^{(n+1)\delta'} \cdot \left[J(\scrC) \times U_{\nu_1,\nu_2}(X_1,\nu_1^0 ; X_2,\nu_2^0) \right] \cdot S'_{\nu_1+\nu_2 ; X_1+X_2 ,\nu_1^0+\nu_2^0}
\end{align*}
with
\[
\delta' = \deg(\scrS \cap X) + |\nu^0| = \deg(\scrS \cap (X_1 + X_2)) + |\nu_1^0 + \nu_2^0| .
\]
Compare this to (\ref{lem_factor_sum2}): showing that $S'_{\nu_1;X_1,\nu_1^0} \cdot S'_{\nu_2;X_2,\nu_2^0} = S'_{\nu_1+\nu_2 ; X_1+X_2 ,\nu_1^0+\nu_2^0}$ over $\tilde W$ would imply the lemma. Let us do this.

By lemma \ref{Gamma_S'} and the definition of $U_{\nu_1,\nu_2}(X_1,\nu_1^0 ; X_2,\nu_2^0)$ above, one can write over $\tilde W$, with some closed
\[
(\scrT, (Z_1,Z_2)) \in U_{\nu_1,\nu_2}(X_1,\nu_1^0 ; X_2,\nu_2^0) ,
\]
\begin{align*}
S'_{\nu_1;X_1,\nu_1^0} \cdot S'_{\nu_2;X_2,\nu_2^0} & = \left[ P_{|Z_1}^{\vee*} \times P_{-\scrST|Z_1}^{n+1} \times \tilde W ,\, \beta_1\left(f\left( \mathbf{r}_1 + \sigma_{|Z_1}^{\scrT}(\mathbf{w}) \right)\right) \right] \\
&\qquad \cdot \left[ P_{|Z_2}^{\vee*} \times P_{-\scrST|Z_2}^{n+1} \times \tilde W ,\, \beta_2\left(f\left( \mathbf{r}_2 + \sigma_{|Z_2}^{\scrT}(\mathbf{w}) \right)\right) \right] \\
& = \Big[ \left(P_{|Z_1}^{\vee*} \times P_{|Z_2}^{\vee*}\right) \times \left(P_{-\scrST|Z_1}^{n+1} \times P_{-\scrST|Z_2}^{n+1} \right) \times \tilde W ,\, \\
&\qquad \beta_1\left(f\left( \mathbf{r}_1 + \sigma_{|Z_1}^{\scrT}(\mathbf{w}) \right)\right) + \beta_2\left(f\left( \mathbf{r}_2 + \sigma_{|Z_2}^{\scrT}(\mathbf{w}) \right)\right) \Big] .
\end{align*}
Since we took $Z_1 \cap Z_2 = \varnothing$ there is an obvious isomorphism of $\C$-algebras
\[
\varphi : P_{|Z_1} \times P_{|Z_2} \overset{\sim}{\longrightarrow} P_{|Z_1+Z_2},
\]
whose inverse is $(\cdot|_{Z_1}, \cdot|_{Z_2})$. It restricts to $P_{-\scrST|Z_1} \times P_{-\scrST|Z_2} \overset{\sim}{\to} P_{-\scrST|Z_1+Z_2}$ (in fact recall that $P_{-\scrST|Z_1} = P_{-\scrST \cap Z_1|Z_1}$ and same for the two other spaces). The dual isomorphism of $\varphi$ is
\[\psi :
\left\{\begin{array}{rcl}
P_{|Z_1}^\vee \times P_{|Z_2}^\vee & \overset{\sim}{\longrightarrow} & P_{|Z_1+Z_2}^\vee \\
(\beta_1, \beta_2) & \longmapsto & \beta_1(\cdot|_{Z_1}) + \beta_2(\cdot|_{Z_2})
\end{array}\right. ,
\]
which restricts to $P_{|Z_1}^{\vee*} \times P_{|Z_2}^{\vee*} \overset{\sim}{\to} P_{|Z_1+Z_2}^{\vee*}$. Indeed write $\gamma = \psi(\beta_1, \beta_2)$: if $\beta_1$ factors through $Z'_1 \subsetneq Z_1$ then $\gamma$ factors through $Z'_1+Z_2 \subsetneq Z_1+Z_2$ (and same for $\beta_2$). Conversely if $\gamma$ factors through a $Z' \subset Z_1+Z_2$, $\beta_i$ factors through $Z_i \cap Z'$ for $i\in\{1,2\}$ and at least of one the $Z_i \cap Z' \subset Z_i$ is strict. We want to make $\varphi$ (extended to $(n+1)$-tuples) and $\psi$ appear in the exponential. For $i\in\{1,2\}$, $\mathbf{r}_i = \varphi(\mathbf{r}_1, \mathbf{r}_2)|_{Z_i}$. We wish to also replace $\sigma_{|Z_i}^{\scrT}(\mathbf{w})$ with $\sigma_{|Z_1+Z_2}^{\scrT}(\mathbf{w})\big|_{Z_i}$: just like in the proof of lemma \ref{Gamma_S'}, they need not be equal, but we can still do so by noting that the translation
\[
\mathbf{r}_i \longmapsto \mathbf{r}_i + \sigma_{|Z_i}^{\scrT}(\mathbf{w}) - \sigma_{|Z_1+Z_2}^{\scrT}(\mathbf{w})\big|_{Z_i}
\]
is an automorphism of $P_{-\scrST|Z_i}^{n+1}$. The class $S'_{\nu_1;X_1,\nu_1^0} \cdot S'_{\nu_2;X_2,\nu_2^0}$ can therefore be written with the exponential
\begin{align*}
& \beta_1\left(f\left( \mathbf{r}_1 + \sigma_{|Z_1+Z_2}^{\scrT}(\mathbf{w})\big|_{Z_1} \right)\right) + \beta_2\left(f\left( \mathbf{r}_2 + \sigma_{|Z_1+Z_2}^{\scrT}(\mathbf{w})\big|_{Z_2} \right)\right) \\
& =  \beta_1\left(\left.f\left( \varphi(\mathbf{r}_1, \mathbf{r}_2) + \sigma_{|Z_1+Z_2}^{\scrT}(\mathbf{w}) \right)\right|_{Z_1}\right) + \beta_2\left(\left.f\left( \varphi(\mathbf{r}_1, \mathbf{r}_2) + \sigma_{|Z_1+Z_2}^{\scrT}(\mathbf{w}) \right)\right|_{Z_2}\right) \\
& = \psi(\beta_1, \beta_2)\left( f\left( \varphi(\mathbf{r}_1, \mathbf{r}_2) + \sigma_{|Z_1+Z_2}^{\scrT}(\mathbf{w}) \right)\right) .
\end{align*}
This justifies that
\begin{align*} \label{lem_factor_pi}
S'_{\nu_1;X_1,\nu_1^0} \cdot S'_{\nu_2;X_2,\nu_2^0} & = \left[ P_{|Z_1+Z_2}^{\vee*} \times P_{-\scrST|Z_1+Z_2}^{n+1} \times \tilde W ,\, \gamma\left( f\left( \mathbf{q} + \sigma_{|Z_1+Z_2}^{\scrT}(\mathbf{w}) \right)\right) \right] \nonumber \\
& = S'_{\nu_1+\nu_2 ; X_1+X_2 ,\nu_1^0+\nu_2^0} .
\end{align*}
\end{proof}

For any integer $i\geq1$ recall that in subsection \ref{part_conf} we denoted by $[i]$ the trivial partition of $i$ that writes as $i=i$, such that $\Sym^{[i]}(\scrC) = \{ ix \,;\, x\in\scrC \} \simeq \scrC$. In particular we consider $S_{[i]} \in K_0\left(\mathrm{ExpVar}_{\calB_\mu \times \scrC}\right)$, that will soon be described in lemma \ref{S_i}. For $x\in\scrC$ we denote
\[
S_{[i],x} := S_{[i]}(x) \in K_0\left(\mathrm{ExpVar}_{\calB_\mu \times \kappa(x)}\right)
\]
its motivic evaluation in $x$.

\begin{proposition} \label{sing_mep_1}
In $K_0\left(\mathrm{ExpVar}_{\calB_\mu}\right)\llbracket T \rrbracket$ we have
\[
\sum_{\delta=0}^\infty \left( \sum_{\nu\vdash\delta} S_\nu \right) T^\delta = \prod_{x\in\scrC} \left( 1 + \sum_{i\geq1} S_{[i],x} T^i \right) .
\]
In particular the singular series can be written as a motivic Euler product in $\widehat{\mathscr{E}xp\mathscr{M}}_{\calB_\mu}$:
\[
\mathfrak{S}^\mu = \prod_{x\in\scrC} \left( 1 + \sum_{i\geq1} S_{[i],x} T^i \right)_{|T = \mathbf{L}^{-(n+1)}} .
\]
\end{proposition}
Technically the motivic Euler products should be indexed by $x \in (\calB_\mu \times \scrC) / \calB_\mu$, but thanks to remark \ref{MEP_short_index} we only write $x\in\scrC$. This shortening will be used until the proof of proposition \ref{sing_prod_final}.

\begin{proof}
The second part is the evaluation of the first one at $\mathbf{L}^{-(n+1)}$, which is possible by the previous paragraph on convergence. For the first identity, by definition of motivic Euler products, the goal is to show that
\[
S_\nu = \Conf_{\calB_\mu \times \scrC / \calB_\mu}^\nu \left((S_{[i]})_{i\geq1}\right)
\]
for every integer partition $\nu = (\nu_i)_{i\geq1}$. First of all let us introduce generalized versions of the objects appearing in the factorisation lemma \ref{lem_fact}: define
\[
Q_\nu \subset \prod_{i\geq1} \left(\Sym^{[i]}(\scrC)\right)^{\nu_i}
\]
as the subvariety made up of the $(ix_{i,j})_{i,j} = (ix_{i,1}, \dots, ix_{i,\nu_i})_{i\geq1}$ such that $x_{i,j} \neq x_{i',j'}$ for all $(i,j)\neq(i',j')$. Set
\[
\pi_{i',j'} :
\left\{\begin{array}{rcl}
Q_\nu & \longrightarrow & \Sym^{[i']}(\scrC) \\
(ix_{i,j})_{i,j} & \longmapsto & i'x_{i',j'}
\end{array}\right.
\]
for all $i\geq1$ and $1 \leq j \leq \nu_i$, as well as
\[
\tau :
\left\{\begin{array}{rcl}
Q_\nu & \longrightarrow & \Sym^\nu(\scrC) \\
(ix_{i,j})_{i,j} & \longmapsto & \sum_{i,j} ix_{i,j}
\end{array}\right. .
\]
$\tau$ is well defined: since the $x_{i,j}$ are pairwise distinct, $\sum_{i,j} ix_{i,j}$ has, for any $i\geq1$, exactly $\nu_i$ points of multiplicity $i$, which means it lies in $\Sym^\nu(\scrC)$. With these notations, iterating the factorisation lemma gives
\begin{equation} \label{sing_mep_1_intermed}
\tau^* S_\nu = \prod_{i,j} \pi_{i,j}^* S_{[i]}
\end{equation}
over $\calB_\mu \times Q_\nu$. Note that this product is equivalent to taking the product of the $S_{[i]}$ over $\calB_\mu$ and to request that they map to pairwise distinct $x \in \scrC$; recalling definition (\ref{def_conf}) we start to see the link with configuration spaces. The symmetric group $\mathfrak{S}_\nu = \prod_{i\geq1} \mathfrak{S}_{\nu_i}$ acts on $Q_\nu$, and $\tau$ induces an isomorphism $\bar\tau : Q_\nu/\mathfrak{S}_\nu \overset{\sim}{\to} \Sym^\nu(\scrC)$. Since lemma \ref{lem_fact} commutes with switching components of same partition, which amounts in this generalized context to permuting the $ix_{i,1}, \dots, ix_{i,\nu_i}$ for any fixed $i$, we can quotient (\ref{sing_mep_1_intermed}) by $\mathfrak{S}_\nu$ to obtain
\[
\bar\tau^* S_\nu = \Conf_{\calB_\mu \times \scrC / \calB_\mu}^\nu \left((S_{[i]})_{i\geq1}\right)
\]
over $\calB_\mu \times \Sym^\nu(\scrC)$. Since $\bar\tau$ is an isomorphism one can drop the $\bar\tau^*$, which concludes.
\end{proof}

We can complete this by showing that in the motivic Euler product obtained for $\mathfrak{S}^\mu$, the evaluation in $T = \mathbf{L}^{-(n+1)}$ commutes with the product. Thanks to the weight bound established for $S_\nu$ in the proof of lemma \ref{sing_cv}, the formal series $\sum_{i\geq0} S_{[i]} T^i$ converges at $T=\mathbf{L}^{-(n+1)}$ in $\widehat{\mathscr{E}xp\mathscr{M}}_{\calB_\mu \times \scrC}$. The result below follows.
\begin{corollary} \label{sing_mep_2}
With the notation (\ref{MEP_noT}) we can write
\[
\mathfrak{S}^\mu = \prod_{x\in\scrC} \left( 1 + \sum_{i\geq1} \mathbf{L}^{-(n+1)i} S_{[i],x} \right) .
\]
in $\widehat{\mathscr{E}xp\mathscr{M}}_{\calB_\mu}$.
\end{corollary}

\begin{proof}
We introduce the companion motivic Euler product, with an extra variable $T_0$,
\[
F(T,T_0) := \prod_{x\in\scrC} \left( 1 + \sum_{i\geq1} S_{[i],x} T^i T_0 \right) \in K_0\left(\mathrm{ExpVar}_{\calB_\mu}\right)\llbracket T, T_0 \rrbracket .
\]
Proposition \ref{sing_mep_1} states that $\mathfrak{S}^\mu = \left( F(T, T_0)_{|T_0=1} \right)_{|T=\mathbf{L}^{-(n+1)}}$, where the order of evaluation matters. By the remark on convergence just above, $F(T,T_0)$ can be evaluated at $T=\mathbf{L}^{-(n+1)}$ first, which means the two evaluations can take place in any order. Thus
\begin{align*}
\mathfrak{S}^\mu & = \left( F(T, T_0)_{|T=\mathbf{L}^{-(n+1)}} \right)_{|T_0=1} \\
& = \prod_{x\in\scrC} \left( 1 + \sum_{i\geq1} \mathbf{L}^{-(n+1)i}S_{[i],x} T_0 \right)_{|T_0=1} ,
\end{align*}
which is precisely the wanted motivic Euler product by notation (\ref{MEP_noT}).
\end{proof}

\subsubsection{Local densities} \label{locdens}

In this paragraph we are going to express the factors of the previous motivic Euler product as \emph{local densities} (like \cite{BB25} in their paragraph 6.3.4). The first step for us is to express $S_{[i]}$ properly. Before that for $p\in|\scrS|$ we define a new section map, in order to replace $\sigma_{|ip}^{\scrT}(\mathbf{w})$, as follows. Recall that $l_p$ is the multiplicity of $p$ in $\scrS$. For $N \geq l_p$ a special case of local Euclidean division exact sequence (\ref{exseq_ed_loc}) is
\[
\xymatrix{
0 \ar[r] & P_{-l_p p|Np} \ar[r]^-{\subset} & P_{|Np} \ar[r]^-{\cdot|_{l_p p}} & P_{|l_p p} \ar[r] & 0
} :
\]
it admits a section
\begin{equation} \label{def_sect_sigma_3}
\sigma_{Np} : P_{|l_p p} \longrightarrow P_{|Np} .
\end{equation}
Next we denote
\begin{equation} \label{def_U}
U := \Sym^{\mu,[1]}(\CmS)
\end{equation}
and $U^0$ its complement in
\[
\Sym^\mu(\CmS) \times (\CmS)
\]
(i.e. the $(\scrT,x)$ such that $x\in\scrT$), as well as
\[
U_p^\scrS := \Sym^\mu(\CmS) \times \kappa(p)
\]
for $p\in|\scrS|$. In the end we have a partition
\[
\Sym^\mu(\CmS) \times \scrC = U \sqcup U^0 \sqcup \bigsqcup_{p\in|\scrS|} U_p^\scrS .
\]
The reason for these shorthands will become clear in the proof below; note that~$i$ doesn't appear in them. For $p\in|\scrS|$ we also put $j_{ip} := \min(l_p, i)$.

\begin{lemma} \label{S_i}
Over $\calB_\mu \times \scrC$ we have, for any closed $x\in\CmS$ and $N \geq \max(l, i)$,
\begin{align*}
S_{[i]} = & \left[ J(\scrC) \times U\right] \cdot \left[ P_{|ix}^{\vee*} \times P_{|ix}^{n+1} \times \tilde W ,\, \beta\left(f\left( \mathbf{r} \right)\right) \right] \\
& + \left[ J(\scrC) \times U^0\right] \cdot \mathbf{L}^{n+1} \cdot \left[ P_{|ix}^{\vee*} \times P_{-x|ix}^{n+1} \times \tilde W ,\, \beta\left(f\left( \mathbf{r} \right)\right) \right] \\
& + \sum_{p\in|\scrS|} \left[J(\scrC) \times U_p^\scrS\right] \cdot \mathbf{L}^{(n+1)j_{ip}} \cdot \left[ P_{|ip}^{\vee*} \times P_{-j_{ip}p|ip}^{n+1} \times \tilde W ,\, \beta\left(f\left( \mathbf{r} + \sigma_{Np}(\mathbf{w}_p)|_{ip} \right)\right) \right]
\end{align*}
(the classes don't depend on the chosen $x$).
\end{lemma}

\begin{proof}
By definition (\ref{def_S}) we have
\[
S_{[i]} = \sum_{(X, \nu^0) \in I_{[i]}} \mathbf{L}^{(n+1)\left(\deg(\scrS \cap X) + |\nu^0|\right)} \cdot \left[J(\scrC) \times U_{[i]}(X,\nu^0)\right] \cdot S'_{[i];X,\nu^0} .
\]
The only subpartitions of $[i]$ are itself and $\varnothing$. Recall as well the definition of $I_\nu$ above (\ref{def_Unu}): in our case its elements are $(X,\nu^0)$ with first $\pt X \leq [i]$ and $|X| \subset |\scrS|$, i.e. either $\varnothing$ or $ip$ with $p \in |\scrS|$; and second $\nu^0 \leq [i] - \pt X$ and $|\nu^0| \leq |\mu|$, i.e. either $\varnothing$ or, as long as $X=\varnothing$ and $\mu\neq\varnothing$, $[i]$. Let's look at the terms corresponding to these three respective cases.
\begin{itemize}

\item $(\varnothing,\varnothing)$. First, $\deg(\scrS \cap \varnothing)+|\varnothing| = 0$ and
\[
U_{[i]}(\varnothing,\varnothing) = \Sym^{\mu,[i]}(\CmS) \simeq U .
\]
As a consequence, going back to the definition (\ref{def_S'}) of the $S'_{\nu;X,\nu^0}$, for any closed $(\scrT,x) \in U$ we have
\[
S'_{[i];\varnothing,\varnothing} = \left[ P_{|ix}^{\vee*} \times P_{|ix}^{n+1} \times \tilde W ,\, \beta\left(f\left( \mathbf{r} + \sigma_{|ix}^{\scrT}(\mathbf{w}) \right)\right) \right]
\]
because $x$ does not meet $\scrST$. Since $\sigma_{|ix}^{\scrT}(\mathbf{w}) \in P_{|ix}^{n+1}$, we can remove it, by invoking a translation in $P_{|ix}^{n+1}$, as in the proof of lemma \ref{maj_intermed_2} and several other times. We obtain the class
\[
S'_{[i];\varnothing,\varnothing} = \left[ P_{|ix}^{\vee*} \times P_{|ix}^{n+1} \times \tilde W ,\, \beta\left(f\left( \mathbf{r} \right)\right) \right] ,
\]
which doesn't feature $\scrT$ anymore: as a consequence this equality holds for any closed $x\in\CmS$.

\item $(\varnothing, [i])$. Here $\deg(\scrS \cap \varnothing)+|[i]| = 1$ and $U_{[i]}(\varnothing,[i])$ corresponds to the $(\scrT,ix)$ supported by $\CmS$ with $x \in \scrT$: in other words it is isomorphic (by forgetting the $i$ in front of~$x$) to $U^0$. If $\mu\neq\varnothing$, $U^0$ is nonempty and for $(\scrT,x) \in U^0$ closed, since $x\in\scrT$ we have
\[
S'_{[i];\varnothing,[i]} = \left[ P_{|ix}^{\vee*} \times P_{-x|ix}^{n+1} \times \tilde W ,\, \beta\left(f\left( \mathbf{r} + \sigma_{|ix}^{\scrT}(\mathbf{w}) \right)\right) \right] .
\]
This time $\sigma_{|ix}^{\scrT}(\mathbf{w})$ restricts to $0$ over $x$ so it belongs to $P_{-x|ix}^{n+1}$ and we can again remove it; the resulting class is the same for any closed $x\in\CmS$. In the case $\mu=\varnothing$, $U^0$ is empty and the second line in the formula of the lemma is zero anyway: we can put the class we want on its right.

\item $(ip,\varnothing), p\in|\scrS|$. For this last one $\deg(\scrS \cap i)+|\varnothing| = \min(l_p,i) = j_{ip}$. We also have
\[
U_{[i]}(ip,\varnothing) = \Sym^\mu(\CmS) \times \{ip\} \simeq U_p^\scrS .
\]
Then recall that $P_{-Z_2|Z_1} = P_{-(Z_1 \cap Z_2)|Z_1}$: for $\scrT \in \Sym^\mu(\CmS)$ closed,
\[
S'_{[i];ip,\varnothing} = \left[ P_{|ip}^{\vee*} \times P_{-j_{ip}p|ip}^{n+1} \times \tilde W ,\, \beta\left(f\left( \mathbf{r} + \sigma_{|ip}^{\scrT}(\mathbf{w}) \right)\right) \right] .
\]
By a translation in $P_{-j_{ip}p|ip}^{n+1}$ we replace $\sigma_{|ip}^{\scrT}(\mathbf{w})$ with $\sigma_{Np}(\mathbf{w}_p)|_{ip}$, which is possible because $N \geq \max(l_p,i)$ and by construction both restrict to $\mathbf{w}_p|_{j_{ip}p}$.

\end{itemize}
\end{proof}

Now, with a view towards the local densities we aim for, let us recall that we defined
\[
\Lambda_N(x) := \left\{\left. \mathbf{u} \in P_{|Nx}^{n+1} \,\right|\, f(\mathbf{u}) = 0 \right\}
\]
in (\ref{def_Lambda}), for any closed $x\in\scrC$ and $N\in\mathbb{N}$. We also set, for $N\geq1$,
\begin{equation} \label{def_Lambda0}
\Lambda_N^0(x) := \left\{\left. \mathbf{u} \in P_{-x|Nx}^{n+1} \,\right|\, f(\mathbf{u}) = 0 \right\} \subset \Lambda_N(x) .
\end{equation}
For any $x,y \in \scrC$, $\Lambda_N(x) \simeq \Lambda_N(y)$ and $\Lambda_N^0(x) \simeq \Lambda_N^0(y)$. Finally if $x = p \in |\scrS|$ and $N \geq l_p$ we define
\[
\Lambda_N^W(p) := \left\{ \mathbf{u} \in P_{|Np}^{n+1} \,\left|\,
\begin{array}{l}
f(\mathbf{u}) = 0 \\
\mathbf{u}|_{l_p p} \in \tilde W_p
\end{array} \right.\right\} \subset \Lambda_N(p) .
\]
Like we did in (\ref{reformulation_M}) just before starting the circle method, but this time locally, we can use the $\sigma_{Np}$ defined in (\ref{def_sect_sigma_3}) to rewrite
\begin{equation} \label{descr_LambdaW}
\Lambda_N^W(p) \simeq \left\{\left. (\mathbf{v}, \mathbf{w}_p) \in P_{-l_p p|Np}^{n+1} \times \tilde W_p\,\right|\, f(\mathbf{v} + \sigma_{Np}(\mathbf{w}_p)) = 0 \right\} .
\end{equation}
In what follows we use a local version of orthogonality lemma \ref{lem_orth} as below, that has the same proof.
\begin{lemma} \label{lem_orth_loc}
For all $u \in P_{|Nx}$, over $\kappa(u)$ we have
\[
\left[ P_{|Nx}^\vee ,\, \gamma(u) \right]_\kappa =
\left\{ \begin{array}{ll}
\mathbf{L}^N & \text{if } u=0 \\
0 & \text{otherwise}
\end{array} \right. .
\]
\end{lemma}

Now we can make the connection between the $S_{[i]}$ and the spaces $\Lambda$.
\begin{lemma}
Placing ourselves over $\calB_\mu \times \scrC$, for $N \geq \max(l,1)$ and any closed $x\in\CmS$ we have
\begin{align*}
1 + \sum_{i=1}^N \mathbf{L}^{-(n+1)i} S_{[i]} = & \left[U\right] \cdot \mathbf{L}^{-nN} \left[ \Lambda_N(x) \right]_{J(\scrC) \times \tilde W} \\
& + \left[U^0\right] \cdot \mathbf{L}^{-nN + n+1} \left[ \Lambda_N^0(x) \right]_{J(\scrC) \times \tilde W} \\
& + \sum_{p\in|\scrS|} \left[U_p^\scrS\right] \cdot \mathbf{L}^{-nN + (n+1)l_p} \left[ \Lambda_N^W(p) \right]_{J(\scrC) \times \tilde W} .
\end{align*}
\end{lemma}

Following remark \ref{rk_not_class_subscript}, here $\left[ \Lambda_N(x) \right]_{J(\scrC) \times \tilde W}$ stands for $\left[ \Lambda_N(x) \times J(\scrC) \times \tilde W\right]$, and same with $\Lambda_N^0(x)$; while $\left[ \Lambda_N^W(p) \right]_{J(\scrC) \times \tilde W}$ is $\left[ \Lambda_N^W(p) \times_{\tilde W_p} \tilde W \times J(\scrC)\right]$.

\begin{proof}
We fix $N \geq \max(l,1)$ and $y\in\scrC$ closed . If $y=p\in|\scrS|$ we write $k = l_p$, otherwise we take $k\in\{0,1\}$. Let $1\leq i\leq N$ and denote $j = \min(k,i)$ and $j' = \max(k,i)$. For $\mathbf{w} \in \tilde W$ we put
\[
\varsigma_{Ny}(\mathbf{w}) =
\left\{ \begin{array}{ll}
\sigma_{Np}(\mathbf{w}_p) & \text{if } y=p\in|\scrS| \\
0 & \text{otherwise}
\end{array} \right.
\]
in $P_{|Ny}^{n+1}$. The class we are going to work with is
\[
C_{iy}^k := \mathbf{L}^{(n+1)j} \left[ P_{|iy}^{\vee*} \times P_{-jy|iy}^{n+1} \times \tilde W ,\, \beta\left(f\left( \mathbf{r} + \varsigma_{Ny}(\mathbf{w})|_{iy} \right)\right) \right] .
\]
Remark that with a closed $x\in\CmS$, lemma \ref{S_i} rewrites as
\begin{equation} \label{Lambda_intermed_1}
S_{[i]} = \left[ J(\scrC) \times U\right] C_{ix}^0 + \left[ J(\scrC) \times U^0\right] C_{ix}^1 + \sum_{p\in|\scrS|} \left[J(\scrC) \times U_p^\scrS\right] C_{ip}^{l_p} .
\end{equation}
With the statement we aim for in mind, we are first going to reformulate $1 + \sum_{i=1}^N \mathbf{L}^{-(n+1)i} C_{iy}^k$.

Let's replace the $P_{-jy|iy} = P_{-ky|iy}$ in $C_{iy}^k$ with $P_{-ky|Ny}$. Since $N \geq j'$, the exact sequence\footnote{Remark that this is a local version of exact sequence \ref{exseq_ed_vanish}. Furthermore when $x\notin|\scrS|$ and $k=0$, this is just exact sequence (\ref{exseq_ed_loc}) with $Z_1=Ny$ and $Z_2 = iy$.}
\[
\xymatrix{
0 \ar[r] & P_{-j'y|Ny} \ar[r]^-{\subset} & P_{-ky|Ny} \ar[r]^-{\cdot|_{ix}} & P_{-ky|iy} \ar[r] & 0
}
\]
is the same as
\[
\xymatrix{
0 \ar[r] & t^{j'} \C[t]/t^N \ar[r]^-{\subset} & t^k \C[t]/t^N \ar[r]^-{\bmod t^i} & t^k \C[t]/t^i \ar[r] & 0
} ;
\]
the only thing to check is the kernel of $\bmod\, t^i$: it is made up of the multiples of both $t^k$ and $t^i$, i.e. of the multiples of $t^{j'}$. Now note that $P_{-j'y|Ny}^{n+1}$ is a vector space of dimension $(n+1)(N-j')$, and apply the exact sequence to replace $\mathbf{r} \in P_{-ky|iy}^{n+1}$ with $\mathbf{v} \in P_{-ky|Ny}^{n+1}$:
\begin{align*}
\mathbf{L}^{-(n+1)i} C_{iy}^k & = \mathbf{L}^{-(n+1)(i-j)} \cdot \mathbf{L}^{-(n+1)(N-j')} \left[ P_{|iy}^{\vee*} \times P_{-ky|Ny}^{n+1} \times \tilde W ,\, \beta\left(f\left( \mathbf{v}|_{iy} + \varsigma_{Ny}(\mathbf{w})|_{iy} \right)\right) \right] \\
& = \mathbf{L}^{-(n+1)(N-k)} \left[ P_{|iy}^{\vee*} \times P_{-ky|Ny}^{n+1} \times \tilde W ,\, \beta\left(f\left( \mathbf{v} + \varsigma_{Ny}(\mathbf{w})\right)|_{iy}\right) \right]
\end{align*}
where we use that $i-j-j' = -k$. Next we also want to change the $\beta \in P_{|iy}^{\vee*}$ for some $\gamma \in P_{|Ny}^\vee$. For $0\leq i\leq N$ call $F_i$ the set of $\gamma \in P_{|Ny}^\vee$ that factor minimally through $iy$: this gives a partition
\[
P_{|Ny}^\vee = \bigsqcup_{i=0}^N F_i .
\]
There is an obvious isomorphism $P_{|iy}^{\vee*} \simeq F_i$ that allows us to rewrite
\[
\mathbf{L}^{-(n+1)i} C_{iy}^k = \mathbf{L}^{-(n+1)(N-k)} \left[ F_i \times P_{-ky|Ny}^{n+1} \times \tilde W ,\, \gamma\left(f\left( \mathbf{v} + \varsigma_{Ny}(\mathbf{w})\right)\right) \right]
\]
for $1\leq i\leq N$. As for $i=0$, $F_0 = \{0\}$ and above $\tilde W$,
\begin{align*}
\left[ F_0 \times P_{-ky|Ny}^{n+1} \times \tilde W ,\, \gamma\left(f\left( \mathbf{v} + \varsigma_{Ny}(\mathbf{w})\right)\right) \right] & = \left[ P_{-ky|Ny}^{n+1} \times \tilde W \right] \\
& = \mathbf{L}^{(n+1)(N-k)}
\end{align*}
so that
\[
1 = \mathbf{L}^{-(n+1)(N-k)} \left[ F_0 \times P_{-ky|Ny}^{n+1} \times \tilde W ,\, \gamma\left(f\left( \mathbf{v} + \varsigma_{Ny}(\mathbf{w})\right)\right) \right] .
\]
Putting this together we obtain
\begin{align*}
1 + \sum_{i=1}^N \mathbf{L}^{-(n+1)i} C_{iy}^k
& = \mathbf{L}^{-(n+1)(N-k)} \sum_{i=0}^N \left[ F_i \times P_{-ky|Ny}^{n+1} \times \tilde W ,\, \gamma\left(f\left( \mathbf{v} + \varsigma_{Ny}(\mathbf{w})\right)\right) \right] \\
& = \mathbf{L}^{-(n+1)(N-k)} \left[ P_{|Ny}^\vee \times P_{-ky|Ny}^{n+1} \times \tilde W ,\, \gamma\left(f\left( \mathbf{v} + \varsigma_{Ny}(\mathbf{w})\right)\right) \right] .
\end{align*}
This expression leads us to introduce
\[
A_N^k(y) := \left\{\left. (\mathbf{v}, \mathbf{w}) \in P_{-ky|Ny}^{n+1} \times \tilde W \,\right|\, f\left( \mathbf{v} + \varsigma_{Ny}(\mathbf{w})\right) = 0 \right\} .
\]
Applying previous lemma \ref{lem_orth_loc} to this in the exact same way as we used the orthogonality lemma \ref{lem_orth} to prove proposition \ref{init}, we obtain, in $\mathscr{E}xp\mathscr{M}_{P_{-ky|Ny}^{n+1} \times \tilde W}$,
\[
\left[ A_N^k(y) \right] = \mathbf{L}^{-N} \left[ P_{|Ny}^\vee \times P_{-ky|Ny}^{n+1} \times \tilde W ,\, \gamma\left(f\left( \mathbf{v} + \varsigma_{Ny}(\mathbf{w})\right)\right) \right] :
\]
thus
\begin{equation} \label{Lambda_intermed_2}
1 + \sum_{i=1}^N \mathbf{L}^{-(n+1)i} C_{iy}^k = \mathbf{L}^{-nN + (n+1)k} \left[ A_N^k(y) \right] .
\end{equation}

To conclude, remark that if $y\in\CmS$, $A_N^0(y) = \Lambda_N(y) \times \tilde W$ and $A_N^1(y) = \Lambda_N^0(y) \times \tilde W$, while for $p\in|\scrS|$, $A_N^{l_p}(p) \simeq \Lambda_N^W(p) \times_{\tilde W_p} \tilde W$. Combine the relationship (\ref{Lambda_intermed_1}) between $S_{[i]}$ and the $C_{iy}^k$, the additional fact that over $\calB_\mu \times \scrC$,
\[
1 = \left[ J(\scrC) \times U\right] \big[\tilde W\big] + \left[ J(\scrC) \times U^0\right] \big[\tilde W\big] + \sum_{p\in|\scrS|} \left[J(\scrC) \times U_p^\scrS\right] \big[\tilde W\big] ,
\]
and (\ref{Lambda_intermed_2}) to obtain the desired formula for $1 + \sum_{i=1}^N \mathbf{L}^{-(n+1)i} S_{[i]}$.
\end{proof}

The identity stated in the lemma can be taken to the limit for $N\to\infty$, since the series on the left-hand side converges (see the remark before corollary \ref{sing_mep_2}). Actually, not only does the right-hand side converge, but each of its terms
\[
\mathbf{L}^{-nN} \cdot \left[ \Lambda_N(x) \right]_{J(\scrC) \times \tilde W}, \quad \mathbf{L}^{-nN + n+1} \cdot \left[ \Lambda_N^0(x) \right]_{J(\scrC) \times \tilde W}, \quad \mathbf{L}^{-nN + (n+1)l_p} \cdot \left[ \Lambda_N^W(p) \right]_{J(\scrC) \times \tilde W}
\]
does so individually: by evaluating the identity at a $(\scrT,x)$ in respectively $U$, $U^0$ or $U_p^\scrS$ we keep only the corresponding term and the convergence holds. This gives
\begin{align*}
1 + \sum_{i\geq1} \mathbf{L}^{-(n+1)i} S_{[i]} = & \left[U\right] \cdot \lim_{N\to\infty} \mathbf{L}^{-nN} \left[ \Lambda_N(x) \right]_{J(\scrC) \times \tilde W} \\
& + \left[U^0\right] \cdot \lim_{N\to\infty} \mathbf{L}^{-nN + n+1} \left[ \Lambda_N^0(x) \right]_{J(\scrC) \times \tilde W} \\
& + \sum_{p\in|\scrS|} \left[U_p^\scrS\right] \cdot \lim_{N\to\infty} \mathbf{L}^{-nN + (n+1)l_p} \left[ \Lambda_N^W(p) \right]_{J(\scrC) \times \tilde W} .
\end{align*}
In particular this is constant over $\CmS$; more precisely if $x\in\CmS$, with $U_x = \Sym^\mu(\scrC \smallsetminus |\scrS+x|) \times \kappa(x)$ and $U_x^0$ its complement in $\Sym^\mu(\CmS) \times \kappa(x)$,
\begin{align} \label{locdens_x}
1 + \sum_{i\geq1} \mathbf{L}^{-(n+1)i} S_{[i],x} = & \left[U_x\right] \cdot \lim_{N\to\infty} \mathbf{L}^{-nN} \left[ \Lambda_N(x) \right]_{J(\scrC) \times \tilde W} \\
& + \left[U_x^0\right] \cdot \lim_{N\to\infty} \mathbf{L}^{-nN + n+1} \left[ \Lambda_N^0(x) \right]_{J(\scrC) \times \tilde W} . \nonumber
\end{align}
On the other hand for $p\in|\scrS|$,
\begin{equation} \label{locdens_p}
1 + \sum_{i\geq1} \mathbf{L}^{-(n+1)i} S_{[i],p} = \left[U_p^\scrS\right] \cdot \lim_{N\to\infty} \mathbf{L}^{-nN + (n+1)l_p} \left[ \Lambda_N^W(p) \right]_{J(\scrC) \times \tilde W} .
\end{equation}
Combining this, corollary \ref{sing_mep_2}, and the fact that we can take a finite number of factors outside the motivic Euler product, we deduce the following.

\begin{corollary} \label{sing_locdens}
In $\widehat{\mathscr{E}xp\mathscr{M}}_{\calB_\mu}$ the singular series is equal to
\[
\mathfrak{S}^\mu = \prod_{x\in\CmS} \left( 1 + \sum_{i\geq1} \mathbf{L}^{-(n+1)i} S_{[i],x} \right) \cdot \prod_{p\in|\scrS|} \left( 1 + \sum_{i\geq1} \mathbf{L}^{-(n+1)i} S_{[i],p} \right) ,
\]
where the factors are local densities via (\ref{locdens_x}) and (\ref{locdens_p}) above.
\end{corollary}
The first product is still a motivic Euler product, but the second one is more simply a finite product.

\subsubsection{One last expression for the singular series}

To finish with the treatment of the major arcs, we can in turn reformulate the classes of the $\Lambda$ spaces and get the following final expression for the singular series, before a quick conclusion. Recall that $\calB_\mu = J(\scrC) \times \Sym^\mu(\CmS) \times \tilde W$.

\begin{proposition} \label{sing_prod_final}
We can rewrite the singular series as
\begin{equation}
\mathfrak{S}^\mu = \mathbf{L}^{l+d|\mu|} \prod_{x\in\CmS}\left( \frac{1-\mathbf{L}^{-1}}{1-\mathbf{L}^{-(n-d+1)}} \cdot \frac{[V]}{\mathbf{L}^{n-1}} \right) \in \widehat{\mathscr{M}}_{\calB_\mu} .
\end{equation}
In this formula, the right-hand side is pulled back via $\calB_\mu \to \Spec\C$; in particular the motivic Euler product writes, over $\C$, as $1$ plus a class of negative weight.
\end{proposition}
It may seem odd that $l$ and $|\mu|$ appear with positive coefficients in the normalisation, but we will see in (\ref{link_zeta_xi}) that this simply balances the negative multiplicity with which these two quantities already appeared in the former normalisation $\zeta_{e,|\mu|}$. Furthermore, observe that
\begin{equation} \label{sing_link_mu_empty}
\mathfrak{S}^\mu = \mathbf{L}^{d|\mu|} \mathfrak{S}^\varnothing
\end{equation}
over $\calB_\mu$ (with $\mathfrak{S}^\varnothing$ originally over $J(\scrC) \times \tilde W$).

\begin{proof}
The starting point is of course corollary \ref{sing_locdens}. During this proof we denote
\begin{equation} \label{sing_factor}
\mathcal{F} := \frac{1-\mathbf{L}^{-1}}{1-\mathbf{L}^{-(n-d+1)}} \cdot \frac{[V]}{\mathbf{L}^{n-1}} \in \widehat{\mathscr{M}}_\C .
\end{equation}
Since $V$ is irreducible, lemma \ref{nb_compo_max} implies that $w(\mathcal{F}-1) < 0$. First take any closed $x\in\scrC$. Using $P_{|Nx} = \C[t]/t^N$ and then lifting, our $\Lambda_N(x)$ as defined in (\ref{def_Lambda}) is isomorphic to the
\[
\Lambda_N(f,0) = \left\{ \mathbf{g} \in \C[t]_{<N}^{n+1} \,\left|\, f(\mathbf{g}) \equiv 0 \pmod{t^N} \right.\right\}
\]
from \cite{BB25}, defined in (1.2) (once more, mind that their $n$ corresponds to our $n+1$). In particular we can directly apply their lemma 8.2, which provides
\begin{equation} \label{lim_Lambda}
\lim_{N\to\infty} \mathbf{L}^{-nN} \left[\Lambda_N(x)\right]_{J(\scrC) \times \tilde W} = \left[J(\scrC) \times \tilde W\right] \cdot \mathcal{F} .
\end{equation}
For the rest of the proof (i.e., what happens when adding the conditions) we reproduce arguments used to show this result, namely from the proof of their lemma 6.13.

Suppose $x \in \CmS$ and take $N\geq d$. By definition (\ref{def_Lambda0}) and through the above correspondance,
\begin{align*}
\Lambda_N^0(x) & \simeq \left\{ \mathbf{g} \in t\C[t]_{<N-1}^{n+1} \,\left|\, f(\mathbf{g}) \equiv 0 \pmod{t^N} \right.\right\} \\
& \simeq \left\{ \mathbf{h} \in\C[t]_{<N-1}^{n+1} \,\left|\, f(\mathbf{h}) \equiv 0 \pmod{t^{N-d}} \right.\right\}
\end{align*}
by writing $\mathbf{g} = t\mathbf{h}$ and using that $f$ is homogeneous of degree $d$. It is clear that in the last space, the coefficients $N-d$ to $N-2$ of $\mathbf{h}$ can be chosen arbitrarily, while its remainder modulo $t^{N-d}$ varies in $\Lambda_{N-d}(x)$. This implies $\Lambda_N^0(x) \simeq \mathbb{A}^{(d-1)(n+1)} \times \Lambda_{N-d}(x)$, and in particular
\[
\mathbf{L}^{-nN + n+1} \big[\Lambda_N^0(x)\big]_{J(\scrC) \times \tilde W} = \mathbf{L}^d \cdot \mathbf{L}^{-n(N-d)} \left[\Lambda_{N-d}(x)\right]_{J(\scrC) \times \tilde W} \xrightarrow[N\to\infty]{} \mathbf{L}^d \left[J(\scrC) \times \tilde W\right] \cdot \mathcal{F} .
\]
Injecting this and (\ref{lim_Lambda}) into (\ref{locdens_x}) shows that the motivic Euler product over $\CmS$ from corollary \ref{sing_locdens} is
\[
\prod_{x\in\CmS}\left( \left(\left[U_x\right] + \left[U_x^0\right] \mathbf{L}^d\right) \left[J(\scrC) \times \tilde W\right] \cdot \mathcal{F} \right) = \left[J(\scrC) \times \tilde W\right] \cdot \prod_{x\in\CmS}\left( \left(\left[U_x\right] + \left[U_x^0\right] \mathbf{L}^d\right) \mathcal{F} \right)
\]
(the motivic Euler product on the left lives over $\calB_\mu$, while the one on the right lives over $\Sym^\mu(\CmS)$). We said under (\ref{sing_factor}) that $w(\mathcal{F}-1) < 0$, as a consequence $\prod_{x\in\CmS} \mathcal{F}$ converges and satisfies $w\left(\prod_x \mathcal{F} - 1\right) < 0$ itself. We then separate the former product into
\[
\left[J(\scrC) \times \tilde W\right] \cdot \prod_{x\in\CmS}\left(\left[U_x\right] + \left[U_x^0\right] \mathbf{L}^d\right) \cdot \prod_{x\in\CmS} \mathcal{F} :
\]
we show in the next lemma \ref{prod_U_outsource} that the product in the middle is $\mathbf{L}_{\Sym^\mu(\CmS)}^{d|\mu|}$. As a consequence we have, for this part,
\begin{equation} \label{sing_prod_final_intermed}
\mathbf{L}_{\calB_\mu}^{d|\mu|} \cdot \prod_{x\in\CmS} \mathcal{F} .
\end{equation}

Finally let us treat the factors for $p\in|\scrS|$ from corollary \ref{sing_locdens}. First, going back to description (\ref{descr_LambdaW}), remark that
\[
\Lambda_{l_p}^W(p) = \tilde W_p .
\]
Now for $N\geq l_p$, like previously we have
\[
\Lambda_N^W(p) \simeq \left(\Lambda_N^W(p)\right)_\mathrm{poly} :=  \left\{ \mathbf{g} \in \C[t]_{<N}^{n+1} \,\left|\,
\begin{array}{l}
f(\mathbf{g}) \equiv 0 \pmod{t^N} \\
\mathbf{g} \bmod t^{l_p} \in \big(\tilde W_p\big)_\mathrm{poly}
\end{array} \right.\right\}
\]
with $\big(\tilde W_p\big)_\mathrm{poly}$ the image of $\tilde W_p$ in $\left(\C[t]/t^{l_p}\right)^{n+1}$. One can consider the surjective morphism
\[
\bmod ~ t^N : \left(\Lambda_{N+1}^W(p)\right)_\mathrm{poly} \longrightarrow \left(\Lambda_N^W(p)\right)_\mathrm{poly} ,
\]
whose fibre above $\mathbf{h} \in \left(\Lambda_N^W(p)\right)_\mathrm{poly}$ is made up of the $\mathbf{h} + \mathbf{a}t^N$ with $\mathbf{a}\in\C^{n+1}$ such that $f(\mathbf{h} + \mathbf{a}t^N) \equiv 0 \pmod{t^{N+1}}$. By Taylor's formula, and then the fact that $f(\mathbf{h}) \equiv 0 \pmod{t^N}$ so that there exists $b\in\C$ such that $f(\mathbf{h}) \equiv b t^N \pmod{t^{N+1}}$,
\begin{align*}
f(\mathbf{h} + \mathbf{a}t^N) & \equiv f(\mathbf{h}) + \nabla f(\mathbf{h}(0)) \cdot \mathbf{a} t^N \pmod{t^{N+1}} \\
& \equiv \left(b + \nabla f(\mathbf{h}(0)) \cdot \mathbf{a}\right) t^N \pmod{t^{N+1}} .
\end{align*}
The $\mathbf{a}$ that compose the fibre are the ones that make the coefficient on the second line zero. Since $\mathbf{h}(0)$ is the constant coefficient of an element of $\big(\tilde W_p\big)_\mathrm{poly}$, it is non-zero. Furthermore by assumption $f$ is non-singular outside zero, so $\nabla f(\mathbf{h}(0)) \neq 0$. This implies that $n$ components of $\mathbf{a}$ can be chosen freely while the last one will be imposed by the vanishing of the above coefficient. The fibre is thus isomorphic to $\mathbb{A}^n$. Going back to the world of sections (and taking the product with $\tilde W$) this leads to
\[
\left[\Lambda_{N+1}^W(p)\right]_{J(\scrC) \times \tilde W} = \mathbf{L}^n \left[\Lambda_N^W(p)\right]_{J(\scrC) \times \tilde W} .
\]
Iterating that,
\begin{align} \label{lim_S}
\mathbf{L}^{-nN + (n+1)l_p} \left[\Lambda_N^W(p)\right]_{J(\scrC) \times \tilde W} & = \mathbf{L}^{-nN + (n+1)l_p} \cdot \mathbf{L}^{n(N-l_p)} \big[\Lambda_{l_p}^W(p)\big]_{J(\scrC) \times \tilde W} \nonumber \\
& = \mathbf{L}^{l_p} \big[\tilde W_p\big]_{J(\scrC) \times \tilde W} \nonumber \\
& = \mathbf{L}^{l_p}
\end{align}
(recall the convention from remark \ref{rk_not_class_subscript}). Remark that this is even constant in $N$ and does not require taking a limit. It follows that the finite product of local densities over $|\scrS|$ from corollary \ref{sing_locdens} is $\prod_{p\in|\scrS|} \mathbf{L}^{l_p} = \mathbf{L}^l$. Putting (\ref{sing_prod_final_intermed}) with this, we obtain
\[
\mathfrak{S}^\mu = \mathbf{L}_{\calB_\mu}^{l+d|\mu|} \prod_{x\in\CmS} \mathcal{F}
\]
as announced.
\end{proof}

In this proof we outsourced the study of a motivic Euler product to the following lemma.
\begin{lemma} \label{prod_U_outsource}
We have
\[
\prod_{x\in\CmS}\left(\left[U_x\right] + \left[U_x^0\right] \mathbf{L}^d\right) = \mathbf{L}^{d|\mu|}
\]
in $K_0\left(\mathrm{Var}_{\Sym^\mu(\CmS)}\right)$.
\end{lemma}
\begin{proof}
During this proof we denote the product on the left-hand side $\mathfrak{P}^\mu$. By construction (\ref{def_U}), $[U]+[U^0] = 1$ over $\Sym^\mu(\CmS) \times (\CmS)$, and the factors are equal to $1 + [U_x^0](\mathbf{L}^d-1)$. After remark \ref{MEP_noT_sum},
\[
\mathfrak{P}^\mu = \sum_{j\geq0} \left[\Conf^{j[1]}\left( U^0 \times (\mathbb{A}^d \smallsetminus \{0\}) \right)\right]
\]
where the subscript to the $\Conf$ is
\[
\Sym^\mu(\CmS) \times (\CmS) \ \big/ \ \Sym^\mu(\CmS) .
\]
The above configuration spaces are made up of $j$-tuples of elements of $U^0 \times (\mathbb{A}^d \smallsetminus \{0\})$ that map to distinct points of $\CmS$ but to a common $\scrT \in \Sym^\mu(\CmS)$, up to permutation. Since $U_\scrT^0 = |\scrT| \subset \CmS$ by definition (\ref{def_U}), the fibre over such a $\scrT$ is
\[
\Conf^{j[1]}\left( U^0 \times (\mathbb{A}^d \smallsetminus \{0\}) \right)_\scrT = \left\{\left. (\mathbf{x},\mathbf{a}) \in |\scrT|^j \times (\mathbb{A}^d \smallsetminus \{0\})^j \,\right|\, \forall k \neq k',\, x_k \neq x_{k'} \right\} / \mathfrak{S}_j .
\]
This fibre is isomorphic to the finite disjoint union
\[
\bigsqcup_{\substack{I \subset |\scrT| \\ \#I = j}} (\mathbb{A}^d \smallsetminus \{0\})^I
\]
by sending $(\mathbf{x},\mathbf{a})$ to
\[
\left\{\begin{array}{rcl}
\{x_1, \dots, x_j\} & \longrightarrow & \mathbb{A}^d \smallsetminus \{0\} \\
x_k & \longmapsto & a_k
\end{array}\right. ;
\]
i.e. it is composed of $\binom{|\mu|}{j}$ copies of $(\mathbb{A}^d \smallsetminus \{0\})^j$. The consequence is that the configuration space is empty as soon as $j>|\mu|$, and for any $\scrT \in \Sym^\mu(\CmS)$ we have, over $\kappa(\scrT)$,
\[
\mathfrak{P}^\mu(\scrT) = \sum_{j=0}^{|\mu|} \binom{|\mu|}{j} \left( \mathbf{L}^d-1 \right)^j = \mathbf{L}^{d|\mu|}
\]
by the binomial identity. Lemma \ref{lift_equal} concludes.
\end{proof}

Let us sum up our results regarding the major arcs. First recall the definition of $\zeta_{e,|\mu|}$ (\ref{def_zeta}), and denote by $R_e$ the remainder of the singular series, that we bounded in lemma \ref{sing_cv}. Proposition \ref{maj_1} yields
\begin{equation} \label{maj_concl}
\Gamma_\mathfrak{M} = \mathbf{L}^{\zeta_{e,|\mu|}} (\mathfrak{S}^\mu - R_e) :
\end{equation}
by proposition \ref{sing_prod_final} this should have $\calB_\mu$-dimension
\begin{align} \label{link_zeta_xi}
\zeta_{e,|\mu|} + l + d|\mu| & = (n-d+1)e - n(g-1) - nl - (n-d+1)|\mu| \\
& = \xi_{e,|\mu|} , \nonumber
\end{align}
as expected since remark \ref{expect_dim_naive}. The coefficients of $\xi_{e,|\mu|}$ in $l$ and $|\mu|$ are still negative: in fact the idea is that when we defined $\zeta_{e,|\mu|}$ we substracted a bit too many $l' = l + |\mu|$ to compensate for the conditions at $\scrST$, and now we get to balance it perfectly by adding back $l + d|\mu|$.

\subsection{Minor arcs}
\label{subsect_min}

Finally we adress the minor arcs and their contribution $\Gamma_\mathfrak{m}$. This is possible thanks to the global version of Weyl's inequality, proposition \ref{concl_Weyl_global}. Just above it we were led to define a quantity $s(\delta)$ (\ref{def_s}) as the maximum of three possible expressions $s_i(\delta)$: in what follows, this disjunction naturally leads up to several quantities appearing as maxima (or minima) of three expressions, each relying on one of the $s_i(\delta)$. Recall
\[
q = \frac{n+1}{2^{d-2}(d-1)} > 4
\]
from (\ref{def_q}).

\begin{proposition} \label{min_prop}
Suppose that $e \geq \max_{1\leq i\leq3} \theta_i$, where we put
\[
\left\{\begin{array}{l}
\theta_1 := 4g-\inddd+2l'\\
\theta_2 := 2dg-2+l' \\
\theta_3 := 2g+1+l'
\end{array}\right. .
\]
Then we have
\[
w_{\calB_\mu}(\Gamma_\mathfrak{m}) \leq 2\xi_{e,|\mu|} + \max_{1\leq i\leq3} \rho_i
\]
with first
\begin{align*}
& \rho_1 := -\frac{d(q-4)}{2}e + \phi_{d,q}^1(g) + (dq-2)l + d(q-2)|\mu| \\
& \text{where}\quad \phi_{d,q}^1(g) :=
4 + q\cdot \left\{\begin{array}{ll}
-1 & \text{if } g=0 \\
d & \text{if } g=1 \\
\frac{3(d+1)}{2}g + \frac{(d-3)}{2} & \text{if } g\geq2
\end{array}\right. ,
\end{align*}
second
\begin{align*}
& \rho_2 := -(q-4)e + \phi_{d,q}^2(g) + (q-6)l + (q-2d-4)|\mu| \\
& \text{where}\quad \phi_{d,q}^2(g) :=
\left\{\begin{array}{ll}
-2(q-4) & \text{if } g=0 \\
(d-1)q & \text{if } g=1 \\
\left(\frac{(3d+1)q}{2}-8\right)g + \frac{(d-5)q}{2}+8 & \text{if } g\geq2
\end{array}\right. ,
\end{align*}
and third
\begin{align*}
& \rho_3 := -((d-1)(q-2)-2)e + \phi_{d,q}^3(g) + ((d-1)q-2)l + ((d-1)(q-2)-2)|\mu| \\
& \text{where}\quad \phi_{d,q}^3(g) :=
4 + (d-1)q\cdot \left\{\begin{array}{ll}
1 & \text{if } g=0 \\
2 & \text{if } g=1 \\
\frac{7g-3}{2} & \text{if } g\geq2
\end{array}\right. .
\end{align*}
\end{proposition}
Before proving this statement, remark that the three possible bounds are of the form $ae + \phi_{d,q}^i(g) + bl + c|\mu|$ with $a,b,c$ constants and $\phi_{d,q}^i(g)$ a piecewise linear function of $g$. All are important in the proposition because $\rho_1$ has the largest coefficients in $l$ and $|\mu|$, $\rho_2$ has the largest coefficient in $e$, and $\rho_3$ has potentially the largest coefficient in $g$ (depending on the value of $d$). Remark also that since $d\geq2$ and $q>4$, all three coefficients for $e$ are negative.

\begin{proof}
Instead of applying Weyl's inequality \ref{concl_Weyl_global} once to $\mathfrak{m}$ as a whole, we will partition it and then apply the bound to each piece. By definition, minor arcs are made up of the $\alpha$ with $\deg\alpha \geq \eta_e+1 = e-l'-2g+2$. Additionally lemma \ref{max_deg} shows that one always has $\deg\alpha \leq \frac{d}{2}e+1$. As a consequence,
\begin{equation} \label{part_min}
\mathfrak{m} = \bigsqcup_{\delta = e-l'-2g+2}^{\left\lfloor\frac{d}{2}e\right\rfloor+1} \mathfrak{m}_\delta \quad\text{with}\quad \mathfrak{m}_\delta := \left\{\left. (L,\scrT,\alpha) \in \tau_d^*\calP_{de}^\vee \,\right|\, \deg\alpha = \delta \right\} .
\end{equation}
In order to apply the proposition the condition $s(\delta) \leq e-l'-2g+1$ needs to be satisfied for all $e-l'-2g+2 \leq \delta \leq \frac{d}{2}e+1$. Since $s(\delta)$ is defined in (\ref{def_s}) as the maximum of three expressions, let us look at these three cases successively. First
\begin{align*}
s_1(\delta) & = \left\lfloor \frac{\delta-e+l'+2g-2}{d-1} \right\rfloor + 1 \\
& \leq \frac{\frac{d}{2}e+1-e+l'+2g-2}{d-1} + 1 \\
& = \frac{\left(\frac{d}{2}-1\right)e+l'+2g-1}{d-1} + 1 .
\end{align*}
This last quantity is smaller than $e-l'-2g+1$ if and only if
\begin{align*}
\mkern-36mu && \left(\frac{d}{2}-1\right)e+l'+2g-1 & \leq (d-1)(e-l'-2g) \\
\Longleftrightarrow \mkern-108mu && d(l'+2g)-1 &\leq \frac{d}{2}e \\
\Longleftrightarrow \mkern-108mu && e & \geq 2l'+4g-\frac{2}{d} .
\end{align*}
Because $d\geq2$ and $e$ is an integer, $-\frac{2}{d}$ can be replaced with $-\inddd$: this is the hypothesis $e\geq\theta_1$. For the next one:
\begin{align*}
s_2(\delta) & = e-l'+1 - \left\lceil \frac{\delta}{d-1} \right\rceil \\
& \leq e-l'+1 - \frac{e-l'-2g+2}{d-1} ,
\end{align*}
with the second line smaller than $e-l'-2g+1$ when
\begin{align*}
\mkern-48mu && -\frac{e-l'-2g+2}{d-1} & \leq -2g \\
\Longleftrightarrow \mkern-144mu && e-l'-2g+2 & \geq 2(d-1)g \\
\Longleftrightarrow \mkern-144mu && e & \geq l'+2dg-2 .
\end{align*}
This justifies the hypothesis $e\geq\theta_2$. Finally
\[
s_3(\delta) = s_3 = \max(2g-1, 2):
\]
on the one hand
\[
2g-1 \leq e-l'-2g+1 \quad\Longleftrightarrow\quad e \geq l'+4g-2,
\]
which is already implied by $e\geq\theta_1$. On the other hand
\[
2 \leq e-l'-2g+1 \quad\Longleftrightarrow\quad e \geq l'+2g+1,
\]
which results in the hypothesis $e\geq\theta_3$. All this justifies that our hypothesis with the $\theta_i$ allows us to apply proposition \ref{concl_Weyl_global}, which we are going to do soon.

Now towards the bound itself: by definition of $\Gamma_\mathfrak{m}$ just above identity (\ref{part_maj_min}), by relation (\ref{link_zeta_xi}), and by definition of $\zeta_{e,l'}$ in (\ref{def_zeta}),
\begin{align*}
& w_{\calB_\mu}(\Gamma_\mathfrak{m}) -2\xi_{e,|\mu|} \\
& = 2(-de+g-1) + w_{\calB_\mu}\left(\left[ \mathfrak{m} \,\diamond\, \calP_{e(-\scrS_\mu)}^{\diamond(n+1)} \times \tilde W ,\, \alpha\left(f\left( \mathbf{v} + \sigma_L^\scrT(\mathbf{w}) \right)\right) \right]\right) - 2(\zeta_{e,l'} + l + d |\mu|) \\
& = - 2(n+1)(e-l'-g+1) - 2l - 2d|\mu| + w_{\calB_\mu}\left(\left[ \mathfrak{m} \,\diamond\, \calP_{e(-\scrS_\mu)}^{\diamond(n+1)} \times \tilde W ,\, \alpha\left(f\left( \mathbf{v} + \sigma_L^\scrT(\mathbf{w}) \right)\right) \right]\right) .
\end{align*}
Partition (\ref{part_min}) then implies that
\[
w_{\calB_\mu}(\Gamma_\mathfrak{m}) -2\xi_{e,|\mu|} \leq \max_{e-l'-2g+2 \leq \delta \leq \frac{d}{2}e+1} B(\delta)
\]
where
\begin{align*}
B(\delta) := & -2(n+1)(e-l'-g+1) - 2l - 2d|\mu| \\
& + w_{\calB_\mu}\left(\left[ \mathfrak{m}_\delta \,\diamond\, \calP_{e(-\scrS_\mu)}^{\diamond(n+1)} \times \tilde W ,\, \alpha\left(f\left( \mathbf{v} + \sigma_L^\scrT(\mathbf{w}) \right)\right) \right]\right)
\end{align*}
is the same expression except that $\mathfrak{m}$ has been replaced with $\mathfrak{m}_\delta$. Now the weight on the second line is at most
\begin{equation} \label{min_intermed}
2\dim_{\calB_\mu}(\mathfrak{m}_\delta \times \tilde W) + w_{\mathfrak{m}_\delta \times \tilde W}\left(\left[ \mathfrak{m}_\delta \,\diamond\, \calP_{e(-\scrS_\mu)}^{\diamond(n+1)} \times \tilde W ,\, \alpha\left(f\left( \mathbf{v} + \sigma_L^\scrT(\mathbf{w}) \right)\right) \right]\right)
\end{equation}
by point 2 of proposition \ref{prop_weight}. On the one hand, $\dim_{\calB_\mu}(\mathfrak{m}_\delta \times \tilde W) = \dim_{\Pic^e(\scrC) \times \Sym^\mu(\CmS)} \mathfrak{m}_\delta$. Consider the vector bundle $\calP_{|\nu}$ over
\[
\Pic^e(\scrC) \times \Sym^\mu(\CmS) \times \Sym^\nu(\scrC)
\]
with fiber $P_{|Z}$ at $(L,\scrT,Z)$. By minimal factorisation of linear forms and by definition of $\mathfrak{m}_\delta$, every $(L,\scrT,\alpha) \in \mathfrak{m}_\delta$ has its $\alpha$ coming from $P_{|Z}^\vee$ for some $Z \in \Sym^\delta(\scrC)$, which means that the map
\[
\left\{\begin{array}{rcl}
\bigsqcup_{\nu\vdash\delta} \calP_{|\nu}^\vee & \longrightarrow & \mathfrak{m}_\delta \\
(L, \scrT, Z, \beta) & \longmapsto & (L, \scrT, \beta(\cdot|_Z))
\end{array}\right.
\]
is surjective (note that $\beta \mapsto \beta(\cdot|_Z)$ is the map $P_{|Z}^\vee \to P_{L^{\otimes d}}^\vee$ from (\ref{exseq_dual})). As a consequence, relatively to $\Pic^e(\scrC) \times \Sym^\mu(\CmS)$, the dimension of $\mathfrak{m}_\delta$ is at most that of $\bigsqcup_{\nu\vdash\delta} \calP_{|\nu}^\vee$. The latter is the maximum of the dimensions of the $\calP_{|\nu}$, and each of those is the sum of $|\nu|$ (for $Z \in \Sym^\nu(\scrC)$) and $\delta$ (for $\beta \in P_{|Z}$), so it is $\leq 2\delta$. Finally,
\[
2\dim_{\calB_\mu}(\mathfrak{m}_\delta \times \tilde W) \leq 4\delta .
\]
On the other hand, as announced we apply proposition \ref{concl_Weyl_global} to the weight in (\ref{min_intermed}), with $A=\mathfrak{m}_\delta$. Thus
\begin{align} \label{min_bound_B}
B(\delta) & \leq 2(n+1)(e-l'-g+1) - 2l - 2d|\mu| + 4\delta \nonumber \\
& \qquad + 2^{-(d-2)}(n+1) \left(s(\delta) + \left(2^{d-1}-1\right)(e-l'-g+1) + \frac{g+1}{2}\mathbbm{1}_{g\geq2} \right) \nonumber \\
& = -2l - 2d|\mu| + 4\delta + 2^{-(d-2)}(n+1) \left(s(\delta) - (e-l'-g+1) + \frac{g+1}{2}\mathbbm{1}_{g\geq2} \right) \nonumber \\
& = -2l - 2d|\mu| + 4\delta + (d-1)q \left(s(\delta) - e + l' + g - 1 + \frac{g+1}{2}\mathbbm{1}_{g\geq2} \right) .
\end{align}

The last step is to inject the three possibilities for $s(\delta)$ in the last formula. First we substitute it with
\[
s_1(\delta) = \left\lfloor \frac{\delta-e+l'+2g-2}{d-1} \right\rfloor + 1 ,
\]
so that this version of expression (\ref{min_bound_B}) is bounded by
\begin{align*}
& -2l - 2d|\mu| + 4\delta + (d-1)q \left(\frac{\delta-e+l'+2g-2}{d-1} + 1 - e + l' + g - 1 + \frac{g+1}{2}\mathbbm{1}_{g\geq2} \right) \\
& = -2l - 2d|\mu| + (q+4)\delta + q\left( -de + dl' + (d+1)g - 2 + (d-1)\frac{g+1}{2}\mathbbm{1}_{g\geq2} \right) \\
& \leq -2l - 2d|\mu| + (q+4)\left(\frac{d}{2}e+1\right) + q\left( -de + dl' + (d+1)g - 2 + (d-1)\frac{g+1}{2}\mathbbm{1}_{g\geq2} \right) \\
& = d\left(\frac{q+4}{2}-q\right)e + dql' - 2l - 2d|\mu| + q\left( (d+1)g - 1 + (d-1)\frac{g+1}{2}\mathbbm{1}_{g\geq2} \right) + 4 \\
& = -\frac{d(q-4)}{2}e + (dq-2)l + d(q-2)|\mu| + \phi_{d,q}^1(g) \\
& = \rho_1 .
\end{align*}
Between the fourth and fifth lines we use that $l' = l+|\mu|$; one checks easily that the last part of the expression is the $\phi_{d,q}^1(g)$ announced in the statement. Next, expression (\ref{min_bound_B}) with $s(\delta)$ replaced by
\[
s_2(\delta) = e-l'+1 - \left\lceil \frac{\delta}{d-1} \right\rceil
\]
is at most
\begin{align*}
& -2l - 2d|\mu| + 4\delta + (d-1)q \left( e-l'+1 - \frac{\delta}{d-1} - e + l' + g - 1 + \frac{g+1}{2}\mathbbm{1}_{g\geq2} \right) \\
& = -2l - 2d|\mu| - (q-4)\delta + (d-1)q\left( g + \frac{g+1}{2}\mathbbm{1}_{g\geq2} \right) \\
& \leq -2l - 2d|\mu| - (q-4)(e-l'-2g+2) + (d-1)q\left( g + \frac{g+1}{2}\mathbbm{1}_{g\geq2} \right) \\
& = -(q-4)e + (q-4)l' - 2l - 2d|\mu| + 2(q-4)(g-1) + (d-1)q\left( g + \frac{g+1}{2}\mathbbm{1}_{g\geq2} \right) \\
& = -(q-4)e + (q-6)l + (q-2d-4)|\mu| + \phi_{d,q}^2(g) \\
& = \rho_2 .
\end{align*}
Finally we exchange $s(\delta)$ for
\[
s_3 = \max(2g-1, 2) = 2 + (2g-3)\mathbbm{1}_{g\geq2}
\]
in (\ref{min_bound_B}), which gives
\begin{align*}
& -2l - 2d|\mu| + 4\delta + (d-1)q \left( 2 + (2g-3)\mathbbm{1}_{g\geq2} - e + l' + g - 1 + \frac{g+1}{2}\mathbbm{1}_{g\geq2} \right) \\
& = -2l - 2d|\mu| + 4\delta + (d-1)q\left( -e + l' + g + 1 + \frac{5}{2}(g-1)\mathbbm{1}_{g\geq2} \right) \\
& \leq -2l - 2d|\mu| + 4\left(\frac{d}{2}e+1\right) + (d-1)q\left( -e + l' + g + 1 + \frac{5}{2}(g-1)\mathbbm{1}_{g\geq2} \right) \\
& = (2d-(d-1)q)e + (d-1)ql' - 2l - 2d|\mu| + (d-1)q\left( g + 1 + \frac{5}{2}(g-1)\mathbbm{1}_{g\geq2} \right) + 4 \\
& = -((d-1)(q-2)-2)e + ((d-1)q-2)l + ((d-1)(q-2)-2)|\mu| + \phi_{d,q}^3(g) \\
& = \rho_3 .
\end{align*}

We can now conclude: these computations show that $B(\delta) \leq \max_{1\leq i\leq3} \rho_i$ for all relevant~$\delta$, therefore we also have
\[
w_{\calB_\mu}(\Gamma_\mathfrak{m}) -2\xi_{e,|\mu|} \leq \max_{1\leq i\leq3} \rho_i .
\]
\end{proof}

\begin{remark} \label{discard_th3}
Looking at the condition on $e$ in the last proposition, it appears that most of the time $\theta_1 \geq \theta_3$; namely, $\theta_1-\theta_3 = 2g - 1 - \inddd + l'$ is only negative in the following cases: $g=l'=0$; or $d=2$, $g=0$ and $l'=1$. The value of $\theta_3$ in these cases is respectively $1$ and $2$, so in general if we assume both $e \geq \theta_1$ and $e \geq 1+\inddd$ we also have $e \geq \theta_3$. As a consequence one can replace the hypothesis $e \geq \max_{1\leq i\leq3} \theta_i$ in the proposition with
\[
e \geq \max\left( \theta_1, \theta_2, 1+\inddd \right) .
\]
\end{remark}

\subsection{Conclusion of the study of the naive spaces}
\label{subsect_concl_naive}

We keep the notations $\theta_1, \theta_2$ and $\rho_1, \rho_2, \rho_3$ from proposition \ref{min_prop} ($\theta_3$ is unnecessary thanks to remark \ref{discard_th3} above), to which we add
\begin{equation} \label{def_rho0}
\rho_0 := -(q-4)e + 2(q-4)(g-1) + ((d+1)q-6)l + ((d+1)(q-2)-2)|\mu| ,
\end{equation}
which comes from the bound in lemma \ref{sing_cv}, as explained in the next proof. We also denote the maximum of all the $\rho_i$ by
\begin{equation} \label{def_rho_max}
\rho := \max_{0 \leq i \leq 3} \rho_i .
\end{equation}
Recall
\[
\xi_{e,|\mu|} = (n-d+1)e - n(g-1) - nl - (n-d+1)|\mu|
\]
from (\ref{def_xi_mu}) and
\[
\calB_\mu = J(\scrC) \times \Sym^\mu(\CmS) \times \tilde W
\]
from (\ref{def_Bmu}). It is now time to deduce the theorem for the naive spaces.

\begin{theorem} \label{th_naive}
Assume $d\geq2$, $n \geq 2^d (d-1)$ and $e \geq \max\left( \theta_1, \theta_2, 1+\inddd \right)$. In $\widehat{\mathscr{M}}_{\calB_\mu}$ we have
\[
\mathbf{L}^{-\xi_{e,|\mu|}} \left[ \mathcal{M}_e^{W,\mu} \right] = \mathbf{L}^{-l-d|\mu|} \mathfrak{S}^\mu + \varepsilon_e^\mu
\]
where $\varepsilon_e^\mu$ is an error term such that
\[
w_{\calB_\mu}(\varepsilon_e^\mu) \leq \rho .
\]
In particular,
\[
\varepsilon_e^\mu \underset{e\to\infty}{\longrightarrow} 0 .
\]
\end{theorem}
We additionally recall from proposition \ref{sing_prod_final} that $\mathbf{L}^{-l-d|\mu|} \mathfrak{S}^\mu \in \widehat{\mathscr{M}}_\C$ with
\[
w\left(\mathbf{L}^{-l-d|\mu|} \mathfrak{S}^\mu - 1\right) < 0 .
\]

\begin{proof}
First remember that to apply the circle method at all we need conditions (\ref{cond_BPF}) and (\ref{cond_CM}) to be satisfied: $e \geq \max(2g, 2g-1+l')$. This maximum is at most $\theta_3$ by definition of the latter in proposition \ref{min_prop}, which is itself bounded above by $\leq \max\left( \theta_1, 1+\inddd \right)$ by remark \ref{discard_th3}: that first hypothesis is verified. Let
\[
R_e = \sum_{\delta=\eta_e+1}^\infty \mathbf{L}^{-(n+1)\delta} \sum_{\nu\vdash\delta} S_\nu
\]
be the remainder of the singular series, whose $\calB_\mu$-weight is bounded by
\[
-(q-4)e + 2(q-4)(g-1) + ((d+1)q-4)l'
\]
following lemma \ref{sing_cv} (this is where we need the hypothesis $n \geq 2^d(d-1)$). The $\rho_0$ defined before the theorem is this bound minus $2(l+d|\mu|)$. Going back to the decomposition (\ref{part_maj_min}) in major and minor arcs and then applying the final major arcs expression (\ref{maj_concl}), we can write
\begin{align*}
\mathbf{L}^{-\xi_{e,|\mu|}} \left[ \mathcal{M}_e^{W,\mu} \right] & = \mathbf{L}^{-\xi_{e,|\mu|}} \left(\Gamma_\mathfrak{M} + \Gamma_\mathfrak{m}\right) \\
& = \mathbf{L}^{-l-d|\mu|} \left( \mathfrak{S}^\mu - R_e \right) + \mathbf{L}^{-\xi_{e,|\mu|}} \Gamma_\mathfrak{m} .
\end{align*}
Putting
\[
\varepsilon_e^\mu := -\mathbf{L}^{-l-d|\mu|} R_e + \mathbf{L}^{-\xi_{e,|\mu|}} \Gamma_\mathfrak{m}
\]
we have the announced identity. As for the weight of this error term,
\[
w_{\calB_\mu}\left(-\mathbf{L}^{-l-d|\mu|} R_e\right) = -2(l+d|\mu|) + w_{\calB_\mu}(R_e) \leq \rho_0 ,
\]
while proposition \ref{min_prop} and remark \ref{discard_th3} provide
\[
w_{\calB_\mu}\left(\mathbf{L}^{-\xi_{e,|\mu|}} \Gamma_\mathfrak{m}\right) \leq \max_{1\leq i\leq3} \rho_i ,
\]
hence the upper bound on $w_{\calB_\mu}(\varepsilon_e^\mu)$. Since the four $\rho_i$ have a negative coefficient in $e$, this weight tends to $-\infty$, in other words $\varepsilon_e^\mu$ itself converges to zero.
\end{proof}

Finally, we want to know more precisely when the weight of the error term becomes negative. For this it suffices to impose $\rho<0$. Let us look at the different $\rho_i$: going back to its definition (\ref{def_rho0}) above, $\rho_0$ is negative when
\begin{equation} \label{def_omega0}
e > \omega_0 := 2(g-1) + \frac{((d+1)q-6)l + ((d+1)(q-2)-2)|\mu|}{q-4} .
\end{equation}
Recall the definition of the other $\rho_i$ from proposition \ref{min_prop}, with the $\phi_{d,q}^i$: likewise, $\rho_1 < 0$ equates to
\begin{equation} \label{def_omega1}
e > \omega_1 := 2\frac{\phi_{d,q}^1(g) + (dq-2)l + d(q-2)|\mu|}{d(q-4)} ,
\end{equation}
$\rho_2 < 0$ becomes
\begin{equation} \label{def_omega2}
e > \omega_2 := \frac{\phi_{d,q}^2(g) + (q-6)l + (q-2d+4)|\mu|}{q-4} ,
\end{equation}
and $\rho_3 < 0$ gives
\begin{equation} \label{def_omega3}
e > \omega_3 := \frac{\phi_{d,q}^3(g) + ((d-1)q-2)l}{(d-1)(q-2)-2} + |\mu|.
\end{equation}
Recall one more time that the $\theta_i$ are defined in proposition \ref{min_prop} as well.

\begin{corollary} \label{corol_naive}
Assume $d\geq2$ and $n \geq 2^d (d-1)$: as soon as
\[
e > \max\left( \theta_1-1, \theta_2-1, \omega_0, \omega_1, \omega_2, \omega_3 \right) ,
\]
the naive space $\mathcal{M}_e^{W,\mu}$ has $\calB_\mu$-dimension 
\[
\xi_{e,|\mu|} = (n-d+1)e - n(g-1) - nl - (n-d+1)|\mu| .
\]
Furthermore under the same hypotheses, $\mathcal{M}_e^{W,\mu}$ has the same number of irreducible components of maximal dimension as $\tilde W$.
\end{corollary}
\begin{proof}
The hypotheses $e > \theta_i-1$ are precisely the hypotheses $e \geq \theta_i$ of theorem \ref{th_naive}. We omitted here $e > \inddd$: in fact by construction $\phi_{d,q}^3(g) \geq (d-1)q+4$ so
\[
\omega_3 \geq \frac{(d-1)q+4}{(d-1)(q-2)-2} > 1 \geq \inddd
\]
and $e > \inddd$ is implied by $e > \omega_3$. The theorem writes as
\[
\mathbf{L}^{-\xi_{e,|\mu|}} \left[ \mathcal{M}_e^{W,\mu} \right] = \mathbf{L}^{-l-d|\mu|} \mathfrak{S}^\mu + \varepsilon_e^\mu ,
\]
with $w_{\calB_\mu}(\varepsilon_e^\mu) < 0$ thanks to the hypotheses $e > \omega_i$, and $\left(\mathbf{L}^{-l-d|\mu|} \mathfrak{S}^\mu\right)$ of weight zero (be it seen over $\C$ or $\calB_\mu$). Therefore $\mathbf{L}^{-\xi_{e,|\mu|}} \big[ \mathcal{M}_e^{W,\mu} \big]$ has weight zero as well, which by lemma \ref{weight_var} means $\dim_{\calB_\mu} \mathcal{M}_e^{W,\mu} = \xi_{e,|\mu|}$.

As for the second part of the corollary, since
\[
\calB_\mu = J(\scrC) \times \Sym^\mu(\CmS) \times \tilde W
\]
we first deduce that the complex dimension of $\mathcal{M}_e^{W,\mu}$ is $\xi_{e,|\mu|} + g + |\mu| + \dim\tilde W$, which we denote $\Xi$ for this proof. Then in $\widehat{\mathscr{M}}_\C$ we have
\begin{equation} \label{nb_compo_max_intermed}
\mathbf{L}^{-\Xi} \left[ \mathcal{M}_e^{W,\mu} \right] = \frac{[J(\scrC)]}{\mathbf{L}^g} \cdot \frac{[\Sym^\mu(\CmS)]}{\mathbf{L}^{|\mu|}} \cdot \frac{[\tilde W]}{\mathbf{L}^{\dim \tilde W}} \cdot \mathbf{L}^{-l-d|\mu|} \mathfrak{S}^\mu + \mathbf{L}^{-(g + |\mu| + \dim\tilde W)} \varepsilon_e^\mu .
\end{equation}
We are going to  make use of lemma \ref{nb_compo_max}. To begin with, on the one hand as an abelian variety the Jacobian $J(\scrC)$ is irreducible, and on the other hand $\CmS$ is irreducible so $\Sym^\mu(\CmS)$ is as well: this implies that the first two factors in (\ref{nb_compo_max_intermed}) are of the form $1$ plus something of negative weight over $\C$. The same is true of $\mathbf{L}^{-l-d|\mu|} \mathfrak{S}^\mu$ by proposition \ref{sing_prod_final}. Let $\gamma$ be the number of irreducible components of $\tilde W$ of maximal dimension: $\frac{[\tilde W]}{\mathbf{L}^{\dim \tilde W}} - \gamma$ is of negative weight. Putting all that together, and adding the error term which has negative weight over $\C$ (because of $w_{\calB_\mu}(\varepsilon_e^\mu) < 0$ and point 1 of proposition \ref{prop_weight}), we obtain
\[
w\left( \mathbf{L}^{-\Xi} \left[ \mathcal{M}_e^{W,\mu} \right] - \gamma \right) < 0 .
\]
An application of lemma \ref{nb_compo_max} in the converse direction concludes. Note that $\tilde W$ and $W$ have the same number $\gamma$ of irreducible components of maximal dimension by the same lemma, since
\[
[\tilde W] = \mathbf{L}^l \left(1-\mathbf{L}^{-1}\right)^{\#|\scrS|} [W]
\]
(see (\ref{class_GmP})).
\end{proof}

\newpage

\section{The space of morphisms}
\label{sect_Hom}

In this section we relate the space of morphisms to naive spaces in order to show the main theorem. We keep the assumption $d\geq2$, but replace $n \geq 2^d (d-1)$ with
\begin{equation} \label{assum_n_Hom}
n \geq 2^{\max(d,3)} (d-1) ,
\end{equation}
i.e. when $d=2$ we request $n\geq8$ instead of $n\geq4$. We also keep the fixed conditions $(\scrS, W)$. In what follows we need to take $e$ satisfying the condition in theorem \ref{th_naive} for $\mu=\varnothing$, i.e. $l'=l$; thus we denote
\[
\left\{\begin{array}{l}
\theta_1^0 := 4g-\inddd+2l\\
\theta_2^0 := 2dg-2+l
\end{array}\right.
\]
and fix a
\begin{equation} \label{assum_e_1}
e \geq \max\left( \theta_1^0, \theta_2^0, 1+\inddd \right) .
\end{equation}
We go back to considering the moduli space $\Hom^e(\scrC,V|W)$ defined in the introduction, and whose construction was made explicit in remark \ref{constr_fibr_Hom}.

\subsection{Motivic inclusion-exclusion}
\label{subsec_MIE}

We are now going to use the theorem of \emph{motivic inclusion-exclusion} proved by \cite{BH21} to precisely link the space of morphisms to naive spaces and express it as the sum of the desired product and an error term; the latter will be bounded in the following subsection. It should be noted that this process of motivic inclusion exclusion first appeared, in a special case, in the section 3 of \cite{VW15}.

In remark (\ref{constr_fibr_Hom}) we gave the description
\[
\Hom^e(\scrC,V|W) = \left\{\left. (L, \bar{\mathbf{u}}) \in (\pi_*\Poinc_e)^{n+1} /\mathbb{G}_m \,\right|\, \bar{\mathbf{u}} \in \Hom^L(\scrC,V|W) \right\} .
\]
Note that as long as we see the $\calP_e$ from (\ref{def_calP}) as over $\Pic^e(\scrC)$ (i.e. take $\mu = \varnothing$), $(\pi_*\Poinc_e)^{n+1}$ is precisely $\calP_e^{\diamond(n+1)}$. Combining this with definitions (\ref{def_Hom_L}) and (\ref{def_Hom_LW}) we can then write
\begin{equation} \label{Hom_H}
\Hom^e(\scrC,V|W) = \calH_e^W / \mathbb{G}_m
\end{equation}
with
\[
\calH_e^W := \left\{ (L, \mathbf{u}) \in \calP_e^{\diamond(n+1)} \,\left|\,
\begin{array}{l}
\forall x\in\scrC, \mathbf{u}(x)\neq0 \\
f(\mathbf{u}) = 0 \\
\mathbf{u}|_\scrS \in \tilde W
\end{array} \right.\right\} .
\]
Now recall from (\ref{def_calM}) that
\[
\calM_e^{W,\mu} = \left\{ (L, \scrT, \mathbf{u}) \in \calP_e^{\diamond(n+1)} \,\left|\,
\begin{array}{l}
f(\mathbf{u}) = 0 \\
\mathbf{u}|_{\scrST} \in \tilde W \times \{0\}
\end{array} \right.\right\} ,
\]
and observe that
\[
\calH_e^W = \left\{\left. (L,\mathbf{u}) \in \calM_e^{W,\varnothing} \,\right|\, \forall x\in\scrC, \mathbf{u}(x)\neq0 \right\} .
\]
More precisely if we consider the projection
\[
\varpi :
\left\{\begin{array}{rcl}
\calM_e^{W,[1]} & \longrightarrow & \calM_e^{W,\varnothing} \\
(L, x, \mathbf{u}) & \longmapsto & (L, \mathbf{u})
\end{array}\right.
\]
we see that
\begin{equation} \label{calH_calM}
\calM_e^{W,\varnothing} = \calH_e^W \sqcup \varpi\left(\calM_e^{W,[1]}\right) .
\end{equation}
To compute the class of $\varpi\left(\calM_e^{W,[1]}\right)$ we will apply the following theorem, as announced. It features the same set of partitions $\mathcal{Q}$ from identity (\ref{inv_zeta_sum}).

\begin{theorem}[Motivic inclusion-exclusion -- \cite{BH21}, theorem 7.2.1]
If $X$ is a variety over a field of characteristic zero and $\varphi : Y \to X$ is a quasi-finite map, then in $K_0\left(\mathrm{Var}_X\right)$,
\[
[\varphi(Y)] = \sum_{\mu\in\mathcal{Q}\smallsetminus\{\varnothing\}} (-1)^{\|\mu\|-1} \left[ \Sym_{/X}^\mu(Y) \right] .
\]
Moreover for $\|\mu\|$ large enough the symmetric powers are empty: in particular the sum is finite.
\end{theorem}
More precisely this is the version over $S$ mentioned in remark 7.2.5 of the article. We are actually going to apply the theorem with
\begin{equation} \label{def_calM*}
\calM_{e,*}^{W,\mu} := \calM_e^{W,\mu} \smallsetminus (\Pic^e(\scrC) \times \Sym^\mu(\CmS) \times \{0\})
\end{equation}
(i.e. we remove $\mathbf{u} = 0$). Remark that as soon as $\scrS \neq \varnothing$, $\calM_{e,*}^{W,\mu} = \calM_e^{W,\mu}$.

\begin{proposition} \label{calH_sum_calM}
In $K_0\left(\mathrm{Var}_{\calM_{e,*}^{W,\varnothing}}\right)$ we have
\[
\left[ \calH_e^W \right] = \sum_{\mu\in\mathcal{Q}} (-1)^{\|\mu\|} \left[ \calM_{e,*}^{W,\mu} \right] .
\]
\end{proposition}

\begin{proof}
By removing $0$, (\ref{calH_calM}) becomes
\[
\calM_{e,*}^{W,\varnothing} = \calH_e^W \sqcup \varpi\left(\calM_{e,*}^{W,[1]}\right) ,
\]
which implies
\begin{equation} \label{calH_sum_calM_intermed}
\left[ \calH_e^W \right] = \left[\calM_{e,*}^{W,\varnothing}\right] - \left[\varpi\left(\calM_{e,*}^{W,[1]}\right)\right] .
\end{equation}
For every $(L,\mathbf{u}) \in \calM_{e,*}^{W,\varnothing}$ (in particular $\mathbf{u} \neq 0$) we have
\[
\varpi^{-1}(L,\mathbf{u}) \simeq \left\{ x \in \CmS \,\big|\, \mathbf{u}(x) = 0 \right\} :
\]
this is finite (and even of cardinality $\leq e$). Therefore the map
\[
\varpi : \calM_{e,*}^{W,[1]} \longrightarrow \calM_{e,*}^{W,\varnothing}
\]
is quasi-finite, and we can apply the motivic inclusion-exclusion theorem, which gives
\[
\left[\varpi\left(\calM_{e,*}^{W,[1]}\right)\right] = \sum_{\mu\in\mathcal{Q}\smallsetminus\{\varnothing\}} (-1)^{\|\mu\|-1} \left[ \Sym_{/ \calM_{e,*}^{W,\varnothing}}^\mu\left( \calM_{e,*}^{W,[1]} \right)\right] .
\]

We claim that for $\mu \neq \varnothing$, the symmetric power is precisely $\calM_{e,*}^{W,\mu}$. In fact its points are the pairwise distinct tuples $(L, x_{i,j}, \mathbf{u})_{i,j}$ with $(L,\mathbf{u}) \in \calM_{e,*}^{W,\varnothing}$ and $\mathbf{u}(x_{i,j}) = 0$, up to permutation along $\mu$. In other words such a point is the data of a $(x_{i,j})_{i,j} = \scrT \in \Sym^\mu(\CmS)$ and a $(L,\mathbf{u}) \in \calM_{e,*}^{W,\varnothing}$ such that $\mathbf{u}|_{|\scrT|} = 0$; i.e. a $\scrT \in \Sym^\mu(\CmS)$ and a $(L,\mathbf{u}) \in \left(\calM_{e,*}^{W,\mu}\right)_\scrT$. Hence the claim.

As a consequence we have
\[
\left[\varpi\left(\calM_{e,*}^{W,[1]}\right)\right] = \sum_{\mu\in\mathcal{Q}\smallsetminus\{\varnothing\}} (-1)^{\|\mu\|-1} \left[\calM_{e,*}^{W,\mu}\right] ,
\]
whose combination with (\ref{calH_sum_calM_intermed}) gives the expected identity.
\end{proof}

\begin{remark} \label{zero_mu_big}
In the fashion of (\ref{exseq_9}), let $\calP_{e(-|\mu|)}$ be the subbundle of $\calP_e$ over
\[
\Pic^e(\scrC) \times \Sym^\mu(\CmS)
\]
whose fibres are the $P_{L(-|\scrT|)} \subset P_L$. By construction, $\calM_e^{W,\mu} \subset \calP_{e(-|\mu|)}$. If $|\mu| > e$ then $\calP_{e(-|\mu|)} = 0$, which forces $\calM_{e,*}^{W,\mu} = \varnothing$: we retrieve the fact that the sum in proposition \ref{calH_sum_calM} is finite.
\end{remark}

Having made our naive spaces appear, we want to apply theorem \ref{th_naive}. The hypothesis to satisfy is $e \geq \max\left( \theta_1, \theta_2, 1+\inddd \right)$, where
\[
\left\{\begin{array}{l}
\theta_1 = 4g - \inddd + 2(l+|\mu|) \\
\theta_2 = 2dg - 2 + l+|\mu|
\end{array}\right. .
\]
Transforming this into an hypothesis on $|\mu|$, we obtain the condition $|\mu| \leq \min\left(\sigma_1, \sigma_2\right)$ with
\[
\left\{\begin{array}{l}
\sigma_1 := \frac{e+\inddd}{2} - 2g - l \\
\sigma_2 := e - 2dg + 2 - l
\end{array}\right. .
\]
Note that for $e$ large enough the minimum will always be $\sigma_1$, namely,
\[
\sigma_1 \leq \sigma_2 \Longleftrightarrow e \geq 4(d-1)g - 4 + \inddd .
\]
To keep computations manageable in the remainder of the article, we suppose this inequality is satisfied: assume
\begin{equation} \label{assum_e_2}
e \geq \tau := 4(d-1)g - 4 + \inddd .
\end{equation}
Under these circumstances, the hypothesis to apply the theorem becomes
\[
|\mu| \leq \sigma
\]
with
\begin{equation} \label{def_sigma}
\sigma := \sigma_1 = \frac{e+\inddd}{2} - 2g - l .
\end{equation}
For these $\mu$, over $\calB_\mu$ we have
\begin{align*}
\left[ \calM_{e,*}^{W,\mu} \right] & = \left[ \calM_e^{W,\mu} \right] - \indl & \text{by formula (\ref{def_calM*})} \\
& = \mathbf{L}^{\xi_{e,|\mu|}} \left(\mathbf{L}^{-l-d|\mu|} \mathfrak{S}^\mu + \varepsilon_e^\mu\right) - \indl & \text{by theorem \ref{th_naive}} \\
& = \mathbf{L}^{\xi_{e,|\mu|}} \left(\mathbf{L}^{-l} \mathfrak{S}^\varnothing + \varepsilon_e^\mu\right) - \indl & \text{by (\ref{sing_link_mu_empty}).}
\end{align*}
By extrapolation for all $\mu \in \mathcal{Q}$ we set an error term
\begin{equation} \label{def_error}
E_e^\mu := \left[ \calM_{e,*}^{W,\mu} \right] - \mathbf{L}^{\xi_{e,|\mu|}-l} \mathfrak{S}^\varnothing \in \widehat{\mathscr{M}}_{\calB_\mu} ;
\end{equation}
if $|\mu| \leq \sigma$ this gives
\begin{equation} \label{error_small_mu}
E_e^\mu = \mathbf{L}^{\xi_{e,|\mu|}} \varepsilon_e^\mu - \indl .
\end{equation}
We postpone the estimation of the errors until the next subsection: for now let us focus on deducing the following result. Note that this statement is over $W$ instead of $\tilde W$.

\begin{proposition} \label{th_main_no_estimate}
In $\widehat{\mathscr{M}}_{J(\scrC) \times W}$ we have
\[
\mathbf{L}^{-(\xi_e-1+l)} \left[ \Hom^e(\scrC,V|W) \right] = \left(1-\mathbf{L}^{-1}\right)^{\#|\scrS|-1} \cdot \prod_{x\in\CmS} \left( \left(1-\mathbf{L}^{-1}\right) \frac{[V]}{\mathbf{L}^{n-1}} \right) + \mathcal{E}_e
\]
where $\mathcal{E}_e$ is an error term such that
\[
w_{J(\scrC) \times W}\left(\mathcal{E}_e\right) \leq -2\xi_e + \sup_{\mu\in\mathcal{Q}} w_{J(\scrC) \times \tilde W}\left(E_e^\mu\right) .
\]
\end{proposition}
This is exactly the main theorem \ref{th_main_intro}, without yet the precise estimation of the error term. Once this proposition is proved, our goals will be:
\begin{itemize}
\item to check that the weight of the error term tends to $-\infty$ as $e\to\infty$, for equidistribution;
\item and to know when it becomes negative, for corollary \ref{corol_main_intro} on dimension and irreducible components.
\end{itemize}
See respectively proposition \ref{error_mininfty} and corollaries \ref{corol_main_dim} and \ref{corol_main_irred}.

\begin{proof}
By definition of $E_e^\mu$, over $J(\scrC) \times \tilde W$ we have
\[
\left[ \calM_{e,*}^{W,\mu} \right] = \mathbf{L}^{\xi_{e,|\mu|}-l} \cdot \left[\Sym^\mu(\CmS)\right] \cdot \mathfrak{S}^\varnothing + E_e^\mu .
\]
The next thing to do is of course to inject this into proposition \ref{calH_sum_calM} to get, over $J(\scrC) \times \tilde W$ again,
\begin{equation} \label{approxH_1}
\left[ \calH_e^W \right] = \sum_{\mu\in\mathcal{Q}} (-1)^{\|\mu\|} \mathbf{L}^{\xi_{e,|\mu|}-l} \cdot \left[\Sym^\mu(\CmS)\right] \cdot \mathfrak{S}^\varnothing + \sum_{\mu\in\mathcal{Q}} (-1)^{\|\mu\|} E_e^\mu ;
\end{equation}
however we need to check that these two sums converge in $\widehat{\mathscr{M}}_{J(\scrC) \times \tilde W}$. It is enough to check the first one. Recalling from (\ref{def_xi_mu}) that
\begin{equation} \label{xi_coef_mu}
\xi_{e,|\mu|} = \xi_e - (n-d+1)|\mu|
\end{equation}
with
\[
\xi_e = (n-d+1)e - n(g-1) - nl ,
\]
for $N \geq 0$ we have
\begin{align*}
& \sum_{|\mu| \leq N} (-1)^{\|\mu\|} \mathbf{L}^{\xi_{e,|\mu|}-l} \cdot \left[\Sym^\mu(\CmS)\right] \cdot \mathfrak{S}^\varnothing \\
& = \mathbf{L}^{\xi_e-l} \cdot \mathfrak{S}^\varnothing \cdot \sum_{|\mu| \leq N} (-1)^{\|\mu\|} \left[ \Sym^\mu(\CmS) \right] \mathbf{L}^{-(n-d+1)|\mu|} .
\end{align*}
By (\ref{inv_zeta_sum}), this last sum is the partial sum of the inverse zeta function
\[
Z_\CmS^{-1}(T) = \sum_{\mu\in\mathcal{Q}} (-1)^{\|\mu\|} \left[ \Sym^\mu(\CmS) \right] T^{|\mu|}
\]
evaluated in $T = \mathbf{L}^{-(n-d+1)}$. Let $r$ be the radius of convergence of this formal series: we claim that $r\leq1$. In fact the coefficient in front of $T^i$ is an alternating sum of $\left[ \Sym^\mu(\CmS) \right]$ with $|\mu| = i$: they have dimension $i$, so the weight of the total coefficient is at most $2i$. As a consequence, $r \leq \limsup_{i\to\infty} \frac{2i}{2i} = 1$. In particular $r < n-d+1$, and $Z_\CmS^{-1}(T)$ converges at $T = \mathbf{L}^{-(n-d+1)}$. Thanks to this, the sums in (\ref{approxH_1}) converge and
\[
\left[ \calH_e^W \right] = \mathbf{L}^{\xi_e-l} \cdot \mathfrak{S}^\varnothing \cdot Z_\CmS^{-1}\left(\mathbf{L}^{-(n-d+1)}\right) + \sum_{\mu\in\mathcal{Q}} (-1)^{\|\mu\|} E_e^\mu .
\]

Then, thanks to (\ref{zeta_prod}) and properties of motivic Euler products, the inverse zeta function can be expressed as $Z_\CmS^{-1}(T) = \prod_{x\in\CmS} (1-T)$. By a renormalising change of variable (which is permitted by \cite{Bil23}, paragraph 3.8.1, property 4) we deduce
\[
Z_\CmS^{-1}\left(\mathbf{L}^{-(n-d+1)}T\right) = \prod_{x\in\CmS} \left(1-\mathbf{L}^{-(n-d+1)}T\right) .
\]
Then we evaluate at $T=1$ using notation (\ref{MEP_noT}): this and proposition \ref{sing_prod_final} entail
\begin{align*}
& \mathbf{L}^{-l} \mathfrak{S}^\varnothing \cdot Z_\CmS^{-1}\left(\mathbf{L}^{-(n-d+1)}\right) \\
& = \prod_{x\in\CmS} \left( \frac{1-\mathbf{L}^{-1}}{1-\mathbf{L}^{-(n-d+1)}} \cdot \frac{[V]}{\mathbf{L}^{n-1}} \right) \cdot \prod_{x\in\CmS} \left(1-\mathbf{L}^{-(n-d+1)}\right) \\
& = \prod_{x\in\CmS} \left( \left(1-\mathbf{L}^{-1}\right) \frac{[V]}{\mathbf{L}^{n-1}} \right)
\end{align*}
and
\[
\left[ \calH_e^W \right] = \mathbf{L}^{\xi_e} \prod_{x\in\CmS} \left( \left(1-\mathbf{L}^{-1}\right) \frac{[V]}{\mathbf{L}^{n-1}} \right) + \sum_{\mu\in\mathcal{Q}} (-1)^{\|\mu\|} E_e^\mu \in \widehat{\mathscr{M}}_{J(\scrC) \times \tilde W} .
\]

The next step is to descend from $\tilde W$ to $W$. Recall that $W = \tilde W \! \big/ \mathbb{G}_{m, P_{|\scrS}}$ with $\mathbb{G}_{m, P_{|\scrS}} = \prod_{p\in|\scrS|} \mathbb{G}_{m, P_{|l_p p}}$. The invertible elements of the $\C$-algebra $P_{|l_p p} \simeq \C[t]/t^{l_p}$ are those with invertible coefficient at order $0$, so $\mathbb{G}_{m, P_{|l_p p}} \simeq \mathbb{G}_m \times \mathbb{A}^{l_p-1}$. Hence
\begin{equation} \label{class_GmP}
\left[\mathbb{G}_{m, P_{|\scrS}}\right] = \prod_{p\in|\scrS|} \left( (\mathbf{L}-1)\mathbf{L}^{l_p-1} \right) = \mathbf{L}^l \left(1-\mathbf{L}^{-1}\right)^{\#|\scrS|} :
\end{equation}
the main term in the estimation gains this factor when seen over $J(\scrC) \times W$ instead. Then, by (\ref{Hom_H}), to recover $\left[ \Hom^e(\scrC,V|W) \right]$ one just has to divide by $[\mathbb{G}_m] = \mathbf{L}\left(1-\mathbf{L}^{-1}\right)$, with eventually
\begin{align*}
\left[ \Hom^e(\scrC,V|W) \right] =\ & \mathbf{L}^{\xi_e-1+l} \left(1-\mathbf{L}^{-1}\right)^{\#|\scrS|-1} \prod_{x\in\CmS} \left( \left(1-\mathbf{L}^{-1}\right) \frac{[V]}{\mathbf{L}^{n-1}} \right) \\
& + \left( \mathbf{L}\left(1-\mathbf{L}^{-1}\right) \right)^{-1} \sum_{\mu\in\mathcal{Q}} (-1)^{\|\mu\|} E_e^\mu
\end{align*}
over $J(\scrC) \times W$. Renormalising by $\mathbf{L}^{\xi_e-1+l}$ yields the expected formula, with
\[
\mathcal{E}_e := \left( \mathbf{L}^{\xi_e+l}\left(1-\mathbf{L}^{-1}\right) \right)^{-1} \sum_{\mu\in\mathcal{Q}} (-1)^{\|\mu\|} E_e^\mu .
\]
The weight of this error term satisfies
\[
w_{J(\scrC) \times W}\left(\mathcal{E}_e\right) \leq -2(\xi_e + l) + \sup_{\mu\in\mathcal{Q}} w_{J(\scrC) \times W}\left(E_e^\mu\right)
\]
by point 1 of proposition \ref{prop_weight}; mind that on the right-hand side we have $W$ but want to go back to $\tilde W$. To change this it suffices to apply point 2 of proposition \ref{prop_weight} as follows:
\[
w_{J(\scrC) \times W}\left(E_e^\mu\right) \leq 2\dim_{J(\scrC) \times W}\left(J(\scrC) \times \tilde W\right) + w_{J(\scrC) \times \tilde W}\left(E_e^\mu\right) .
\]
The dimension is $\dim_W \tilde W = \dim \mathbb{G}_{m, P_{|\scrS}} = l$, which concludes.
\end{proof}

\subsection{Bounding the weight of the $E_e^\mu$}

The goal of this subsection is to establish an upper bound on the weight of the error term $E_e^\mu$ defined in (\ref{def_error}), which following proposition \ref{th_main_no_estimate} will provide a weight bound for the main error term $\mathcal{E}_e$. We will work over $J(\scrC) \times \tilde W$, where we have
\begin{equation} \label{E_JW}
E_e^\mu = \left[ \calM_{e,*}^{W,\mu} \right] - \mathbf{L}^{\xi_{e,|\mu|}-l} \cdot \left[\Sym^\mu(\CmS)\right] \cdot \mathfrak{S}^\varnothing .
\end{equation}
Our treatment of $E_e^\mu$ will depend on the value of $|\mu|$: as mentioned in (\ref{error_small_mu}), if $|\mu| \leq \sigma$ we apply the theorem for naive spaces. Moreover when $|\mu|>e$, remark \ref{zero_mu_big} yields $\calM_{e,*}^{W,\mu} = \varnothing$, which leaves only the second half of $E_e^\mu$. We also have to treat the case $\sigma < |\mu| \leq e$. We study them successively.

Two notes: in this subsection we will consider the $\rho_i$ from (\ref{def_rho0}) and proposition \ref{min_prop} as well as their maximum $\rho$ from (\ref{def_rho_max}) for different values of $|\mu|$, so we will write them with this number as an argument: $\rho_i(|\mu|)$ and $\rho(|\mu|)$. Furthermore the $q = \frac{n+1}{2^{d-2}(d-1)}$ defined in (\ref{def_q}) is going to appear a lot again; in this section the assumption (\ref{assum_n_Hom}) on $n$ translates as
\begin{equation} \label{assum_q_Hom}
\left\{ \begin{array}{ll}
q>4 & \text{if } d\geq3 \\
q>8 & \text{if } d=2
\end{array} \right. .
\end{equation}

\subsubsection{Case $|\mu| \leq \sigma$}

First suppose $|\mu| \leq \sigma$. In this case we saw in (\ref{error_small_mu}) that
\[
E_e^\mu = \mathbf{L}^{\xi_{e,|\mu|}} \varepsilon_e^\mu - \indl \in \widehat{\mathscr{M}}_{\calB_\mu} :
\]
by theorem \ref{th_naive} we have
\[
w_{\calB_\mu}\left(E_e^\mu\right) \leq \max\left( 2\xi_{e,|\mu|} + \rho(|\mu|) ,\ 0 \right) ,
\]
where the $0$ is relevant only in the case $l=0$. Recall that
\[
\calB_\mu = J(\scrC) \times \Sym^\mu(\CmS) \times \tilde W :
\]
we apply point 2 of proposition \ref{prop_weight} then formula (\ref{xi_coef_mu}) to get
\begin{align*}
w_{J(\scrC) \times \tilde W}\left(E_e^\mu\right) & \leq 2|\mu| + \max\left( 2\xi_{e,|\mu|} + \rho(|\mu|) ,\ 0 \right) \\
& = \max\big( 2\left(\xi_e-(n-d)|\mu|\right) + \rho(|\mu|) ,\ 2|\mu| \big) .
\end{align*}
Next we take the maximum for $0 \leq |\mu| \leq \lfloor\sigma\rfloor$: for $2|\mu|$ it is obviously $2\lfloor\sigma\rfloor$. Let us examine the coefficient in front of $|\mu|$ in $2\left(\xi_e-(n-d)|\mu|\right) + \rho_i(|\mu|)$, depending on $0\leq i\leq3$. We use that $n = 2^{d-2}(d-1)q - 1$ (\ref{def_q}):
\begin{itemize}
\item $i=0$: $-2(n-d) + (d+1)(q-2)-2 = -\big((2^{d-1}-1)(d-1) - 2\big)q - 2$;
\item $i=1$: $-2(n-d) + d(q-2) = -\big((2^{d-1}-1)(d-1) - 1\big)q + 2$;
\item $i=2$: $-2(n-d) + q-2d-4 = -\big(2^{d-1}(d-1) - 1\big)q - 2 < 0$;
\item $i=3$: $-2(n-d) + (d-1)(q-2)-2 = -(2^{d-1}-1)(d-1)q + 2 < 0$.
\end{itemize}
The last two are always negative, while the first two are negative as long as $d\geq3$ (which gives a maximum for $|\mu|=0$) but positive for $d=2$ (respectively equal to $q-2$ and $2$; maximum for $|\mu|=\lfloor\sigma\rfloor$). As a consequence, for $|\mu| \leq \sigma$, $w_{J(\scrC) \times \tilde W}\left(E_e^\mu\right)$ is bounded by
\begin{equation} \label{mu_small_d3}
\max\big( 2\xi_e + \rho(0) ,\ 2\lfloor\sigma\rfloor \big) \leq \max\big( 2\xi_e + \rho(0) ,\ 2\sigma \big)
\end{equation}
for $d\geq3$, and by
\begin{align} \label{mu_small_d2}
& \max\left( 2\xi_e + \rho_0(0) + (q-2)\lfloor\sigma\rfloor ,\ 2\xi_e + \rho_1(0) + 2\lfloor\sigma\rfloor ,\ 2\xi_e + \max_{i\in\{2,3\}} \rho_i(0) ,\ 2\lfloor\sigma\rfloor \right) \nonumber \\
& \leq \max\left( 2\xi_e + \rho_0(0) + (q-2)\sigma ,\ 2\xi_e + \rho_1(0) + 2\sigma ,\ 2\xi_e + \max_{i\in\{2,3\}} \rho_i(0) ,\ 2\sigma \right)
\end{align}
if $d=2$ (for $i\in\{0,1\}$ we extract the coefficient of $|\mu|=\lfloor\sigma\rfloor$ since we computed it already). In both cases, the last term $2\sigma$ can be forgotten if $l\neq0$.

\subsubsection{Case $\sigma < |\mu| \leq e$}

Now take $\sigma < |\mu| \leq e$. We look at the two terms of
\[
E_e^\mu = \left[ \calM_{e,*}^{W,\mu} \right] - \mathbf{L}^{\xi_{e,|\mu|}-l} \cdot \left[\Sym^\mu(\CmS)\right] \cdot \mathfrak{S}^\varnothing \in \widehat{\mathscr{M}}_{J(\scrC) \times \tilde W}
\]
independently. Recalling from proposition \ref{sing_prod_final} that $w_{J(\scrC) \times \tilde W}\left( \mathbf{L}^{-l} \mathfrak{S}^\varnothing \right) = 0$, the second one has weight
\begin{align} \label{mu_mid_B3}
w_{J(\scrC) \times \tilde W}\left( \mathbf{L}^{\xi_{e,|\mu|}-l} \cdot \left[\Sym^\mu(\CmS)\right] \cdot \mathfrak{S}^\varnothing \right) & = 2\left( \xi_{e,|\mu|} + |\mu| \right) \nonumber \\
& = 2\left(\xi_e-(n-d)|\mu|\right) \nonumber \\
& \leq 2\left(\xi_e-(n-d)(\lfloor\sigma\rfloor+1)\right) \nonumber \\
& \leq 2\left(\xi_e-(n-d)\left(\sigma+1/2\right)\right) ;
\end{align}
to get to the last line we use that $\sigma \in \frac{1}{2}\mathbb{Z}$ and $\lfloor\sigma\rfloor \geq \sigma-1/2$.

As for the first term $\left[ \calM_{e,*}^{W,\mu} \right]$, let $\mu'$ be the partition $\lfloor\sigma\rfloor \cdot [1] = [1, \dots, 1]$. For any $\scrT \in \Sym^\mu(\CmS)$, since $|\mu'| = \lfloor\sigma\rfloor < |\mu|$ there exists a $\scrT' \subset \scrT$ of partition $\mu'$, and by construction of the naive spaces (\ref{def_calM}) we have $\big(\calM_e^{W,\mu}\big)_\scrT \hookrightarrow \big(\calM_e^{W,\mu'}\big)_{\scrT'}$. As a consequence,
\begin{align*}
\dim_{\calB_\mu} \calM_{e,*}^{W,\mu} & \leq \dim_{\calB_\mu}  \calM_e^{W,\mu} \\
& = \sup_{\scrT \in \Sym^\mu(\CmS)} \dim_{J(\scrC) \times \kappa(\scrT) \times \tilde W} \left(\calM_e^{W,\mu}\right)_\scrT \\
& \leq \sup_{\scrT' \in \Sym^{\mu'}(\CmS)} \dim_{J(\scrC) \times \kappa(\scrT') \times \tilde W} \big(\calM_e^{W,\mu'}\big)_{\scrT'} \\
& = \dim_{\calB_{\mu'}} \calM_e^{W,\mu'} .
\end{align*}
From this and lemma \ref{weight_var} we deduce
\[
w_{\calB_\mu}\left(\calM_{e,*}^{W,\mu}\right) \leq w_{\calB_{\mu'}}\big(\calM_e^{W,\mu'}\big) .
\]
The theorem for naive spaces \ref{th_naive} can be applied with $|\mu'| = \lfloor\sigma\rfloor$:
\[
\left[ \mathcal{M}_e^{W,\mu'} \right] = \mathbf{L}^{\xi_{e,\lfloor\sigma\rfloor}} \left(\mathbf{L}^{-l} \mathfrak{S}^\varnothing + \varepsilon_e^{\mu'}\right) .
\]
The two previous lines, together with $w_{\calB_{\mu'}}\left(\mathbf{L}^{-l} \mathfrak{S}^\varnothing\right) = 0$, provide
\[
w_{\calB_\mu}\left(\calM_{e,*}^{W,\mu}\right) \leq 2\xi_{e,\lfloor\sigma\rfloor} + \max\big( 0 , \rho(\lfloor\sigma\rfloor) \big) .
\]
To obtain a weight bound over $J(\scrC) \times \tilde W$ it suffices to add $2|\mu|$ to the above (which doesn't depend on $|\mu|$); then taking the maximum for $\lfloor\sigma\rfloor+1 \leq |\mu| \leq e$ gives
\begin{equation}
w_{J(\scrC) \times \tilde W} \left(\calM_{e,*}^{W,\mu}\right) \leq 2\left(\xi_{e,\lfloor\sigma\rfloor}+e\right) + \max\big( 0 , \rho(\lfloor\sigma\rfloor) \big) .
\end{equation}
We would like to replace $\lfloor\sigma\rfloor$ with either $\sigma-1/2$ or $\sigma$ depending on the sign of its coefficient. For the first option in the $\max$,
\begin{align} \label{mu_mid_B4}
2\left(\xi_{e,\lfloor\sigma\rfloor}+e\right) & = 2\left(\xi_e+e - (n-d+1)\lfloor\sigma\rfloor\right) \nonumber \\
& \leq 2\left(\xi_e+e - (n-d+1)\left(\sigma-1/2\right)\right) .
\end{align}
For the second option, let $c_i$ be the coefficient in front of $\lfloor\sigma\rfloor$ in $2\left(\xi_{e,\lfloor\sigma\rfloor}+e\right) + \rho_i(\lfloor\sigma\rfloor)$, that is $-2(n-d+1)$ plus the coefficient of $|\mu|$ in $\rho_i$. We retrieve the coefficients computed for the case $|\mu|\leq\sigma$, minus $2$, that is
\begin{equation} \label{value_ci}
\left\{\begin{array}{ll}
c_0 = -\big((2^{d-1}-1)(d-1) - 2\big)q - 4 \\
c_1 = -\big((2^{d-1}-1)(d-1) - 1\big)q & \leq 0 \\
c_2 = -\big(2^{d-1}(d-1) - 1\big)q - 4 & \leq 0\\
c_3 = -(2^{d-1}-1)(d-1)q & \leq 0
\end{array}\right.
\end{equation}
with $c_0\leq0$ as soon as $d\geq3$ but $c_0 = q-4 > 0$ if $d=2$. As a consequence we have
\[
2\left(\xi_{e,\lfloor\sigma\rfloor}+e\right) + \rho(\lfloor\sigma\rfloor) = 2\left(\xi_e+e\right) + \max_{0\leq i\leq3} (\rho_i(0) + c_i\lfloor\sigma\rfloor) ,
\]
which for $d\geq3$ is at most
\[
2\left(\xi_e+e\right) + \max_{0\leq i\leq3} (\rho_i(0) + c_i(\sigma-1/2)) ,
\]
and for $d=2$ is at most
\[
2\left(\xi_e+e\right) + \max\left( \rho_0(0) + c_0\sigma ,\ \max_{1\leq i\leq3} (\rho_i(0) + c_i(\sigma-1/2)) \right) .
\]

Putting this together with (\ref{mu_mid_B3}) and (\ref{mu_mid_B4}), as long as $\sigma < |\mu| \leq e$ the weight of $E_e^\mu$ is bounded by
\begin{align} \label{mu_mid}
2\xi_e + \max\bigg( & -2(n-d)\left(\sigma+1/2\right) , \nonumber \\
& 2\left(e - (n-d+1)\left(\sigma-1/2\right)\right) , \\
& 2e + \max_{0\leq i\leq3} (\rho_i(0) + c_i(\sigma-1/2)) \bigg) \nonumber
\end{align}
for $d\geq3$, and the same thing with $c_0\sigma$ instead of $c_0(\sigma-1/2))$ when $d=2$.

\subsubsection{Case $e < |\mu|$}

As mentioned at the start of this subsection and based on remark \ref{zero_mu_big}, as soon as $|\mu| > e$ we have
\[
E_e^\mu = - \mathbf{L}^{\xi_{e,|\mu|}-l} \cdot \left[\Sym^\mu(\CmS)\right] \cdot \mathfrak{S}^\varnothing \in \widehat{\mathscr{M}}_{J(\scrC) \times \tilde W} ,
\]
so that, similarly to (\ref{mu_mid_B3}),
\begin{align} \label{mu_large}
w_{J(\scrC) \times \tilde W}\left(E_e^\mu\right) & = 2\left(\xi_e-(n-d)|\mu|\right) \nonumber \\
& \leq 2\left(\xi_e-(n-d)(e+1)\right) .
\end{align}

\subsubsection{Putting cases together}
\label{subsub_Bj}

In order to sum up the previous computations, we take the following shorthands for the bounds after removing $2\xi_e$. Let's define them for $d\geq3$ first. Coming from $|\mu|\leq\sigma$, by (\ref{mu_small_d3}) and (\ref{mu_small_d2}):
\begin{itemize}

\item $B^{(1)} := \max_{0\leq i\leq3} B_i^{(1)}$ with $B_i^{(1)} := \rho_i(0)$
\footnote{For $d\geq3$, $B^{(1)} = \rho(0)$.};

\item $B^{(2)} := 2(-\xi_e+\sigma)$ if $l=0$, and $B^{(2)} := -\infty$ otherwise.

\end{itemize}
Coming from $\sigma < |\mu| \leq e$, by (\ref{mu_mid}):
\begin{itemize}

\item $B^{(3)} := -2(n-d)\left(\sigma+1/2\right)$;

\item $B^{(4)} := 2\left(e - (n-d+1)\left(\sigma-1/2\right)\right)$;

\item $B^{(5)} := \max_{0\leq i\leq3} B_i^{(5)}$ with $B_i^{(5)} := 2e + \rho_i(0) + c_i(\sigma-1/2)$.

\end{itemize}
Finally, coming from $e < |\mu|$, by (\ref{mu_large}):
\begin{itemize}

\item $B^{(6)} := -2(n-d)(e+1)$.

\end{itemize}
If $d=2$ we take the same notations except for $B_0^{(1)} := \rho_0(0)+(q-2)\sigma$, $B_1^{(1)} := \rho_1(0)+2\sigma$ and $B_0^{(5)} := 2e + \rho_0(0) + c_0\sigma$. We proved the following bound on the error term $\mathcal{E}_e$ from proposition \ref{th_main_no_estimate}:
\begin{equation} \label{weight_Bj_1}
w_{J(\scrC) \times W}\left(\mathcal{E}_e\right) \leq \max_{1\leq j\leq 6} B^{(j)} .
\end{equation}
We can check that three of the $B^{(j)}$ are not relevant.
\begin{lemma}
The bounds $B^{(2)}$, $B^{(3)}$ and $B^{(6)}$ are smaller or equal to $B^{(4)}$.
\end{lemma}

\begin{proof}
First recall the definition (\ref{def_sigma}) of $\sigma$: one has $\sigma < e$, so
\[
B^{(6)} = -2(n-d)(e+1) < -2(n-d)\left(\sigma+1/2\right) = B^{(3)} .
\]
Next simply compute
\begin{align*}
B^{(4)} - B^{(3)} & = 2\left(e - (n-d+1)\left(\sigma-1/2\right)\right) + 2(n-d)\left(\sigma+1/2\right) \\
& = 2(e-\sigma) + 2(n-d)+1 \\
& > 0 .
\end{align*}
Then we rewrite
\begin{align} \label{expr_B4}
B^{(4)} & = 2\left(e - (n-d+1)\left(\sigma-1/2\right)\right) \nonumber \\
& = 2\left(e - (n-d+1)\left(\frac{e+\inddd-1}{2} - 2g - l\right)\right) \nonumber \\
& = -(n-d-1)e + (n-d+1)\left(4g+1-\inddd+2l\right)
\end{align}
(we will use that expression again later). For proving that $B^{(2)} \leq B^{(4)}$, we can assume $l=0$:
\begin{align*}
B^{(2)} & = 2(-\xi_e+\sigma) \\
& = 2\left(-(n-d+1)e+n(g-1) + \frac{e+\inddd}{2} - 2g \right) \\
& = -(2n-2d+1)e + 2(n-2)g - 2n + \inddd
\end{align*}
and
\[
B^{(4)} - B^{(2)} = (n-d+2)e + 2(n-2d+4)g + (n-d+1)(1-\inddd)+2n-\inddd .
\]
The coefficient of $e$ and the constant one are clearly positive. As for the one of $g$,
\[
n-2d+4 \geq 2^d(d-1)-2d+4 = \left(2^d-2\right)(d-1)+2 > 0 ,
\]
which concludes.
\end{proof}

After these computations, (\ref{weight_Bj_1}) rewrites as the following lemma.
\begin{lemma} \label{weight_Bj_2}
\[
w_{J(\scrC) \times W}\left(\mathcal{E}_e\right) \leq \max\left(B^{(1)}, B^{(4)}, B^{(5)}\right) .
\]
\end{lemma}
Recall from the start of this paragraph \ref{subsub_Bj} that $B^{(1)}$ and $B^{(5)}$ are themselves maxima of $4$ terms: we actually have a total of $9$ terms. For the estimation to be relevant we need all of them to have negative coefficient in $e$ (for the convergence), or even to be negative (for the corollary on dimension). That last part is going to translate as lower bounds on $e$.

\subsection{Lower bounds on $e$}

Before making the conditions on $e$ that arise from the $B^{(j)}$ explicit, we recall the ones we already imposed. As seen in (\ref{assum_e_1}) at the start of this section, we suppose
\[
e \geq \max\left( \theta_1^0, \theta_2^0, 1+\inddd \right)
\]
with
\[
\left\{\begin{array}{l}
\theta_1^0 = 4g-\inddd+2l\\
\theta_2^0 = 2dg-2+l
\end{array}\right. .
\]
We also assumed, in (\ref{assum_e_2}),
\[
e \geq \tau = 4(d-1)g - 4 + \inddd .
\]
We thus already have $4$ fairly simple lower bounds on $e$. For what follows recall that we write
\[
q = \frac{n+1}{2^{d-2}(d-1)} ,
\]
so that $q>4$ when $d\geq3$ and $q>8$ when $d=2$.

\subsubsection{Lower bounds common to $d\geq3$ and $d=2$}
 
Let us treat the lower bounds on $e$ coming from the $B^{(j)}$, themselves defined at the start of paragraph \ref{subsub_Bj}. Only three of them ($B_0^{(1)}$, $B_1^{(1)}$ and $B_0^{(5)}$) differ in the cases $d\geq3$ and $d=2$: we focus on the other ones first.

Start with $B_2^{(1)} = \rho_2(0)$ and $B_3^{(1)} = \rho_3(0)$. Just like in (\ref{def_omega2}) and (\ref{def_omega3}), $\rho_2(0)$ is negative if and only if
\begin{equation} \label{def_omegap2}
e > \omega'_2 := {\omega_2}_{\big||\mu|=0} = \frac{\phi_{d,q}^2(g) + (q-6)l}{q-4} ,
\end{equation}
while $\rho_3(0)$ gives the condition
\begin{equation} \label{def_omegap3}
e > \omega'_3 := {\omega_3}_{\big||\mu|=0} = \frac{\phi_{d,q}^3(g) + ((d-1)q-2)l}{(d-1)(q-2)-2} .
\end{equation}
Next, expression (\ref{expr_B4}) gave $B^{(4)} = -(n-d-1)e + (n-d+1)\left(4g+1-\inddd+2l\right)$, which is negative if and only if
\begin{equation} \label{def_vartheta}
e > \vartheta := \left(1+\frac{2}{n-d-1}\right)\left(4g+1-\inddd+2l\right) .
\end{equation}
Remark that $\vartheta = \left(1+\frac{2}{n-d-1}\right)\left(\theta_1^0+1\right)$.

Now let us look at $B^{(5)}$: for $1\leq i\leq3$,
\[
B_i^{(5)} = 2e + \rho_i(0) + c_i(\sigma-1/2)
\]
with $c_i$ given by (\ref{value_ci}). First, recalling from proposition \ref{min_prop} and (\ref{def_sigma}) that $\rho_1(0) = -\frac{d(q-4)}{2}e + \phi_{d,q}^1(g) + (dq-2)l$ and $\sigma-1/2 = \frac{e-\inddt}{2} - 2g - l$, we have
\begin{align*}
B_1^{(5)} & = \frac{4 - d(q-4) + c_1}{2}e + \phi_{d,q}^1(g) - c_1\left(2g+\frac{\inddt}{2}\right) + (dq-2-c_1)l \\
& = -\left((2^{d-2}q-2)(d-1)-4\right)e + \frac{q}{2}\psi_d^1(g) + 4 + 2nl
\end{align*}
where $\psi_d^1(g)$ is taken so that\footnote{We write it this way because the constant coefficient with respect to $q$ is $4$ for any value of $g$.}
\[
\phi_{d,q}^1(g) - c_1\left(2g+\frac{\inddt}{2}\right) = \frac{q}{2}\psi_d^1(g) + 4 ,
\]
that is
\begin{equation}
\psi_d^1(g) := \left\{\begin{array}{ll}
\Big. \inddt \left((2^{d-1}-1)(d-1)-1\right) -2 & \text{if } g=0 \\
\Big. (4+\inddt) \left((2^{d-1}-1)(d-1)-1\right) +2d & \text{if } g=1 \\
\Big. \left((2^{d+1}-1)(d-1)+2\right)g + \inddt \left((2^{d-1}-1)(d-1)-1\right) +d-3 & \text{if } g\geq2
\end{array}\right. .
\end{equation}
In the expression of $B_1^{(5)}$ above, when $d\geq3$ the coefficient of $e$ is negative; for $d=2$ it is $-(q-6)$, which is again negative with the reinforced hypothesis $q>8$ in this case. All this gives the condition
\begin{equation} \label{def_omegapp1}
e > \omega''_1 := \frac{q\psi_d^1(g)+8 + 4nl}{(2^{d-1}q-4)(d-1)-8} .
\end{equation}
Similarly, using that $\rho_2(0) := -(q-4)e + \phi_{d,q}^2(g) + (q-6)l$,
\begin{align*}
B_2^{(5)} & = \left(2-(q-4)+\frac{c_2}{2}\right)e + \phi_{d,q}^2(g) - c_2\left(2g+\frac{\inddt}{2}\right) + (-c_2+q-6)l \\
& = \underbrace{-\left(\left(2^{d-2}(d-1)+\frac{1}{2}\right)q-4\right)}_{<0} e + \frac{q}{2}\psi_d^2(g) + 8+2\cdot\inddt + 2nl
\end{align*}
with $\psi_d^2(g)$ such that
\[
\phi_{d,q}^2(g) - c_2\left(2g+\frac{\inddt}{2}\right) = \frac{q}{2}\psi_d^2(g) + 8+2\cdot\inddt ,
\]
i.e.
\begin{equation}
\psi_d^2(g) :=
\left\{\begin{array}{ll}
\Big. \inddt\left(2^{d-1}(d-1)-1\right)-4 & \text{if } g=0 \\
\Big. (4+\inddt)\left(2^{d-1}(d-1)-1\right)+2(d-1) & \text{if } g=1 \\
\Big. (2^{d+1}+3)(d-1)g + \inddt\left(2^{d-1}(d-1)-1\right)+d-5 & \text{if } g\geq2
\end{array}\right. .
\end{equation}
We get the lower bound
\begin{equation} \label{def_omegapp2}
e > \omega''_2 := \frac{q\psi_d^2(g)+16+4\cdot\inddt + 4nl}{\left(2^{d-1}(d-1)+1\right)q-8} .
\end{equation}
Finally, because $\rho_3(0) := -((d-1)(q-2)-2)e + \phi_{d,q}^3(g) + ((d-1)q-2)l$,
\begin{align*}
B_3^{(5)} & = \left(2-((d-1)(q-2)-2)+\frac{c_3}{2}\right)e + \phi_{d,q}^3(g) - c_3\left(2g+\frac{\inddt}{2}\right) + (-c_3+(d-1)q-2)l \\
& = \underbrace{-\left(\left(\left(2^{d-2}+\frac{1}{2}\right)q-2\right)(d-1)-4\right)}_{<0} e + \frac{q}{2}\psi_d^3(g) + 4 + 2nl
\end{align*}
with
\begin{equation}
\psi_d^3(g) := (d-1) \cdot
\left\{\begin{array}{ll}
\Big. \inddt(2^{d-1}-1)+2 & \text{if } g=0 \\
\Big. (4+\inddt)(2^{d-1}-1)+4 & \text{if } g=1 \\
\Big. (2^{d+1}+3)g + \inddt(2^{d-1}-1)-3 & \text{if } g\geq2
\end{array}\right. ,
\end{equation}
which gives rise to
\begin{equation} \label{def_omegapp3}
e > \omega''_3 := \frac{q\psi_d^3(g)+8 + 4nl}{\left((2^{d-1}+1)q-4\right)(d-1)-8} .
\end{equation}

\subsubsection{Lower bounds that depend on the value of $d$}

Remain the conditions coming from $B_0^{(1)}$, $B_1^{(1)}$ and $B_0^{(5)}$. As long as $d\geq3$,
\begin{align*}
B_0^{(1)} & = \rho_0(0) \\
& = -(q-4)e + 2(q-4)(g-1) + ((d+1)q-6)l
\end{align*}
is negative for
\begin{equation} \label{def_omegap0_d3}
e > \omega'_0 \stackrel{d\geq3}{:=} {\omega_0}_{\big||\mu|=0} = 2(g-1) + \frac{(d+1)q-6}{q-4}l
\end{equation}
with $\omega_0$ from (\ref{def_omega0}). However for $d=2$,
\begin{align*}
B_0^{(1)} & = \rho_0(0) + (q-2)\sigma \\
& = -\left(\frac{q}{2}-3\right)e - 4g-\frac{3q}{2}+7 + 2(q-2)l .
\end{align*}
Since $q>8$ in this case the coefficient of $e$ is negative and we obtain the condition
\begin{equation} \label{def_omegap0_d2}
e > \omega'_0 \stackrel{d=2}{:=} \frac{-8g-3q+14 + 4(q-2)l}{q-6} .
\end{equation}
Let us do the same with $B_1^{(1)}$: if $d\geq3$,
\begin{align*}
B_1^{(1)} & = \rho_1(0) \\
& = -\frac{d(q-4)}{2}e + \phi_{d,q}^1(g) + (dq-2)l
\end{align*}
gives us the lower bound
\begin{equation} \label{def_omegap1_d3}
e > \omega'_1 \stackrel{d\geq3}{:=} {\omega_1}_{\big||\mu|=0} = 2\frac{\phi_{d,q}^1(g) + (dq-2)l}{d(q-4)}
\end{equation}
with $\omega_1$ from (\ref{def_omega1}). For $d=2$,
\begin{align*}
B_1^{(1)} & = \rho_1(0) + 2\sigma \\
& = -(q-5)e + \phi_{2,q}^1(g)-4g+1 + 2(q-2)l
\end{align*}
gets us
\begin{equation} \label{def_omegap1_d2}
e > \omega'_1 \stackrel{d=2}{:=} \frac{\phi_{2,q}^1(g)-4g+1 + 2(q-2)l}{q-5} .
\end{equation}

At last, using the definitions from the start of paragraph \ref{subsub_Bj}, remark that we can write $B_0^{(5)} = 2e + \rho_0(0) + c_0\left(\sigma-\frac{\inddt}{2}\right)$ for any $d$. But
\[
\sigma-\frac{\inddt}{2} = \frac{e+\inddd-\inddt}{2} - 2g - l = \frac{e+1}{2} - \inddt - 2g - l ,
\]
hence
\begin{align*}
B_0^{(5)} & = \left(2-(q-4)+\frac{c_0}{2}\right)e + 2(q-4)(g-1) - c_0\left(2g+\inddt-\frac{1}{2}\right)+ ((d+1)q-6-c_0)l \\
& = -\left(\left(2^{d-2}-\frac{1}{2}\right)(d-1)q-4\right)e + \frac{q}{2}\psi_d^0(g) + 6+4\cdot\inddt + 2nl
\end{align*}
where we put
\begin{equation}
\psi_d^0(g) := 4\left((2^{d-1}-1)(d-1)-1\right)g + (-1+2\cdot\inddt)\left((2^{d-1}-1)(d-1)-2\right) - 4 .
\end{equation}
The coefficient in front of $e$ is negative when $d\geq3$, and equal to $-\left(\frac{q}{2}-4\right)$ when $d=2$, which is negative as well since in this case we assume $q>8$ (this is actually the reason we need this hypothesis for $d=2$). We obtain our last condition on $e$,
\begin{equation} \label{def_omegapp0}
e > \omega''_0 := \frac{q\psi_d^0(g)+12+8\cdot\inddt + 4nl}{(2^{d-1}-1)(d-1)q-8} .
\end{equation}

\subsubsection{Redundant lower bounds and conclusion}

The first consequence of the above is that the coefficients of $e$ in $B^{(1)}, B^{(4)}, B^{(5)}$ are all negative under our hypotheses, which together with lemma \ref{weight_Bj_2} implies the following.
\begin{proposition} \label{error_mininfty}
If $d\geq3$ and $n \geq 2^d (d-1)$, or $d=2$ and $n\geq8$, then the error term $\mathcal{E}_e$ from proposition \ref{th_main_no_estimate} converges to $0$ as $e\to\infty$.
\end{proposition}
Combining proposition \ref{th_main_no_estimate} with this one proves the equidistribution theorem \ref{th_main_intro}. Furthermore the computations of the previous pages show that for the error term $\mathcal{E}_e$ to effectively be of negative weight, it is enough to have $e$ (strictly) larger than all of the following:
\begin{itemize}
\item $\theta_1^0-1$, $\theta_2^0-1$, $\inddd$, $\tau-1$;
\item $\omega'_i$ for $0\leq i\leq3$;
\item $\vartheta$;
\item $\omega''_i$ for $0\leq i\leq3$.
\end{itemize}
Some of them are always smaller than others and can be removed from this list. As a start, because of the definition (\ref{def_omegap3}) of $\omega'_3$, we have $\omega'_3 > \inddd$ like in the proof of corollary \ref{corol_naive}. Also, as we already remarked when defining $\vartheta$ in (\ref{def_vartheta}), $\vartheta = \left(1+\frac{2}{n-d-1}\right)\left(\theta_1^0+1\right) > \theta_1^0-1$.

Using asymptotics (in $g$ and $l$) as well as computations in Sage and Wolfram Alpha, one checks that for $d\geq3$, apart from the two last ones, none of the lower bounds is systematically smaller than other ones, which means all 11 of them are necessary.
\begin{remark} \label{redund_d3}
For $d\geq3$, actually as soon as $g\geq1$ one has $\omega''_3 \leq \omega''_1$.
\end{remark}

When $d=2$ however, it is clear that
\[
\tau-1 \leq \theta_2^0-1 \leq \theta_1^0-1 \leq \vartheta ,
\]
which allows us to remove $\tau-1$ and $\theta_2^0-1$ from the list. Furthermore thanks to computations with the aforementioned softwares, one always has $\omega'_2 \leq \omega''_1$ and $\omega''_2 \leq \vartheta$ : $\omega'_2$ and $\omega''_2$ can be omitted too. We are left with the lower bounds $\omega'_0, \omega'_1, \omega'_3, \vartheta, \omega''_0, \omega''_1, \omega''_3$.
\begin{remark} \label{redund_d2}
With $d=2$ still, when $g=1$ we also find $\omega'_1 \leq \omega''_1$, $\omega'_3 \leq \omega''_1$ and $\omega''_3 \leq \vartheta$, which leaves $\omega'_0, \vartheta, \omega''_0, \omega''_1$. When $g\geq2$, $\omega'_1, \omega'_3, \vartheta, \omega''_3$ are all smaller than $\omega''_1$, and the only remaining lower bounds are $\omega'_0, \omega''_0, \omega''_1$. In both cases no other systematic inequality holds.
\end{remark}

With this work done, we can finally state the precise version of the dimensional part of corollary \ref{corol_main_intro}. The notations are to be found in all of the current subsection.
\begin{corollary} \label{corol_main_dim}
Assume either
\begin{itemize}
\item $d\geq3$, $n \geq 2^d (d-1)$ and
\[
e > \max\left( \theta_2^0-1, \tau-1, \omega'_0, \omega'_1, \omega'_2, \omega'_3, \vartheta, \omega''_0, \omega''_1, \omega''_2, \omega''_3\right) ;
\]
\item or $d=2$, $n\geq8$ and
\[
e > \max\left( \omega'_0, \omega'_1, \omega'_3, \vartheta, \omega''_0, \omega''_1, \omega''_3\right) .
\]
\end{itemize}
Then the moduli space $\Hom^e(\scrC,V|W)$ has $(J(\scrC) \times W)$-dimension 
\[
\xi_e-1+l = (n-d+1)e - n(g-1) - 1 - (n-1)l .
\]
\end{corollary}
Obviously, the lower bounds demanded for $e$ can be simplified depending on the value of $g$, following remarks \ref{redund_d3} and \ref{redund_d2} above. Above $\C$, the dimension becomes
\[
(n-d+1)e - (n-1)(g-1) - ((n-1)l - \dim W) .
\]
To finish, let us state again the part of corollary \ref{corol_main_intro} concerned with irreductibilty, and prove it.
\begin{corollary} \label{corol_main_irred}
Let $(\mathcal{W}_i)_i$ be a decomposition of $W$ into irreducible components. Under the same assumptions as the previous corollary, $\big(\Hom^e(\scrC,V|\mathcal{W}_i)\big)_i$ is a decomposition of $\Hom^e(\scrC,V|W)$ into irreducible components.
\end{corollary}
Note that the $\mathcal{W}_i$ are subvarieties of $\Hom(\scrS,V)$ for the same $\scrS$ as $W$ (and therefore the same $l = \deg\scrS$). In particular the components of $W$ of dimension $\delta$ correspond to the components of $\Hom^e(\scrC,V|W)$ of dimension
\[
(n-d+1)e - (n-1)(g-1) - ((n-1)l - \delta) .
\]

\begin{proof}
We start with the case $W = \{\phi\}$, i.e. we prove that $\Hom^e(\scrC,V|\{\phi\})$ is irreducible. After the previous corollary, since $\dim J(\scrC) = g$, the complex dimension of $\Hom^e(\scrC,V|\{\phi\})$ is
\begin{equation} \label{dim_Hom_phi}
(n-d+1)e - (n-1)(g-1) - (n-1)l .
\end{equation}
As for corollary \ref{corol_naive}, with the help of lemma \ref{nb_compo_max} it can be deduced from proposition \ref{th_main_no_estimate} and corollary \ref{corol_main_dim} that $\Hom^e(\scrC,V|\{\phi\})$ has a unique irreducible component of maximal dimension, that is of dimension (\ref{dim_Hom_phi}). But the second point of theorem II.1.7 from \cite{Kol96}, which is an avatar of (\ref{expect_dim}), states that every irreducible component of this moduli space has local dimension around $\varphi$ at least $\chi(\scrC, \varphi^* T_V \otimes I_\scrS)$, where $\chi$ is the Euler characteristic, $T_V$ the tangent bundle of $V$, and $I_\scrS$ the sheaf of ideals associated to $\scrS$. As in the proof of theorem II.1.2 of the same book, this is equal to
\[
\deg(\varphi^* T_V \otimes I_\scrS) + \dim V \cdot \chi(\mathcal{O}_\scrC) = \deg(\varphi^* T_V \otimes I_\scrS) - (n-1)(g-1).
\]
Because $\varphi^* T_V$ is a vector bundle of rank $\dim V = n-1$,
\begin{align*}
\deg(\varphi^* T_V \otimes I_\scrS) & = \deg(\varphi^* T_V) + (n-1) \deg I_\scrS \\
& = (n-d+1)e - (n-1)l ,
\end{align*}
and the lower bound on the local dimension is precisely (\ref{dim_Hom_phi}). It does not depend on $\varphi$ so this is a lower bound on the dimension of the irreducible components. Since we said at the start of this proof that in our case this was an upper bound and that only one irreducible component had this dimension, the irreductibility of $\Hom^e(\scrC,V|\{\phi\})$ follows.

For the general case, consider the map
\[
\cdot|_\scrS : \Hom^e(\scrC,V|W) \longrightarrow W :
\]
it is open since it appears when taking quotients of the map $\cdot|_\scrS : P_L^{n+1} \to P_{|\scrS}^{n+1}$, which can be seen as a projection and is therefore itself open. Its fibres are the $\Hom^e(\scrC,V|\{\phi\})$, which are all irreducible by the first step of this proof. Tag 037A of \cite{SP} then implies that $\cdot|_\scrS$ induces a bijection between the irreducible components of $\Hom^e(\scrC,V|W)$ and the $\mathcal{W}_i$; in particular the former are the $\Hom^e(\scrC,V|\mathcal{W}_i)$.
\end{proof}

\subsection{Corollary about the Hodge-Deligne polynomial}
\label{subsec_HD}

In this last subsection we simply justify the corollary \ref{corol_HD_intro} concerning the Hodge-Deligne polynomial of $\Hom^e(\scrC,V|W)$. Beforehand let us state a general lemma that is a consequence of theorem 1.4.1 and lemma 1.4.2 of \cite{BL04}.

\begin{lemma}
Let $A$ be a complex abelian variety of dimension $g$. Its Hodge numbers are
\[
h^{p,q}(A) = \binom{g}{p}\binom{g}{q} ,
\]
and its Hodge-Deligne polynomial is
\[
\HD(A) = (u+1)^g (v+1)^g .
\]
\end{lemma}
The second formula follows from the first one via $\sum_{p=0}^g \binom{g}{p} u^p = (u+1)^g$.

\begin{proof}[Proof of corollary \ref{corol_HD_intro}]
Over $\C$, the formula from proposition \ref{th_main_no_estimate} becomes
\[
\left[ \Hom^e(\scrC,V|W) \right] = \mathbf{L}^{\xi_e-1+l} \left(1-\mathbf{L}^{-1}\right)^{\#|\scrS|-1} [J(\scrC)] \cdot [W] \cdot \mathscr{P} + \mathbf{L}^{\xi_e-1+l} \mathcal{E}_e
\]
with
\[
\mathscr{P} = \prod_{x\in\CmS} \left( \left(1-\mathbf{L}^{-1}\right) \frac{[V]}{\mathbf{L}^{n-1}} \right) .
\]
The lemma above provides $\HD(J(\scrC)) = (u+1)^g (v+1)^g$, so applying the motivic measure $\HD$ to the main term we obtain the polynomial in corollary \ref{corol_HD_intro}. The error term $\mathbf{L}^{\xi_e-1+l} \mathcal{E}_e$ has complex weight
\[
2(\xi_e-1+l) + 2(g+\dim W) + w_{J(\scrC) \times W}\left(\mathcal{E}_e\right) = \dim \Hom^e(\scrC,V|W) + w_{J(\scrC) \times W}\left(\mathcal{E}_e\right) ,
\]
and the terms in its Hodge-Deligne polynomial are for $p+q$ at most this weight, hence the claim regarding which part of $\HD(\Hom^e(\scrC,V|W))$ is common with the main term.
\end{proof}

Following remark \ref{MEP_noT_sum},
\begin{equation}
\mathscr{P} = \sum_{j\geq0} \Conf_{\CmS/\C}^{j[1]}\left( \left(\left(1-\mathbf{L}^{-1}\right) \frac{[V]}{\mathbf{L}^{n-1}} - 1\right)[\CmS] \right) :
\end{equation}
in particular the lacking information in corollary \ref{corol_HD_intro} is that of the Hodge-Deligne polynomials of these configuration spaces.

\newpage
\nocite{*}
\bibliographystyle{alpha}
{\small \bibliography{biblio}}

\end{document}